\documentclass[a4paper]{amsart}
\usepackage{graphics}
\usepackage{color}
\usepackage{amssymb}
\usepackage[all]{xy}
\usepackage{amsmath}
\usepackage{stmaryrd}
\CompileMatrices
\newtheorem*{introthm}{Theorem}

\newtheorem{theorem}{Theorem}[section]
\newtheorem{lemma}[theorem]{Lemma}
\newtheorem{proposition}[theorem]{Proposition}
\newtheorem{corollary}[theorem]{Corollary}
\theoremstyle{definition}
\newtheorem{definition}[theorem]{Definition}
\newtheorem{example}[theorem]{Example}

\newtheorem{construction}[theorem]{Construction}

\newtheorem{remark}[theorem]{Remark}
\theoremstyle{remark}

\numberwithin{equation}{section}
\def\Chi{{\mathbb X}}

\def\rel{{\rm rlv}}
\def\div{{\rm div}}

\def\quot{/\!\!/}
\def\mal{\! \cdot \!}

\def\rq#1{\widehat{#1}}
\def\t#1{\widetilde{#1}}
\def\b#1{\overline{#1}}
\def\bangle#1{\langle #1 \rangle}

\def\KK{{\mathbb K}}
\def\TT{{\mathbb T}}
\def\ZZ{{\mathbb Z}}

\def\QQ{{\mathbb Q}}
\def\PP{{\mathbb P}}

\def\WDiv{\operatorname{WDiv}}
\def\PDiv{\operatorname{PDiv}}
\def\id{{\rm id}}
\def\Eff{{\rm Eff}}
\def\Mov{{\rm Mov}}

\def\Ample{{\rm Ample}}
\def\SAmple{{\rm SAmple}}

\def\orb{{\rm orbit}}

\def\cov{{\rm cov}}

\def\Cl{\operatorname{Cl}}

\def\Pic{\operatorname{Pic}}

\def\Hom{{\rm Hom}}

\def\grad{{\rm grad}}

\def\Supp{{\rm Supp}}
\def\Spec{{\rm Spec}}

\def\Sing{{\rm Sing}}
\def\Star{{\rm star}}

\def\cone{{\rm cone}}
\def\lin{{\rm lin}}

\def\rank{\operatorname{rank}}

\def\topto#1{\stackrel{{\scriptscriptstyle #1}}{\longrightarrow}}
\def\faces{{\rm faces}}

\def\cov{{\rm cov}}

\newcounter{itemnumber}

\begin{document}
\title[Cox rings and combinatorics II]%
{Cox rings and combinatorics II}
\author[J.~Hausen]{J\"urgen Hausen} 
\address{Mathematisches Institut, Universit\"at T\"ubingen,
Auf der Morgenstelle 10, 72076 T\"ubingen, Germany}
\email{hausen@mail.mathematik.uni-tuebingen.de}
\begin{abstract}
We study varieties with a finitely generated 
Cox ring. 
In a first part, we generalize a combinatorial 
approach developed in earlier work for varieties
with a torsion free divisor class group to 
the case of torsion.
Then we turn to modifications, e.g., blow ups, 
and the question how the Cox ring changes 
under such maps. 
We answer this question for a certain class of 
modifications induced from modifications of
ambient toric varieties.
Moreover, we show that every variety with 
finitely generated Cox ring 
can be explicitly constructed in a finite series 
of toric ambient modifications from a 
combinatorially minimal one. 
\end{abstract}

\dedicatory{Dedicated to Ernest Borisovich Vinberg
on the occasion of his 70th birthday}

\maketitle

\section{Introduction}
\label{sec:intro}

In~\cite{Cox}, D.~Cox associated to any non-degenerate,
e.g. complete, toric variety $Z$ a multigraded homogeneous 
coordinate ring $\mathcal{R}(Z)$. His construction, 
meanwhile a standard tool in toric geometry,
allows generalization to certain 
non-toric varieties: for any normal variety $X$ having 
only constant invertible global functions and 
a finitely generated divisor class group $\Cl(X)$,
one can define a $\Cl(X)$-graded Cox ring as
\begin{eqnarray*}
\mathcal{R}(X) 
& := &
\bigoplus_{\Cl(X)} \Gamma(X, \mathcal{O}(D)),
\end{eqnarray*}
see Section~\ref{sec:crembvar} for more details.  
We are interested in the case of a finitely generated 
Cox ring.
Then one can define the total coordinate space 
$\b{X} := \Spec \, \mathcal{R}(X)$; it comes with an 
action of the diagonalizable group 
$H := \Spec \, \KK[\Cl(X)]$.
Moreover, $X$ turns out to be the good quotient 
of an open $H$-invariant subset 
$\rq{X} \subseteq \b{X}$ by the action of $H$.
The quotient map $p \colon \rq{X} \to X$ generalizes 
the concept of a universal torsor. 
Cox rings und universal torsors occur in different
settings, for some recent work, see~\cite{Br} 
and~\cite{Der}.

In~\cite{BeHa1}, a combinatorial approach 
to varieties with torsion free divisor class 
group and finitely generated Cox ring was 
presented.
The main observation is that in many cases, 
e.g. for projective $X$, 
the variety $X$ is characterized by its total 
coordinate space $\b{X}$ and combinatorial 
data living in the grading group  $\Cl(X)$.
This leads to the concept of a ``bunched ring'' 
as defining data for $X$. The task then is, 
as in toric geometry, to read off geometric 
properties of the variety $X$ from its defining
combinatorial data.
A first part of the present article
generalizes the approach of~\cite{BeHa1}
to the case of divisor class groups with torsion
and obtains similar descriptions
of the Picard group, the effective, moving,
semiample and ample cones, singularities and the 
canonical divisor.

In a second part, we consider proper modifications 
$X_1 \to X_0$, e.g. blow ups, of varieties $X_0$
with finitely generated total coordinate ring.
In general, it is a delicate task to 
describe the Cox ring of $X_1$ explicitly 
in terms of that of $X_0$; note that even finite 
generation may be lost.
Our aim is to provide some systematic insight 
by investigating modifications $X_1 \to X_0$ 
induced from stellar subdivisions of 
toric ambient varieties.
We figure out a large class of such modifications 
preserving finite generation of the Cox ring 
and show how to compute the Cox ring of $X_1$ 
from that of $X_0$ in these cases.

Let us make our approach a little more precise.
A basic observation is that many varieties $X_0$
with a finitely generated Cox ring $\mathcal{R}(X_0)$, 
e.g. all projective ones, admit a {\em neat\/} 
embedding into a toric variety $Z_0$ such that 
there is a pullback isomorphism $\Cl(Z_0) \to \Cl(X_0)$.
Then the toric universal torsor $\rq{Z}_0 \to Z_0$
admits a simple description in terms of fans.
Moreover, the inverse image $\rq{X}_0$ of $X_0$ 
under $\rq{Z}_0 \to Z_0$ is a reasonable candidate 
for a universal torsor over $X_0$, and, similarly, 
the closure $\b{X}_0$ of $\rq{X}_0$ in 
the affine closure $\b{Z}_0$ of $\rq{Z}_0$
is a candidate for the total coordinate space of 
$X_0$, see Corollary~\ref{cor:qfact2tosor}
for the accurate statements.

If $Z_1 \to Z_0$ is the toric modification 
arising from a stellar subdivision of the 
defining fan of $Z_0$, then one may define a 
strict transform $X_1 \subseteq Z_1$ and ask 
for its universal torsor and its Cox ring. 
The idea is to consider, as before,
the toric universal torsor $\rq{Z}_1 \to Z_1$, 
the inverse image $\rq{X}_1 \subseteq \rq{Z}_1$ 
of $X_1$ and the closure $\b{X}_1 \subseteq \b{Z}_1$.
These data are linked to the corresponding data 
for $X_0$ via a commutative diagram
$$
\xymatrix{
&
{\b{Y}_1}  \ar[dr]^{\quot \KK^*} \ar[dl]_{/ C}
& 
\\
{\b{X}_1}
&
& 
{\b{X}_0}
\\
{\rq{X}_1} \ar[d]_{/ H_1} \ar[u]
&
& 
{\rq{X}_0} \ar[u] \ar[d]^{/ H_0}
\\
X_1 
\ar[rr]
& &
X_0
}
$$
where $\b{Y}_1 \subseteq \b{Z}_1$ is a variety coming
with actions of a finite cyclic group $C$ and the 
multiplicative group $\KK^*$; these data can be 
explicitly computed in terms of the stellar subdivision 
and $\b{X}_0$, see Lemma~\ref{liftmapprops}.
Now a couple of technical difficulties show up.
Firstly, one has to ensure that $X_1$ is again 
neatly embedded; this is done by restricting to
{\em neat\/} ambient modifcations in the 
sense of~\ref{def:ambmod}.
Next, to ensure finite generation of the Cox ring 
$\mathcal{R}(X_1)$, we require the neat ambient 
modifcation $Z_1 \to Z_0$ to be {\em controlled\/} in
the sense that $H_1$ permutes transitively 
the components of the minus cell 
$\b{Y}^-_1 \subseteq \b{Y}_1$. 
In Theorem~\ref{controlled2small}, we then obtain 
the following.

\begin{introthm}
Let $Z_1,Z_0$ be $\QQ$-factorial toric varieties and
$\pi \colon Z_1 \to Z_0$ 
a neat controlled ambient modification
for $X_1 \subseteq Z_1$ and $X_0 \subseteq Z_0$.
Suppose that
$X_0 \subseteq Z_0$ is neatly embedded 
and $\b{X}_0 \setminus \rq{X}_0$ 
is of codimension at least two in $\b{X}_0$.
\begin{enumerate}
\item
The complement $\b{X}_1  \setminus \rq{X}_1$
is of codimension at least two in $\b{X}_1$.
If $\rq{X}_1$ is normal, then 
$\rq{X}_1 \to X_1$ is a universal torsor of $X_1$.
\item
The total coordinate space of $X_1$ is 
the normalization of $\b{X}_1$ together
with the induced action of the torus 
$H_1$.
\end{enumerate}
\end{introthm}

A similar statement holds for contractions, 
see Theorem~\ref{controlled2small2}.
The above result shows in particular, that the Cox ring 
stays finitely generated under neat controlled 
ambient modifications. 
Moreover, starting with an embedding such that $\b{X}_0$ 
is the total coordinate space of $X_0$, 
see Corollary~\ref{veryneatemb}, it gives 
a way to obtain defining equations for the total 
coordinate space $\b{X}_1$ of $X_1$ 
from those of $\b{X}_0$ 
and the data of the stellar subdivision.
In particular, for Cox rings with a single 
defining relation, Proposition~\ref{singlerel} 
gives a very explicit statement.
As Cox rings without torsion in the grading group
are always factorial, this leads 
as well to a new construction of multigraded
UFDs out of given ones, compare e.g.~\cite{SchS}
for other constructions.

As an application of the technique of 
ambient modifications, we show how to 
reduce $\QQ$-factorial projective varieties
$X$ with finitely generated Cox ring,
to ``combinatorially minimal'' 
ones, i.e., varieties $X$ that do not 
admit a class $[D] \in \Cl(X)$ such that 
$[D]$ generates an extremal ray 
of the effective cone of $X$ and, 
for some representative $D$, 
all vector spaces 
$\Gamma(X,\mathcal{O}(nD))$,
where $n > 0$ are of dimension one.
It turns out that a variety is combinatorially 
minimal if and only if its moving cone 
is the whole effective cone.
Our result is the following, see
Theorem~\ref{combcontr}.

\begin{introthm}
Every $\QQ$-factorial projective variety 
$X$ with finitely generated Cox ring arises 
from a combinatorially minimal one $X_0$
via a finite sequence
$$ 
\xymatrix{
X = X_n' \ar@{-->}[r]
&
X_{n}  \ar[r]
&
X'_{n-1}  \ar@{-->}[r]
&
X_{n-1}  \ar[r]
& 
\quad \ldots \quad
\ar[r]
& 
X_0'  
=
X_0
}
$$
where $\xymatrix@1@!{X_i' \ar@{-->}[r] & X_{i}}$
is a small birational transformation and 
$\xymatrix@1@!{X_i \ar[r] & X_{i-1}'}$
comes from a neat controlled ambient modification
of $\QQ$-factorial projective toric varieties.
\end{introthm}

The effect of each reduction step on 
the Cox ring can be explicitly calculated.
Moreover, this theorem shows that the 
class of modifications arising from
neat controlled ambient modifications is 
reasonably large.
As an example for contraction, we treat a singular 
del Pezzo surface, and show that it arises 
via neat controlled ambient modification 
from the projective plane, 
see Example~\ref{ex:delpezzo2}.
Going in the other direction, we show that 
the singularity of this surface can be 
resolved by means of neat ambient modifications 
and thereby obtain the Cox ring of its 
minimal resolution, 
see Example~\ref{ex:delpezzo3}.

\tableofcontents

\section{The Cox ring of an embedded variety}
\label{sec:crembvar}

Here, we introduce the Cox ring of a normal 
variety $X$ with finitely generated divisor class 
group and show how to compute the Cox ring 
if $X$ is neatly embedded into a toric variety.
In contrast to ~\cite{BeHa1}, we define the 
Cox ring also for the case of torsion in the 
divisor class group. For this, we closely 
follow the lines of~\cite[Section~3]{BeHa3}, 
where an analogous 
construction in terms of lines bundles 
instead of divisors is performed.

In the whole paper, we work in the 
category of (reduced) algebraic 
varieties over an algebraically closed 
field $\KK$ of characteristic zero.
By a point, we always mean a 
closed point.

Consider a normal variety  $X$  
with $\Gamma(X,\mathcal{O}^*) = {\KK^*}$
and finitely generated 
divisor class group $\Cl(X)$.  
Let $\mathfrak{D} \subseteq \WDiv(X)$ 
be a finitely generated subgroup of the group of 
Weil divisors mapping onto 
$\Cl(X)$ and consider the sheaf of 
$\mathfrak{D}$-graded algebras
$$
\mathcal{S}
\ :=  \
\bigoplus_{D \in \mathfrak{D}} \mathcal{S}_D,
\qquad\qquad
\mathcal{S}_D
\ :=  \
\mathcal{O}_X(D),
$$
where multiplication is defined via  
multiplying homogeneous sections as 
rational functions on $X$. 
Let $\mathfrak{D}^0 \subseteq \mathfrak{D}$
be the kernel of $\mathfrak{D} \to \Cl(X)$.
We fix a  {\em shifting family\/}, i.e.,
a family of $\mathcal{O}_X$-module 
isomorphisms 
$\varrho_{D^0} \colon \mathcal{S} \to \mathcal{S}$,
where $D^0 \in \mathfrak{D}^0$, such that 
\begin{itemize}
\item
$\varrho_{D^0}(\mathcal{S}_D) = \mathcal{S}_{D+D^0}$
for all $D \in \mathfrak{D}$, $D^0 \in \mathfrak{D}^0$,
\item
$\varrho_{D^0_1+D^0_2} = \varrho_{D^0_2} \circ \varrho_{D^0_1}$
for all $D^0_1,D^0_2 \in \mathfrak{D}^0$,
\item
$\varrho_{D^0}(fg) = f \varrho_{D^0}(g)$ for all $D^0   \in \mathfrak{D}^0$
and any two homogeneous $f,g$. 
\end{itemize}
To obtain such a family, 
take a basis $D_1, \ldots, D_r$ of $\mathfrak{D}$ 
such that $\mathfrak{D}^0$ is spanned by $D^0_i := a_iD_i$,
where $1 \le i \le s$ with $s \le r$, choose isomorphisms 
$\varrho_{D^0_i} \colon \mathcal{S}_0 \to \mathcal{S}_{D^0_i}$,
and, for $D^0 = b_1D^0_1 + \cdots + b_s D^0_s$, define
$\varrho_{D^0}$ on a homogeneous $f$ as   
\begin{eqnarray*}
\varrho_{D^0}(f) 
& := & 
\varrho_{D^0_1}(1)^{b_1}
\cdots
\varrho_{D^0_s}(1)^{b_s} 
f.
\end{eqnarray*}

The shifting family $(\varrho_{D^0})$
defines a quasicoherent sheaf  of 
ideals of $\mathcal{S}$, namely the
sheaf $\mathcal{I}$ generated by all 
sections of the form $f - \varrho_{D^0}(f)$, 
where $f$ is homogeneous  
and $D^0$ runs through $\mathfrak{D}^0$.
Note that $\mathcal{I}$ is a homogeneous 
ideal with respect to the coarsified grading 
$$
\mathcal{S}
\ = \ 
\bigoplus_{[D] \in \Cl(X)} \mathcal{S}_{[D]},
\qquad \qquad 
\mathcal{S}_{[D]}
\ = \ 
\bigoplus_{D' \in D + \mathfrak{D}^0} \mathcal{O}_X(D').
$$
Moreover, it turns out that $\mathcal{I}$ is
a sheaf of radical ideals.
Dividing the $\Cl(X)$-graded $\mathcal{S}$ by the
homogeneous ideal $\mathcal{I}$, 
we obtain a quasicoherent 
sheaf of $\Cl(X)$-graded $\mathcal{O}_X$-algebras, 
the {\em Cox sheaf}: set $\mathcal{R} :=  \mathcal{S}/\mathcal{I}$
let $\pi \colon \mathcal{S} \to \mathcal{R}$ 
be the projection and define the grading by 
$$ 
\mathcal{R}
\ = \ 
\bigoplus_{[D] \in \Cl(X)} \mathcal{R}_{[D]},
\qquad \qquad 
\mathcal{R}_{[D]}
\ := \ 
\bigoplus_{[D] \in \Cl(X)} \pi(\mathcal{S}_{[D]}).
$$

\begin{remark}
\label{rem:secgradvr}
Let $s \colon  \Cl(X)
\to \mathfrak{D}$ be any set-theoretical section of the 
canonical map $\mathfrak{D} \to \Cl(X)$.
Then the projection $\mathcal{S} \to \mathcal{R}$ 
restricts to a canonical isomorphism of sheaves of
$\Cl(X)$-graded vector 
spaces:
\begin{eqnarray*}
\bigoplus_{[D] \in \Cl(X)} \mathcal{S}_{s([D])}
& \cong &
\bigoplus_{[D] \in \Cl(X)} \mathcal{R}_{[D]}.
\end{eqnarray*}
\end{remark}

One can show that, up to isomorphism,
the graded $\mathcal{O}_X$-algebra
$\mathcal{R}$ does not depend on the 
choices of $\mathfrak{D}$ and the shifting family,
compare~\cite[Lemma~3.7]{BeHa3}.
We define the {\em Cox ring\/} $\mathcal{R}(X)$ 
of $X$, also called the {\em total coordinate ring\/} 
of $X$, to be the $\Cl(X)$-graded 
algebra of global sections of the Cox sheaf:
$$
\mathcal{R}(X)
\ := \
\Gamma(X,\mathcal{R})
\ \cong \ 
\Gamma(X,\mathcal{S}) /\Gamma(X,\mathcal{I}).
$$
The sheaf $\mathcal{R}$ defines a
universal torsor $p \colon \rq{X} \to X$
in the following sense.
Suppose that $\mathcal{R}$ is locally of finite type;
this holds for example, if $X$ is locally factorial
or if $\mathcal{R}(X)$ is finitely generated.
Then we may consider the relative spectrum 
\begin{eqnarray*}
\rq{X}
& := & 
\Spec_X(\mathcal{R}).
\end{eqnarray*}
The $\Cl(X)$-grading of the sheaf $\mathcal{R}$ 
defines an action of the diagonalizable group
$H := \Spec(\KK[\Cl(X)])$ on $\rq{X}$, 
and the canonical morphism $p \colon \rq{X} \to X$
is a good quotient, i.e., it is an $H$-invariant
affine morphism satisfying
\begin{eqnarray*}
\mathcal{O}_X 
& = & (p_*\mathcal{O}_{\rq X})^H.
\end{eqnarray*}

In the sequel, we mean by a 
{\em universal torsor\/} for $X$
more generally any good quotient 
$q \colon \mathcal{X} \to X$ for an action 
of $H$ on a variety $ \mathcal{X}$ 
such that there is an equivariant 
isomorphism 
$\imath \colon  \mathcal{X} \to \rq{X}$ with 
$q = p \circ \imath$.
If the Cox ring $\mathcal{R}(X)$ 
is finitely generated, then we define the 
{\em total coordinate space\/}
of $X$ to be any affine $H$-variety 
equivariantly isomorphic to 
$\b{X} = \Spec(\mathcal{R}(X))$
endowed with the $H$-action defined 
by the $\Cl(X)$-grading of $\mathcal{R}(X)$.

We discuss basic properties of 
these constructions.
As usual, we say that a Weil divisor 
$\sum a_D D$, where $D$ runs through 
the irreducible 
hypersurfaces, on an $H$-variety
is $H$-invariant if $a_D = a_{h \mal D}$ 
holds for all $h \in H$. 
Moreover, we say that a closed subset 
$B \subseteq Y$ of a variety $Y$ is small,
in $Y$ if it is of codimension at least two in 
$Y$.

\begin{proposition}
\label{gencoxprops}
Let $X$ be a normal variety with 
$\Gamma(X, \mathcal{O}^*) = \KK^*$ and
finitely generated divisor class group.
Let $\mathcal{R}$ be a sheaf of 
$\Cl(X)$-graded algebras as constructed 
before, assume that the Cox ring 
$\mathcal{R}(X)$ is  finitely generated,
denote by $p \colon \rq{X} \to X$ 
the associated universal torsor
and by $X' \subseteq X$ the set of 
smooth points.
\begin{enumerate}
\item
The Cox ring $\mathcal{R}(X)$ is normal and 
every homogeneous invertible $f \in \mathcal{R}(X)$
is constant.
\item 
The canonical morphism $\rq{X} \to \b{X}$ 
is an open embedding; in particular, $\rq{X}$ is
quasiaffine.
\item
The complement $\b{X} \setminus p^{-1}(X')$ is 
small in $\b{X}$,
and the group $H$ acts freely on the set 
$p^{-1}(X') \subseteq \rq{X}$.
\item
Every $H$-invariant Weil divisor on the total coordinate
space $\b{X}$ is principal. 
\end{enumerate}
\end{proposition}

\begin{lemma}
\label{smoothpull}
Suppose that $X$ is a smooth variety, 
let $\mathfrak{D} \subseteq \WDiv(X)$
be a finitely generated subgroup
and consider the sheaf of 
$\mathfrak{D}$-graded algebras
$$
\mathcal{S}
\ := \ 
\bigoplus_{D \in \mathfrak{D}} \mathcal{S}_D
\qquad\qquad
\mathcal{S}_D
\ := \ 
\mathcal{O}_X(D).
$$
Set $\t{X} := \Spec_X(\mathcal{S})$ and let 
$q \colon \t{X} \to X$ be the canonical morphism.
Then, given $D \in \mathfrak{D}$ and a global 
section $f \in \Gamma(X,\mathcal{O}(D))$, we obtain
\begin{eqnarray*}
q^*(D) 
& = & 
\div(f)  - q^*(\div(f)),
\end{eqnarray*}
where on the right hand side $f$ is firstly
viewed as a homogeneous function on $\t{X}$, 
and secondly as a rational function on $X$.
In particular, $q^*(D)$ is principal.
\end{lemma}

\begin{proof}
On suitable open sets $U_i \subseteq X$,
we find defining equations $f_i^{-1}$ for $D$ 
and thus may write $f = h_i f_i$, where 
$h_i \in \Gamma(U_i,\mathcal{S}_0) =  
\Gamma(U_i,\mathcal{O})$ and
$f_i \in \Gamma(U_i,\mathcal{S}_D)$. 
Then, on $q^{-1}(U_i)$, we have $q^*(h_i) = h_i$
and the function $f_i$ 
is homogeneous of degree $D$ and invertible.
Thus, we obtain
\begin{eqnarray*} 
q^*(D)
& = &
q^*(\div(f)+D) - q^*(\div(f))
\\ 
& = & 
q^*(\div(h_i))  - q^*(\div(f))
\\
& = & 
\div(h_i)  - q^*(\div(f))
\\
& = &  
\div(h_if_i)  - q^*(\div(f))
\\
& = & 
\div(f)  - q^*(\div(f)).
\end{eqnarray*}
\end{proof}

\begin{proof}[Proof of Proposition~\ref{gencoxprops}]
For~(i), let $f \in \mathcal{R}(X)^*$ be homogeneous of
degree $[D]$.
Then its inverse $g \in \mathcal{R}(X)^*$ is 
homogeneous of degree $-[D]$, and
we have $fg \in \mathcal{R}(X)^*_0 = \KK^*$.
According to Remark~\ref{rem:secgradvr}, we may 
view $f$ and $g$ as global sections of sheaves
$\mathcal{O}_X(D)$ and $\mathcal{O}_X(-D)$ 
respectively and thus obtain
$$
0 
\ = \ 
\div(fg)
\ = \ 
(\div(f) + D) + (\div(g) -D).
$$
Since the divsors $(\div(f) + D)$ and $(\div(g) -D)$ 
are both nonnegative, we can conclude 
$D = \div(f^{-1})$ and hence $[D]=0$.
This implies 
$f \in \mathcal{R}(X)^*_0 = \KK^*$ as wanted.
To proceed, note that 
$X$ and $X'$ have the same Cox ring.  
Thus, normality of $\mathcal{R}(X) = \mathcal{R}(X')$ 
follows from~\cite[Prop.~6.3]{BeHa3}.

We prove most of~(iii).
For any affine $U \subseteq X$, 
the set $p^{-1}(U \cap X')$ has the same 
functions as the affine set $p^{-1}(U)$. 
Thus,
$p^{-1}(U) \setminus p^{-1}(U \cap X')$
must be small. 
Consequently also $\rq{X} \setminus p^{-1}(X')$ 
is small.
The fact that $H$ acts freely on $p^{-1}(X')$,
is due to the fact that on $X'$ the homogeneous 
components $\mathcal{R}_{[D]}$ are invertible 
and hence each $H$-orbit in $p^{-1}(X')$
admits an invertible function of any degree
$w \in K$. 

Before proceeding, consider the group 
$\mathfrak{D} \subseteq \WDiv(X)$ and 
the associated sheaf $\mathcal{S}$ as 
used in the definition of $\mathcal{R}$, 
and set $\t{X} := \Spec_X(\mathcal{S})$.
Then the projection $\mathcal{S}\to \mathcal{R}$
defines a commutative diagram 
\begin{equation}
\label{RSdiag} 
\xymatrix{
{\rq{X}} 
\ar[rr]
\ar[dr]_p
& & 
{\t{X}} 
\ar[dl]^q
\\
& X & 
}
\end{equation}

We turn to~(ii). 
Cover $X$ by affine open subsets  
$X_f := X \setminus \Supp(D+\div(f))$,
where $D \in \mathfrak{D}$ and 
$f \in \Gamma(X, \mathcal{O}(D))$.
Then each $f$ defines an element in 
$\mathcal{R}(X)$, and hence a function 
on $\rq{X}$ and as well on $\b{X}$. We claim
\begin{eqnarray*}
p^{-1}(X_f) & = & \rq{X}_f.
\end{eqnarray*}
The inclusion ``$\subseteq$'' follows from the 
observation that the function $f$ is invertible 
on $p^{-1}(X_f)$.
Using Lemma~\ref{smoothpull} and the above commutative 
diagram, we obtain 
\begin{eqnarray*}
p^{-1}(X_f) \cap p^{-1}(X')
& = & 
\rq{X}_f \cap p^{-1}(X').
\end{eqnarray*}
Thus, the complement $\rq{X}_f \setminus p^{-1}(X_f)$ 
is small. Since $p^{-1}(X_f)$ is affine, we 
obtain the desired equality.
Consequently, $\rq{X}$ can be covered by open 
affine subsets $\rq{X}_f$.
The corresponding morphisms $\rq{X}_f \to \b{X}_f$
are isomorphisms and glue together to the desired
open embedding.

Now we easily can finish the proof of~(iii).
Knowing $\rq{X} \subseteq \b{X}$ and that 
both varieties have the same global functions,
we obtain that $\b{X} \setminus \rq{X}$ is 
small. By what we saw before, then also 
$\b{X} \setminus p^{-1}(X')$ is small. 

We prove~(iv). According to~(iii), we may 
assume that $X$ is smooth and we only have 
to show that 
every invariant Weil divisor
$\rq{D}$ on $\rq{X}$ is trivial.
Since $H$ acts freely, we have $\rq{D} = p^*(D)$ 
with a Weil divisor $D$ on $X$.
Thus, we have to show that all pullback
divisors $p^*(D)$ are trivial.
For this, it suffices to 
treat effective divisors
$D$ on $X$, and we may assume that $D$ 
belongs to $\mathfrak{D}$.
In view of the commutative diagram~(\ref{RSdiag}), 
it suffices to show that $q^*(D)$ is principal. 
This in turn follows immediately from
Lemma~\ref{smoothpull}.
\end{proof}

As an important example, we briefly discuss 
the case of a toric variety, where due to a 
combinatorial description, Cox ring and universal
torsor can be constructed explicitly. 
First, recall that a toric 
variety is normal variety $Z$ together
with an action of an algebraic torus 
$T$ and a base point $z_0 \in Z$ such that 
the orbit map $T \mapsto Z$, 
$t \mapsto t \mal z_0$ is an open embedding.

Any toric variety $Z$ arises from a 
fan in a lattice $N$, 
i.e., a finite collection
$\Sigma$ of strictly convex, polyhedral
cones in $N_\QQ = \QQ \otimes_\ZZ N$ 
such that for any $\sigma \in \Sigma$ also every
face $\sigma_{0} \preceq \sigma$ belongs 
to $\Sigma$,
and for any two 
$\sigma_{1}, \sigma_{2} \in \Sigma$
we have 
$\sigma_{1} \cap \sigma_{2} \preceq \sigma_{i}$.
The acting torus of $Z$ is
$T = \Spec(\KK[M])$, 
where $M := \Hom(N,\ZZ)$ denotes the dual 
lattice,
and $Z$ is the equivariant gluing of the 
affine toric varieties
$$
Z_\sigma
\ := \ 
\Spec(\KK[\sigma^\vee \cap M]),
\quad
\sigma \in \Sigma,
$$
where $\sigma^\vee \subseteq M_\QQ$ is the dual
cone, and the $T$-action on $Z_\sigma$ 
is given by the $M$-grading of the semigroup 
algebra $\KK[\sigma^\vee \cap M]$.
There is a bijection from the fan $\Sigma$ 
onto the set of $T_Z$-orbits in~$Z$,
sending a cone $\sigma \in \Sigma$ to 
$\orb(\sigma)$, the unique closed
$T$-orbit in $Z_\sigma$.
If $\varrho_1, \ldots, \varrho_r$
denote the rays of $\Sigma$, 
then the $T$-invariant prime divisors 
of $Z$ are precisely the orbit closures
\begin{eqnarray*}
D_Z^i & := & \b{\orb(\varrho_i)}.
\end{eqnarray*}

Now we recall Cox's construction~\cite{Cox}.  
Suppose that the fan $\Sigma$ is nondegenerate,
i.e., the primitive lattice vectors 
$v_i \in \varrho_i \cap N$ generate 
$N_\QQ$ as a vector space; 
this just means that we have 
$\Gamma(Z, \mathcal{O}^*) = \KK^*$.
Set $F := \ZZ^r$ and consider the 
linear map $P \colon F \to N$ sending the 
$i$-th canonical base vector $e_i \in F$
to $v_i \in N$. 
There is a fan $\rq{\Sigma}$ in $F$ 
consisting of certain faces of the positive 
orthant 
$\delta \subseteq F_{\QQ}$, namely
$$ 
\rq{\Sigma}
\ := \
\{\rq{\sigma} \preceq \delta; \; P(\rq{\sigma}) \subseteq \sigma
\text{ for some } \sigma \in \Sigma\}.
$$
The fan $\rq{\Sigma}$ defines an open toric 
subvariety $\rq{Z}$ of  
$\b{Z} = \Spec(\KK[\delta^\vee \cap E])$,
where $E := \Hom(F,\ZZ)$.
As $P \colon F \to N$ is a map of the fans 
$\rq{\Sigma}$ and $\Sigma$, i.e., sends cones
of $\rq{\Sigma}$ into cones of $\Sigma$,
it defines a morphism $p_Z \colon \rq{Z} \to Z$
of toric varieties. 
Note that for the unions 
$W_{Z} \subseteq Z$ and $W_{\b{Z}} \subseteq \b{Z}$ 
of all at most one codimensional orbits of the 
respective acting tori, we have 
$W_{\b{Z}} = p_{Z}^{-1}(W_{Z}) \subseteq \rq{Z}$.

Let us briefly explain, why  
$p_Z \colon \rq{Z} \to Z$ 
is a universal torsor for $Z$.
First, observe that $P \colon F \to N$
and its dual map $P^* \colon M \to E$ 
give rise to a pair of exact sequences 
of abelian groups:
$$ 
\xymatrix{
0 \ar[r]
&
L \ar[r]^{Q^*}
&
F \ar[r]^{P}
&
N 
&
\\
0
&
K \ar[l]
&
E \ar[l]^{Q}
&
M \ar[l]^{P^*}
&
0  \ar[l]
}
$$
The latter sequence has two geometric interpretations.
Firstly, it relates the character groups of big tori 
$T = \Spec(\KK[M])$ of $Z$
and $\rq{T} = \Spec(\KK[E])$ of $\rq{Z}$
to that of $H = \Spec(\KK[K])$. 
Secondly, the lattices $M$ and $E$ represent 
certain groups of invariant divisors, and
$K$ is the divisor class group of $Z$.
Alltogether, one has the following
commutative diagram. 
\begin{equation}
\label{diag:coxisos} 
\xymatrix{
&
{\Chi(H)} \ar[dl]
&
{\Chi(\rq{T})} \ar[l]
&
{\Chi(T)} \ar[l]^{p_Z^*}
&
\\
0
&
K \ar[l] \ar[u]^\cong \ar[d]_\cong
&
E \ar[l]^{Q}  \ar[u]_\cong^{e \mapsto \chi^e} \ar[d]^\cong_{e_i \mapsto D_i}
&
M \ar[l]^{P^*}  \ar[u]_\cong^{u \mapsto \chi^u} \ar[d]^\cong_{u \mapsto \div(\chi^u)}
&
0  \ar[l] \ar[ul] \ar[dl]
\\
&
{\Cl(Z)}  \ar[ul]
&
{\WDiv^T(Z)} \ar[l] 
&
{\PDiv^T(Z)}\ar[l]
&
}
\end{equation}

To obtain the Cox ring, consider the groups of divisors 
$\mathfrak{D} := \WDiv^T(Z)$ and $\mathfrak{D}^0 := \PDiv^T(Z)$,
and the sheaf of $\mathcal{O}_Z$-algebras $\mathcal{S}$ associated 
to $\mathfrak{D}$. Then there is a {\em canonical shifting family\/}:
for any $\div(\chi^u) \in  \mathfrak{D}^0$ it sends a section
$\chi^{u'}$ to $\chi^{u'-u}$.
Let $\mathcal{I}$ denote the associated sheaf of ideals and set
$\mathcal{R} := \mathcal{S}/\mathcal{I}$. 

Choose any set theoretical section $s \colon K \to E$ 
for $Q \colon E \to K$, and, for $w \in K$, let
$D(w) \in \WDiv^T(Z)$ denote the divisor corresponding to 
$s(w)$.
Moreover, given a maximal $T$-invariant affine open subset 
$W \subseteq Z$, set $\rq{W} := p_Z^{-1}(W)$.
Then, for each $w \in K$ the assignment
$$
\Gamma(W,\mathcal{O}_Z(D(w)))
\ \to \
\Gamma(\rq{W}, \mathcal{O}_{\rq{Z}})_w,
\qquad
\chi^u 
\ \mapsto \ 
\chi^{(P^*(u) + s(w))}
$$
defines an isomorphism of vector spaces.
Using this and Remark~\ref{rem:secgradvr},
we obtain isomorphisms of sheaves 
of graded vector spaces, 
the composition of 
which is even an isomorphism of 
sheaves of graded 
algebras:
\begin{equation} 
\label{toriccoxisos}
\bigoplus_{w \in K} 
\mathcal{R}_{[D(w)]}
\ \cong \ 
\bigoplus_{w \in K} 
\mathcal{O}_Z(D(w))
\ \cong \ 
\bigoplus_{w \in K} 
(p_Z)_*(\mathcal{O}_{\rq{Z}})_w
\end{equation}
Passing to the relative spectra, we 
obtain an $H$-equivariant isomorphism 
$\Spec_Z (\mathcal{R}) \cong \rq{Z}$ 
over $Z$ and thus see that
$p_Z \colon \rq{Z} \to Z$ 
is a universal torsor.
Consequently, the Cox ring of $Z$ is 
the polynomial ring  
$\mathcal{O}(\rq{Z}) = \mathcal{O}(\b{Z})
= \KK[E \cap \delta^{\vee}]$ with the 
$K$-grading given by $\deg (\chi^e) = Q(e)$.

Finally, we note that $H = \Spec(\KK[K])$ acts freely 
on the union $W_{\b{Z}} \subseteq \b{Z}$
of all at most one-codimensional torus orbits.
Indeed, $W_{\b{Z}}$ is covered by the 
affine sets 
$\b{Z}_i := \Spec(\KK[\varrho_i^\vee \cap E])$,
where $\varrho_i = \QQ_{\ge 0} \mal P(e_i)$.
Since $P(e_i) = v_i$ is primitive,
the images $Q(e_j) \in K$, where 
$j \ne i$, generate $K$ as an abelian group.
Since the corresponding $\chi^{e_j}$ 
are invertible functions on $\b{Z}_i$, 
we may construct  invertible  
$f \in \mathcal{O}(\b{Z}_i)_w$
for any $w \in K$.
This gives freeness of the $H$-action on 
$W_{\b{Z}} \subseteq \b{Z}$.

Now we use Cox's construction to 
determine universal torsor and total
coordinate space of certain 
``neatly embedded'' subvarieties
$X \subseteq Z$ of our toric variety $Z$.
To make this precise, 
we need the following pullback 
construction for invariant Weil 
divisors of $Z$.

\begin{remark}
\label{weilpull}
Let $X$ be an irreducible variety,
which is smooth in codimension one,
and let $\varphi \colon X \to Z$  
a morphism to a toric variety $Z$ 
such that 
$\varphi(X)$ intersects the big torus 
orbit of $Z$ and
$X \setminus \varphi^{-1}(W_Z)$ contains no divisors 
of $X$.
Then the usual pullback of invariant 
Cartier divisors over $W_Z$ canonically
defines a homomorphism
$$ 
\varphi^* \colon \WDiv^T(Z) \ \to \ \WDiv(X).
$$
Note that $\varphi^*$ takes nonnegative divisors 
to nonnegative ones.
Moreover, one has canonical pullback homomorphisms
$\varphi^*\mathcal{O}(D) \to \mathcal{O}(\varphi^*D)$
for any  $D \in \WDiv^T(Z)$.
Finally, since $\Cl(Z)$ 
is obtained as $\WDiv^T(Z)/\PDiv^T(Z)$,
there is an induced pull back
homomorphism for the divisor class groups
$$
\varphi^{*} \colon \Cl(Z) \ \to \ \Cl(X).
$$
Note that analogous pullback constructions can be performed 
for any dominant morphism  $\varphi \colon X \to Y$ of
varieties $X$ and $Y$, which are smooth in codimension one, 
provided that for the smooth locus $Y' \subseteq Y$ the 
inverse image $\varphi^{-1}(Y')$ has 
a small complement in $X$.
\end{remark}

\begin{definition}
\label{def:neatemb}
Let $Z$ be a toric variety  with
acting torus $T$ and invariant prime 
divisors $D_Z^i = \b{T \mal z_i}$,
where $1 \le i \le r$.
Moreover, let $X \subseteq Z$ be an irreducible 
closed subvariety and suppose that $X$ is 
smooth in codimension one.
We call $X \subseteq Z$ a 
{\em neat embedding\/} if 
\begin{enumerate}
\item
each $D_X^i := D_Z^i \cap X$ is an irreducible hypersurface 
in $X$ intersecting the toric orbit $T \mal z_i$,
\item
the pull back homomorphism  
$\imath^{*} \colon \Cl(Z) \to \Cl(X)$
defined by the inclusion $\imath \colon X \to Z$
is an isomorphism.
\end{enumerate}
\end{definition}

As we shall see in Section~\ref{sec:bring}
any normal variety with finitely generated Cox
ring that admits a closed embedding into a 
toric variety, admits even a neat embedding into a 
toric variety.

Now consider a toric variety $Z$ with 
$\Gamma(Z, \mathcal{O}^*) = \KK^*$
and toric universal torsor 
$p_Z \colon \rq{Z} \to Z$ as constructed 
before.
Let $X \subseteq Z$ be an irreducible closed 
subvariety.
Set $\rq{X} := p_{Z}^{-1}(X)$,
let $p_X \colon \rq{X} \to X$ denote the 
restriction of $p_Z$, and 
define $\b{X}$ to be the closure
of $\rq{X}$ in $\b{Z}$. 
Then everything fits into a commutative 
diagram: 
$$
\xymatrix{
{\b{X}} \ar[r]^{\b{\imath}} 
&
{\b{Z}} 
\\
{\rq{X}} \ar[r]^{\rq{\imath}} \ar[d]_{p_{X}} \ar[u]
&
{\rq{Z}} \ar[d]^{p_{Z}} \ar[u]
\\
X \ar[r]_{\imath}
&
Z
}
$$
where $\imath, \rq{\imath}$ and $\b{\imath}$ denote
the respective inclusions.
The subvarieties $\rq{X} \subseteq \rq{Z}$ 
and $\b{X} \subseteq \b{Z}$ 
are invariant with respect to the
action of $H$ on $\rq{Z}$ and $\b{Z}$.
In particular, the rings of functions
$\Gamma(\mathcal{O}, \b{X})$
and
$\Gamma(\mathcal{O},\rq{X})$
are graded by the abelian group $K \cong \Cl(Z)$. 
Set 
$$
W_X \ := \ X \setminus \bigcup_{i \ne j} D_X^i \cap D_X^j,
\qquad\qquad
W_{\b{X}} \ := \ p_X^{-1}(W_X).
$$

\begin{theorem}
\label{ringcalc}
Let $X$ be a normal variety and
$X \subseteq Z$ neat embedding,
suppose that $X$ and $Z$ only 
have constant globally invertible 
regular functions and that $X = W_X$ holds. 
Then there is an isomorphism of $K$-graded
$\mathcal{O}_X$-algebras
\begin{eqnarray*}
\mathcal{R}
& \cong & 
(p_X)_* \mathcal{O}_{\rq{X}},
\end{eqnarray*}
where $\mathcal{R}$ denotes the Cox sheaf on $X$.
Moreover, the variety $\rq{X} = W_{\b X}$ 
is normal and the 
restricion $p_X \colon \rq{X} \to X$ 
is a universal torsor for $X$.
\end{theorem}

\begin{proof}
We may assume that $Z = W_Z$ holds.
We define the Cox ring of $Z$ as before, using the 
group $\mathfrak{D}_Z := \WDiv^T(Z)$.
Then we have $\mathfrak{D}^0_Z = \PDiv^T(Z)$
and the associated Cox sheaf
$\mathcal{R}_Z = \mathcal{S}_Z / \mathcal{I}_Z$ 
of $\Cl(Z)$-graded algebras is built up by means 
of the canonical shifting family.

Let $\mathfrak{D}_X := \imath^* \mathfrak{D}_Z$.
By the definition of a neat embedding, 
we have $\mathfrak{D}_X^0 = \imath^*  \mathfrak{D}_Z^0$,
and we may define the Cox sheaf
$\mathcal{R}_X = \mathcal{S}_X / \mathcal{I}_X$ on $X$
via the pullback shifting family.
This gives a commutative diagram
$$ 
\xymatrix{
{\imath^* \mathcal{S}_{Z}}
\ar[r]
\ar[d]
&
{\imath^* \mathcal{R}_{Z}}
\ar[d]
\\
{\mathcal{S}_{X}}
\ar[r]
&
{\mathcal{R}_{X}}
}
$$
Note that the downwards arrows 
of the above diagram are even isomorphisms,
because by $Z = W_Z$ all divisors of 
$\mathfrak{D}_Z$ are locally principal.
Now, choose any set theoretical section 
$s \colon \Cl(Z) \to  \mathfrak{D}$.
Then the isomorphisms~(\ref{toriccoxisos}) define a commutative 
diagram providing the desired isomorphism:
$$ 
\xymatrix{
{\imath^*\mathcal{R}_{Z}}
\ar[d]_\cong
\ar[r]^{\cong \qquad \quad}
&
{\imath^* \bigoplus \mathcal{O}_{Z}(s([D]))}
\ar[r]^{\quad \cong}
\ar[d]
&
{\imath^*(p_Z)_* \mathcal{O}_{\rq{Z}}}
\ar[d]^\cong
\\
{\mathcal{R}_{X}}
\ar[r]
&
{\bigoplus \mathcal{O}_{X}(\imath^*s([D]))}
\ar[r]
&
{(p_X)_* \mathcal{O}_{\rq{X}}}
}
$$

Finally, to see that $\rq{X}$ is normal,
note that the group $H$ acts freely 
on $\rq{Z} = W_{\b{Z}}$.
By Luna's slice theorem, $\rq{X}$ looks 
locally, in etale topology,
like $H \times X$.
Since $X$ is normal, and normality is preserved 
under etale maps, we obtain that $\rq{X}$
is normal.
\end{proof}

We give two applications.
Recall that a variety 
$X$ is said to be $\QQ$-factorial 
if it is normal and for any Weil divisor 
$D$ on $X$ some positive multiple $nD$ 
is Cartier.
A toric variety $Z$ is $\QQ$-factorial if
and only if it arises from a simplicial
fan $\Sigma$, i.e., each cone of $\Sigma$
is generated by a linearly independent 
set of vectors. 
A toric variety $Z$ with Cox 
construction $\rq{Z} \to Z$ is 
$\QQ$-factorial if and only if
$\rq{Z} \to Z$ is a geometric quotient,
i.e., its fibers are precisely the 
$H$-orbits.

\begin{corollary}
\label{cor:qfact2tosor}
Let $X \subseteq Z$ be a neat embedding, where $X$ 
is irreducible and smooth in codimension one and
$Z$ is a $\QQ$-factorial toric variety having
only constant globally invertible functions.
\begin{enumerate}
\item
Suppose that $\Gamma(X, \mathcal{O}^*) = \KK^*$ 
holds and that $\rq{X}$ is normal.
Then $X$ is normal,
$p_X \colon \rq{X} \to X$ 
is a universal torsor for $X$, 
and the Cox ring of $X$ is given by 
$$
\mathcal{R}(X) 
\ = \
\bigoplus_{w \in K} \Gamma(X,\mathcal{R})_{[\imath^*D(w)]}
\ \cong \
\bigoplus_{w \in K} \Gamma(\rq{X},\mathcal{O})_w
\ = \
\Gamma(\rq{X},\mathcal{O}).
$$
\item 
Let $W_X$ be normal with 
$\Gamma(W_X, \mathcal{O}^*) = \KK^*$ 
and let $X^{\rm {nor}}$ be the 
normalization of $X$.
If $\b{X} \setminus \rq{X}$ is small, then
the Cox ring $\mathcal{R}(X^{\rm {nor}})$ 
is finitely generated and $X^{\rm {nor}}$ 
has the $H$-equivariant 
normalization $\b{X}^{\rm {nor}}$ of $\b{X}$ 
as its total coordinate space.
\end{enumerate}
\end{corollary}

\begin{proof}
Since $Z$ is $\QQ$-factorial, each fiber of
$p_Z \colon \rq{Z} \to Z$ consists of a 
single $H$-orbit, having an at most finite 
isotropy group.  
Consequently, the map $p_X \colon \rq{X} \to X$ 
is equidimensional.
Since $X \setminus W_X$ is small in $X$,
we can conclude that $\rq{X} \setminus W_{\b{X}}$ is 
small in $\rq{X}$.

We prove~(i). Since $\rq{X}$ is normal, the 
quotient space $X = \rq{X} \quot H$ is normal 
as well. 
Moreover, by the above consideration, the 
isomorphism of $K$-graded sheaves
of Theorem~\ref{ringcalc} extends from 
$W_X$ to the whole $X$ and the assertion follows.

We turn to~(ii). First note that we have 
$W_X \subseteq  X^{\rm {nor}}$ and that the 
complement $X^{\rm {nor}} \setminus W_X$ is 
small.
Since $\b{X} \setminus \rq{X}$ is small
in $\b{X}$,
the above consideration gives that  
$\b{X} \setminus W_{\b{X}}$ is small
in $\b{X}$.
Consequently, 
$\b{X}^{\rm {nor}} \setminus W_{\b{X}}$ is small
in $\b{X}^{\rm {nor}}$. 
The assertion then follows from
$$
\mathcal{R}(X^{\rm {nor}})
\ = \ 
\Gamma(X^{\rm {nor}},\mathcal{R}) 
\ = \ 
\Gamma(W_X,\mathcal{R}) 
\ = \ 
\Gamma(W_{\b{X}},\mathcal{O}) 
\ = \ 
\Gamma(\b{X}^{\rm {nor}},\mathcal{O}).
$$
\end{proof}

\begin{corollary}
\label{cor:smooth2tosor}
Let $X \subseteq Z$ be a neat embedding of 
a normal variety in a smooth toric variety $Z$, 
and suppose that $X$ as well as $Z$ only admit
constant globally invertible functions.
Then $p_X \colon \rq{X} \to X$ 
is a universal torsor for $X$, 
and the Cox ring of $X$ is given by 
$$
\mathcal{R}(X) 
\ = \
\bigoplus_{w \in K} \Gamma(X,\mathcal{R})_{[\imath^*D(w)]}
\ \cong \
\bigoplus_{w \in K} \Gamma(\rq{X},\mathcal{O})_w
\ = \
\Gamma(\rq{X},\mathcal{O}).
$$ 
Moreover, if $\b{X} \setminus \rq{X}$ 
contains no divisors of $\b{X}$, 
then $\mathcal{R}(X)$ is finitely 
generated and $X$ has the $H$-equivariant 
normalization of $\b{X}$ as its total coordinate
space.
\end{corollary}

\section{Bunched rings and their varieties}
\label{sec:bring}

In this section, we generalize the construction
of varieties with 
a prescribed finitely generated Cox
ring provided in~\cite{BeHa1}
to the case of divisor class groups with torsion,
and we present first basic features of this generalized 
construction.
We begin with discussing the class of rings,
which we will consider; 
the concept~\ref{def:factgrad} as well as 
statements~\ref{factgradchar} 
and~\ref{factgrad2fact} are basically due to 
I.V.~Arzhantsev, see~\cite{Ar2}.

\begin{definition}
\label{def:factgrad}
Let $K$ be a finitely generated abelian 
group and $R = \bigoplus_{w \in K} R_w$ 
any $K$-graded integral $\KK$-algebra with
$R^* = \KK^*$. 
\begin{enumerate}
\item
We say that a nonzero nonunit $f \in R$ 
is {\em $K$-prime} if it is homogeneous and 
$f \vert gh$ with homogeneous $g,h \in R$
always implies $f \vert g$ or $f \vert h$.
\item
We say that an ideal $\mathfrak{a} \subset R$
is {\em $K$-prime\/} if it is homogeneous
and for any two homogeneous 
$f,g \in R$ with $fg \in \mathfrak{a}$ 
one has $f \in \mathfrak{a}$ or $g \in \mathfrak{a}$.
\item
We say that a homogeneous prime ideal 
$\mathfrak{a} \subset R$ has 
{\em $K$-height $d$\/} if $d$ is maximal 
admitting a chain 
$\mathfrak{a}_0 \subset \mathfrak{a}_1 \subset \cdots 
\subset \mathfrak{a}_d = \mathfrak{a}$
of $K$-prime ideals.
\item
We say that the ring $R$ is {\em factorially graded} if 
every  $K$-prime ideal of $K$-height 
one is principal.
\end{enumerate}
\end{definition}

We will briefly indicate the geometric meaning
of these notions.
By an {\em $H$-prime divisor\/} on an $H$-variety 
we mean a sum $\sum a_D D$ 
with prime divisors $D$ such that
always $a_D = 0,1$ holds and the
$D$ with $a_D =1$ are transitively 
permuted  by the group $H$.
Note that every $H$-invariant divisor
is a unique sum of $H$-prime divisors.

\begin{proposition}
\label{factgradchar}
Let $K$ be a finitely generated 
abelian group, $R$ a $K$-graded 
normal affine $\KK$-algebra with 
$R^* = \KK^*$, and consider the 
action of $H := \Spec(\KK[K])$ on  
$\b{X} := \Spec(R)$.
Then the following statements are 
equivalent.
\begin{enumerate}
\item
The ring $R$ is factorially graded.
\item
Every invariant Weil divisor of $\b{X}$ 
is principal.
\item
Every homogeneous $0 \ne f \in R \setminus R^*$ 
is a product of $K$-primes.
\end{enumerate}
Moreover, if one of these statements holds,
then a homogeneous nonzero nonunit
$f \in R$ is $K$-prime 
if and only if $\div(f)$ is 
$H$-prime,
and every $H$-prime divisor is of the form
$\div(f)$ with a $K$-prime $f \in R$.
\end{proposition}

\begin{proof}
Assume that~(i) holds, and let $D$ be an invariant
Weil divisor on $\b{X}$. 
Write $D = D_1 + \cdots + D_r$ with $H$-prime divisors 
$D_i$. 
Then the vanishing ideal $\mathfrak{a}_i$ of $D_i$ 
is $K$-prime of $K$-height one and thus of the form 
$\mathfrak{a}_i = \bangle{f_i}$. 
We obtain $D_i = \div(f_i)$ and thus 
$D = \div(f_1 \cdots f_r)$,
which verifies~(ii).

Assume that~(ii) holds. Given a homogeneous 
$0 \ne f \in R \setminus R^*$, write 
$\div(f) = D_1 + \cdots + D_r$ with $H$-prime divisors 
$D_i$. 
Then we have $D_i = \div(f_i)$, and one directly 
verifies that the $f_i$ are homogeneous $K$-primes.
This gives a decomposition 
$f = a \cdot f_1  \cdots  f_r$ with some 
$a \in R^* = \KK^*$ as required in (i).

If~(iii) holds and $\mathfrak{a}$ is a K-prime ideal
of height one, then we take any homogeneous
$0 \ne f \in \mathfrak{a}$ and find a $K$-prime 
factor $f_1$ of $f$ with $f_1 \in \mathfrak{a}$.
Thus, we have inclusions
$0 \subset \bangle{f_1} \subseteq \mathfrak{a}$
of $K$-prime ideals, which implies 
$\mathfrak{a} = \bangle{f_1}$,
verifying~(i).
\end{proof}

\begin{remark}
\label{factgrad2fact}
A normal factorially $K$-graded affine
$\KK$-algebra $R$ with a free finitely 
generated abelian group $K$ correponds 
to a normal affine variety $\b{X}$ with 
an action of a torus $H$. 
In this case, every Weil divisor on $\b{X}$ 
is $H$-linearizable and hence linear equivalent 
to an $H$-invariant one and thus principal,
i.e., $R$ is even factorial. 
However, if $K$ has torsion, then
there may exist normal affine 
algebras which are factorially 
$K$-graded but not factorial,
see~\cite[Example~4.2]{Ar2}.
\end{remark}

For a finitely generated abelian group $K$,
we denote by $K^t \subseteq K$ the torsion part,
set $K^0 := K/K^t$ and write $w^0 \in K^0$ 
for the class of $w \in K$.
Given a free abelian group $E$ 
and a homomorphism $Q \colon E \to K$, 
we denote by $Q^0 \colon E_\QQ \to K^0_\QQ$
the induced linear map sending 
$a \otimes e$ to $a \otimes Q(e)^0$. 
Moreover, the relative interior of a cone
$\sigma \subseteq V$ in a rational 
vector space $V$ is denoted by $\sigma^\circ$.

\begin{definition}
\label{bringdef}
Consider a finitely generated abelian group $K$ 
and a normal factorially $K$-graded affine $\KK$-algebra 
$R$.
Let 
$$
\mathfrak{F} 
\ = \ 
\{f_{1}, \dots, f_{r}\} 
\ \subset \ 
R
$$ 
be a system of homogeneous pairwise non 
associated $K$-prime generators for $R$.
The {\em projected cone associated\/} 
to $\mathfrak{F} \subset R$ is 
$$(E \topto{Q} K, \gamma),$$
where $ E := \ZZ^{r}$, the homomorphism
$Q \colon E \to K$ sends
$e_{i} \in E$ to  
$w_i := \deg(f_{i}) \in K$
and $\gamma \subseteq  E_{\QQ}$ is the 
convex cone generated by $e_{1}, \dots, e_{r}$.
\begin{enumerate}    
\item
We say that $\mathfrak{F} \subset R$ is {\em admissible},
if, for each facet $\gamma_{0} \preceq \gamma$,
the image $Q(\gamma_{0} \cap E)$ generates 
the abelian group $K$.
\item \label{fface}    
A face $\gamma_{0} \preceq \gamma$ is called an
{\em $\mathfrak{F}$-face\/} if the product
over all $f_{i}$ with $e_{i} \in \gamma_{0}$
does not belong to the ideal
$\sqrt{\bangle{f_{j}; \; e_{j} \not\in \gamma_{0}}} \subseteq R$.
\item \label{fbunch}
An {\em $\mathfrak{F}$-bunch\/} is  a nonempty 
collection $\Phi$ of projected $\mathfrak{F}$-faces 
$Q^0(\gamma_0) \subseteq K_\QQ$ with the 
following properties:
\begin{enumerate}
\item 
a projected $\mathfrak{F}$-face $\tau \subseteq K_\QQ$ 
belongs to $\Phi$ if 
and only if for each $\tau \neq \sigma \in \Phi$ we have
$\emptyset \neq \tau^{\circ} \cap \sigma^{\circ} \neq 
\sigma^{\circ}$,
\item 
for each facet $\gamma_{0} \prec \gamma$,
there is a cone $\tau \in \Phi$ such that 
$\tau^{\circ} \subseteq Q^0(\gamma_{0})^{\circ}$ 
holds.
\end{enumerate}
\end{enumerate}
A {\em bunched ring\/} is a triple 
$(R,\mathfrak{F},\Phi)$, where 
$\mathfrak{F} \subset R$ is an
admissible system of generators
and $\Phi$ is an $\mathfrak{F}$-bunch.
\end{definition}

\begin{construction}
\label{bringvarconstr}
Let $(R,\mathfrak{F},\Phi)$ be a bunched ring 
and $(E \topto{Q} K,\gamma)$ its projected cone.
The associated {\em collection of relevant faces\/} 
and {\em covering collection\/} are
\begin{eqnarray*}
\rel(\Phi)
& := & 
\{
\gamma_0 \preceq \gamma; \; 
\gamma_0 \text{ an $\mathfrak{F}$-face with }
\tau^\circ \subseteq Q^0(\gamma_0)^\circ 
\text{ for some } \tau \in \Phi
\},
\\
\cov(\Phi)
& := & 
\{
\gamma_0 \in \rel(\Phi); \; 
\gamma_0 \text{ minimal} \}.
\end{eqnarray*}
Consider the action of the diagonalizable group 
$H := \Spec(\KK[K])$ on $\b{X} := \Spec(R)$,
and the open $H$-invariant subsets of $\b{X}$ 
defined by
$$ 
\rq{X}(R,\mathfrak{F},\Phi) 
\ := \ 
\bigcup_{\gamma_0 \in \rel(\Phi)} \b{X}_{\gamma_0},
\qquad
 \b{X}_{\gamma_0} \ := \  \b{X}_{f_1^{u_1} \cdots f_r^{u_r}}
\text{ for some }
(u_1, \ldots, u_r) \ \in \ \gamma_0^\circ.
$$
Then the $H$-action on 
$\rq{X} := \rq{X}(R,\mathfrak{F},\Phi)$ 
admits a good quotient; we denote the 
quotient variety by
$$ 
X
\ := \ 
X(R,\mathfrak{F},\Phi)
\ := \
\rq{X}(R,\mathfrak{F},\Phi) \quot H.
$$
The subsets 
$\b{X}_{\gamma_0} \subseteq \rq{X}$,
where $\gamma_0 \in \rel(\Phi)$,
are saturated with respect to  
the quotient map $p_X \colon \rq{X} \to X$, 
and $X$ is covered by the affine open subsets 
$X_{\gamma_0} := p_X(\b{X}_{\gamma_0})$.
\end{construction}

Before proving the statements made in this 
construction, we present its basic features. 
Recall that an {\em $A_2$-variety\/} is a variety, 
in which
any two points admit a common affine neighbourhood.

\begin{theorem}
\label{thm:bring2var}
Let $(R,\mathfrak{F},\Phi)$ be a bunched ring
with projected cone $(E \topto{Q} K,\gamma)$,
and denote
$\rq{X} := \rq{X}(R,\mathfrak{F},\Phi)$
and 
$X := X(R,\mathfrak{F},\Phi)$.
Then $X$ is a normal $A_2$-variety with 
$$
\dim(X) \ = \ \dim(R) - \dim(K^0_\QQ),
\qquad\qquad
\Gamma(X,\mathcal{O}^*) = \KK^*,
$$
and $X$ admits a neat embedding 
into a toric variety.
Moreover, there is a canonical isomorphism 
$\Cl(X) \cong K$, the Cox ring 
$\mathcal{R}(X)$ is isomorphic
to $R$ and $p_X \colon \rq{X} \to X$ 
is a universal torsor for $X$.
\end{theorem}

We come to the proofs of the assertions.
The key is a combinatorial 
understanding of the variation of good 
quotients for diagonalizable group actions 
on affine varieties.
So, let $K$ be a finitely generated abelian group,
$R$ a $K$-graded affine $\KK$-algebra,
and consider the action of $H := \Spec(\KK[K])$ 
on $\b{X} := \Spec(R)$. 

We say that an open subset $W \subseteq \b{X}$ is 
{\em ($H$,2)-maximal\/} if it admits a good quotient
$W \to W \quot H$ with an $A_2$-variety $W \quot H$
and there is no $W' \subseteq \b{X}$ 
with the same property comprising $W$ as a proper subset
and saturated under $W' \to W' \quot H$.
Our task is to describe the ($H$,2)-maximal sets.
We define 
the {\em weight cone\/} $\omega_H(\b{X}) \subseteq K^0_\QQ$
and, for any $z \in \b{X}$,
its {\em orbit cone\/} $\omega_H(z) \subseteq K^0_\QQ$ 
as
\begin{eqnarray*}
\omega_H(\b{X})
& := &
\cone(w^0; \; w \in K \text{ with } R_w \ne 0),
\\
\omega_H(z)
& := &
\cone(w^0; \; w \in K \text{ with } f(z) \ne 0 \text{ for some } f \in R_w).
\end{eqnarray*}

\begin{remark}
For a bunched ring $(R,\mathfrak{F},\Phi)$ 
with projected cone $(E \topto{Q} K, \gamma)$,
the weight cone of the action of 
$H = \Spec(\KK[K])$ on $\b{X} = \Spec(R)$
is $\omega_H(\b{X}) = Q^0(\gamma)$, and 
the orbit cones are precisely the cones $Q^0(\gamma_0)$, 
where $\gamma_0 \preceq \gamma$ is an 
$\mathfrak{F}$-face.
\end{remark}

The orbit cone of a point $z \in \b{X}$ describes the 
set of $H$-orbits in the orbit closure
$\b{H \mal z}$. 
Note that this set becomes partially ordered by 
writing $H \mal z_1 < H \mal z_0$ if 
$H \mal z_1$ is contained in the closure
of $H \mal z_0$.

\begin{proposition}
\label{orbitcones}
For any point $z \in \b{X}$, 
there is an 
order preserving bijection
of finite sets
\begin{eqnarray*}
\{H \text{-orbits in } \b{H \mal z}\}
& \longrightarrow & 
\faces(\omega_H(z)),
\\
H \mal z_0
& \mapsto & 
\omega_H(z_0),
\end{eqnarray*}
and the inverse map is obtained as follows:
given $\omega_0 \preceq \omega_H(z)$, take 
any $w^0 \in \omega_0^\circ$ admitting an 
$f \in R_w$ with $f(z) \ne 0$ and assign
the unique closed orbit 
$H \mal z_0$ in $\b{H \mal z} \cap \b{X}_f$
to $\omega_0$.
\end{proposition}

\begin{proof}
Split $H$ as $H = H^0 \times G$ with the 
component of identity $H^0 = \Spec(\KK[K^0])$ 
and a finite abelian group $G$.
The $H^0$ acts on $\b{X} / G$ 
such that the quotient map 
$\pi \colon \b{X} \to \b{X} / G$ 
becomes $H^0$-equivariant.
For any $z \in \b{X}$, the orbit cones
$\omega_H(z)$ and $\omega_{H^0}(\pi(z))$ coincide,
and we have $H \mal z_0 \subseteq H \mal z$
if and only if 
$H^0 \mal \pi(z_0) \subseteq H^0 \mal \pi(z)$.
Thus, we may assume that $H=H^0$ holds,
i.e. that $H$ is a torus. 
Now observe that it suffices to prove the 
assertions in the case $\b{X} = \b{H \mal z}$.
But then they are basic facts on (not necessarily 
normal) affine toric varieties.
\end{proof}

The set $\Omega_H(\b{X})$ of all orbit cones is finite.
By a {\em 2-maximal\/} collection of orbit cones,
we mean a subset $\Psi \subseteq \Omega_H(\b{X})$ being
maximal with the property that for any two 
$\omega_1,\omega_2 \in \Psi$ one has
$\omega_1^\circ \cap \omega_2^\circ \ne \emptyset$.
Moreover, for any $w \in K$, we define its 
{\em GIT-cone\/} to be the convex polyhedral cone 
$$
\lambda(w) 
\ := \ 
\bigcap_{w^0 \in \omega \in \Omega_H(\b{X})}
\omega
\ \subseteq \
K^0_\QQ.
$$

\begin{proposition}
\label{H2quots}
Let $K$ be a finitely generated abelian group,
$R$ a factorially $K$-graded affine $\KK$-algebra,
and consider the action of $H := \Spec(\KK[K])$ 
on $\b{X} := \Spec(R)$.
Then there are mutually inverse bijections
\begin{eqnarray*}
\{\text{2-maximal collections } \in \Omega_H(\b{X})\}
& \longleftrightarrow &
\{(H,2) \text{-maximal subsets of } \b{X}\},
\\
\Psi 
& \mapsto & 
\{z \in \b{X}; \; \omega_0 \preceq \omega(z) 
                   \text{ for some } \omega_0 \in \Psi\},
\\
\{\omega(z); \; H \mal z \text{ closed in } W\}
& \mapsfrom & 
W.
\end{eqnarray*}
Moreover, the collection of GIT-cones 
$\Lambda(\b{X}) = \{\lambda(w); \; w \in K\}$ is a 
fan in $K^0_\QQ$ having the weight cone 
$\omega_H(\b{X})$ as its support, 
and there is a canonical injection
\begin{eqnarray*}
\Lambda 
& \longrightarrow & 
\{\text{2-maximal collections } \in \Omega_H(\b{X})\},
\\
\lambda 
& \mapsto &
\{\omega \in \Omega_H(\b{X}); \; \lambda^\circ \subseteq \omega^\circ \}
\end{eqnarray*}
The ($H$,2)-maximal open subset of 
$\rq{X}(\lambda) \subseteq \b{X}$ 
arising in this way from $\lambda \in \Lambda(\b{X})$ 
is the set of semistable points
of any $w \in K$ with $w^0 \in \lambda^{\circ}$:
$$ 
\rq{X}(\lambda)
\ = \ 
\b{X}^{ss}(w) 
\ := \ 
\bigcup_{f \in \oplus_{n>0} R_{n \mal w}} \b{X}_f
\ = \ 
\{z \in \b{X}; \; \lambda \subseteq \omega(z)\}.
$$
In particular, the ($H$,2)-maximal subsets arising 
from the $\lambda \in \Lambda$ are precisely those 
having a quasiprojective quotient space.
Moreover, one has
\begin{eqnarray*}
\lambda \ \preceq \ \lambda'
& \iff & 
\rq{X}(\lambda) \ \supseteq \ \rq{X}(\lambda').
\end{eqnarray*}
Finally, the quotient spaces $\rq{X}(\lambda) \quot H$ are
all projective over $\Spec(R_0)$; in particular, they are 
projective if and only if $R_0 = \KK$ holds.
\end{proposition}

The fan $\Lambda(\b{X})$ of all GIT-cones is called the 
{\em GIT-fan\/} of the $H$-variety $\b{X}$, and its cones
are also referred to as {\em chambers\/}.  

\begin{proof}[Proof of Proposition~\ref{H2quots}]
As earlier, choose a splitting $H = H^0 \times G$ 
into the component of identity $H^0 = \Spec(\KK[K^0])$ 
and a finite abelian subgroup $G \subseteq H$.
Then $\b{X}/G$ comes with an 
induced $H^0$-action making 
the quotient map $\pi \colon \b{X} \to \b{X}/G$
equivariant.
The $(H^0,2)$-maximal subsets of $\b{X}/G$ are in
bijection with the ($H$,2)-maximal subsets of 
$\b{X}$ via $V \mapsto \pi^{-1}(V)$.
Moreover, for any $z \in \b{X}$, the orbit cones
$\omega_H(z)$ and $\omega_{H^0}(\pi(z))$ coincide.

By these considerations it suffices to prove 
the assertion for the action of $H^0$ on $\b{X}/G$.
This is done in~\cite[Sec.~1]{ArHa}.
Note that the assumption of factoriality 
posed in \cite[Prop.~1.8]{ArHa} can be weakened to 
requiring that for every $H^0$-invariant Weil 
divisor $D$ on $\b{X}/G$ some multiple is 
principal, which in turn is easily verified by 
considering the pull back divisor $\pi^*(D)$ 
on $\b{X}$.
\end{proof}

\begin{lemma}
\label{bunch22max}
In the notation of Construction~\ref{bringvarconstr},
the open subset 
$\rq{X} \subseteq \b{X}$ equals 
the ($H$,2)-maximal open subset
$W(\Psi) \subseteq \b{X}$ defined by the 
2-maximal collection
$$
\Psi
\ := \
\{\omega \in \Omega_H(\b{X}); \;  
\tau^\circ \subseteq \omega^\circ \text{ for some } \tau \in \Phi
\}
\ = \
\{Q^0(\gamma_0); \; \gamma_0 \in \rel(\Phi)\}.
$$
Moreover, an orbit $H \mal z \subseteq \rq{X}$ 
is closed in $\rq{X}$ if and only if we have 
$\omega_H(z) = Q^0(\gamma_0)$ for some 
$\gamma_0 \in \rel(\Phi)$.
\end{lemma}

\begin{proof}
Let $z \in \b{X}_{\gamma_0}$ for some $\gamma_0 \in \rel(\Phi)$.
Since 
$\b{X}_{\gamma_0} = \b{X}_{f_1^{u_1} \cdots f_r^{u_r}}$ 
does not change under variation of $u \in \gamma_0^\circ$, we obtain 
$Q^0(\gamma_0) \subseteq \omega(z)$. 
Thus, $Q^0(\gamma_0)^\circ \subseteq \omega(z_0)^\circ$
holds for some $z_0$ in the closure of $H \mal z$.
By 2-maximality, we have $\omega(z_0) \in \Psi$.
This implies $z \in W(\Psi)$.
Conversely, given $z \in W(\Psi)$,
the orbit cone $\omega(z)$ comprises some 
$\tau \in \Phi$.
Let $z_0$ be a point in the closure of 
$H \mal z$ in $\b{X}$ with 
$\tau^\circ \subseteq \omega_H(z_0)^\circ$.
Then $\gamma_0 := \cone(e_i; \; f_i(z_0) \ne 0)$
is a relevant face with $z \in \b{X}_{\gamma_0}$.

Finally, let $H \mal z$ be any orbit in $\rq{X}$.
If $H \mal z$ is closed in $\rq{X}$, then we 
have $\omega_H(z) \in \Psi$ and hence 
$\omega_H(z) = Q^0(\gamma_0)$ 
for some $\gamma_0 \in \rel(\Phi)$.
Conversely, let $\omega_H(z) = Q^0(\gamma_0)$ 
hold for some $\gamma_0 \in \rel(\Phi)$.
Then we have $\omega_H(z) \in \Psi$.
If  $H \mal z$ were not closed in $\rq{X}$, 
we had $\omega_H(z_0) \prec \omega_H(z)$
with some $z_0 \in \rq{X}$ such that 
$H \mal z_0$ is closed in $\rq{X}$,
which implies $\omega_H(z_0) \in \Psi$.
Thus, we obtain 
$\omega_H(z_0)^\circ \cap \omega_H(z)^\circ \ne \emptyset$,
a contradiction. 
\end{proof}

\begin{proof}[Proof of Construction~\ref{bringvarconstr} 
and Theorem~\ref{thm:bring2var}, part~1]
By Lemma~\ref{bunch22max} and Proposition~\ref{H2quots},
the good quotient space $X = X(R,\mathfrak{F},\Phi)$ 
exists and is an $A_2$-variety.
It inherits normality from $\rq{X}$,
and we have
$$
\Gamma(X,\mathcal{O}^*)
\ = \ 
\Gamma(\rq{X},\mathcal{O}^*)^H
\ = \ 
\KK^*.
$$

In order to verify the assertion on the dimension,
it suffices to show that the generic fiber of 
$p_X \colon \rq{X} \to X$ is a single $H$-orbit 
of dimension $\dim(H)$.
For this, note the generic 
orbit cone equals $Q^0(\gamma)$.
In particular, the generic orbits are 
of full dimension.
Property~\ref{bringdef}~(iii) of a bunched ring
guarantees $Q^0(\gamma) \in \rel(\Phi)$.
Lemma~\ref{bunch22max} then says that the generic 
orbits are closed in $\rq{X}$
und thus they show up as fibers of
$p_X \colon \rq{X} \to X$.

Finally, we have to show that each 
$\b{X}_{\gamma_0} \subseteq \rq{X}$,
where $\gamma_0 \in \rel(\Phi)$,
is saturated with respect to the quotient
map $p_X \colon \rq{X} \to X$.
For this, it suffices to show that 
given any  $\gamma_0 \in \rel(\Phi)$
and $H \mal z \subseteq \b{X}_{\gamma_0}$, 
which is closed in $\b{X}_{\gamma_0}$, 
then $H \mal z$ is already closed in $\rq{X}$.
Since $H \mal z \subseteq \rq{X}$ 
holds, there must be a face 
$\omega_0 \preceq \omega_H(z)$ with 
$\omega_0 \in \rel(\Phi)$,
see Lemma~\ref{bunch22max}.
Thus, we have 
\begin{eqnarray*}
\omega_0^\circ \ \cap \ Q^0(\gamma_0)^\circ 
& \ne & 
\emptyset.
\end{eqnarray*}
Moreover, denoting by $\b{H \mal z}$ the closure of 
$H \mal z$ in $\b{X}$,
Proposition~\ref{orbitcones} tells us 
that the orbit $H \mal z_0 \subseteq \b{H \mal z}$
corresponding to $\omega_0 \preceq \omega_H(z)$
is the unique closed orbit in 
$\b{H \mal z} \cap
\b{X}_{f_1^{u_1} \cdots f_r^{u_r}}$.
This implies $H \mal z_0 = H \mal z$ and hence
$\omega_H(z) \in \rel(\Phi)$.
Lemma~\ref{bunch22max} thus ensures
that $H \mal z$ is closed in $\rq{X}$.
\end{proof}

The next step is to provide an 
explicit isomorphism 
$K \to \Cl(X)$ for the variety
$X = X(R,\mathfrak{F},\Phi)$.
We need the following preparation.

\begin{construction}
\label{constr:divisorsDi}
Let $(R,\mathfrak{F},\Phi)$ be a bunched ring
with $(E \topto{Q} K,\gamma)$ as its projected 
cone, consider the action of $H = \Spec(\KK[K])$ 
on $\b{X} = \Spec(R)$ and set
$$
\rq{X} 
\ := \ 
\rq{X}(R,\mathfrak{F},\Phi),
\qquad\qquad
X 
\ := \ 
X(R,\mathfrak{F},\Phi).
$$
Then, for every member $f_i \in \mathfrak{F}$ 
its zero set $V(\rq{X},f_i)$ is an $H$-prime 
divisor and thus its image 
$D^i_X \subseteq X$ is a prime divisor.
Moreover, by~\ref{bringdef}~(iii), 
part~(b), we have 
$$ 
W_{\b{X}}
\ := \ 
\b{X} \setminus \bigcup_{i \ne j} V(\b{X}; f_i, f_j)
\ = \ 
\bigcup_{\gamma_0 \text{ facet of } \gamma} \b{X}_{\gamma_0}
\ \subseteq \ 
\rq{X},
$$
and hence the complement $\b{X} \setminus W_{\b{X}}$ 
is small in $\b{X}$.
By~\ref{bringdef}~(i), the group $H$ 
acts freely on $W_{\b{X}}$.
Moreover, as a union of 
saturated subsets $W_{\b{X}}$ is saturated under 
$p_X \colon \rq{X} \to X$.
Its image $W_{X} := p_X(W_{\b{X}})$ is open 
with small complement in $X$ and is given by 
\begin{eqnarray*}
W_X 
& =  &
X 
\setminus 
\bigcup_{i \ne j} D_X^i \cap D_X^j .
\end{eqnarray*}
\end{construction}

\begin{proposition}
\label{isodescr}
In the notation of~\ref{constr:divisorsDi}, 
cover $W_{\b{X}}$
by $H$-invariant open subsets $W_{j}$ 
admitting $w$-homogeneous functions
$h_{j} \in \Gamma(\mathcal{O}^{*},W_{j})$.
For every $w \in K$, fix a $w$-homogeneous 
function $h_{w} \in \KK(\rq{X})^{*}$.
\begin{enumerate}
\item
For any $w \in K$ there is a unique divisor
$D(h_w)$ on $X$ with  $p_X^*(D(h_w)) = \div(h_w)$;
on each $ p_X(W_j)$ it is given by 
$D(h_w) = \div(h_{w}/h_{j})$.
\item
Let $U \subseteq X$ be open.
Then, for every $w \in K$, there is an
isomorphism of $\KK$-vector spaces
$$ 
\Gamma(U,\mathcal{O}(D(h_w))) 
\ \to \  
\Gamma(p_X^{-1}(U), \mathcal{O})_w,
\qquad
g \ \mapsto \ p_X^*(g) h_w.
$$
\item
The assignment $w \mapsto D(h_w)$ induces an isomorphism 
from $K$ onto $\Cl(X)$, not depending on the choice of 
$h_{w}$:
$$
\b{D}_X \colon 
K \ \to \ \Cl(X),
\qquad
w \ \mapsto \ \b{D}(w) := [D(h_w)].
$$
\item
For $h_w = f_i \in \mathfrak{F}$, 
we obtain $D(f_i)= D_X^i$.
In particular, $\Cl(X)$ is generated by the 
classes $\b{D}(w_i)$ of the $D_X^i$,
where $1 \le i \le r$ and $w_i = \deg(f_i)$. 
\end{enumerate}
\end{proposition}

\begin{proof}
The $H$-invariant local equations $h_w/h_j$ 
define a Cartier divisor on $W_X = p_X(W_{\b{X}})$. 
Since $W_X$ has a small complement in $X$, this
Cartier divisor extends in a unique manner
to a Weil divisor $D(h_w)$ on $X$.
This shows~(i).

For~(ii), note that on each $U_j := p_X(W_j) \cap U$ the 
section $g$ is  given as $g = g_j'h_j/h_w$ 
with a regular function 
$g_j' \in \mathcal{O}(U_j)$.
Consequently, the function $p_X^*(g) h_w$ is 
regular on $W_{\b{X}} \cap p_X^{-1}(U)$, and thus, 
by normality, on $p_X^{-1}(U)$.
In particular, the assignment is a well defined
homomorphism.
Moreover, $f \mapsto f/h_w$ defines an inverse 
homomorphism.

We turn to~(iii). To see that the class 
$\b{D}(w)$ does not 
depend on the choice of $h_w$, consider a further 
$w$-homogeneous $g_w \in \KK(X)^*$. 
Then $f := h_w/g_w$ is an
invariant rational function descending 
to $X$, where we obtain
\begin{eqnarray*}
D(h_w) \ - \ D(g_w) &  = & \div(f).
\end{eqnarray*}
Thus, $\b{D} \colon K \to \Cl(X)$ is well defined 
and, by construction, it is homomorphic.
To verify injectivity, let $D(h_w) = \div(f)$ 
for some $f \in \KK(X)^*$. 
Then we obtain $\div(h_w) = \div(p_X^*(f))$.
Thus, $h_w/p_X^*(f)$ is an invertible 
function on $\b{X}$.
By $R^* = \KK^*$, it must be constant.
This implies $w = 0$.
For surjectivity, let any 
$D \in \WDiv(X)$ be given.
Then $p_X^*(D)$ is an invariant divisor
on $\b{X}$.
Since $R$ is factorially graded, we obtain 
$p_X^*(D) = \div(f)$ with some rational
function $h$ on $\b{X}$, wich
is homogeneous of some degree
$w$.
This means $D = D(h)$.

Finally, for~(iv), 
we only have to show that $D_{X}^{i}$ 
equals $D(f_i)$. 
By construction, these divisors 
have the same support, and we 
have 
$p_X^*(D(f_i)) = \div(f_i)$.
Since $\div(f_i)$ is $H$-prime 
and $D_X^i$ is prime, the assertion
follows.
\end{proof}

We indicate now, how to embed the variety $X$ 
arising from a given bunched ring 
$(R,\mathfrak{F},\Phi)$ into a toric variety.
Let $(E \topto{Q} K, \gamma)$ be the associated 
projected cone.
Then, with $M :=\ker(Q)$, we have mutually
dual exact sequences of rational vector spaces
$$ 
\xymatrix{
0 \ar[r]
&
L_\QQ \ar[r]
&
F_\QQ \ar[r]^{P}
&
N_\QQ \ar[r]
&
0
\\
0
&
K^0_\QQ \ar[l]
&
E_\QQ \ar[l]^{Q^0}
&
M_\QQ \ar[l]
&
0  \ar[l]
}
$$
Recall that $\gamma = \cone(e_1, \ldots, e_r)$ 
holds with the canonical basis $e_1, \ldots, e_r$ 
for the lattice $E = \ZZ^r$.
We denote the dual basis in $F = \Hom(E,\ZZ)$ 
again by $e_1, \ldots, e_r$.
Then $\delta := \cone(e_1, \ldots, e_r) \subseteq F_\QQ$ 
is the dual cone of $\gamma \subseteq E_\QQ$.

In order to relate bunches of cones to fans and vice 
versa, we need to transform collections of faces of 
$\gamma$ to collections of faces of $\delta$ and vice 
versa. This is done element-wise via the face 
correspondence
\begin{eqnarray*}
\faces(\delta) & \longleftrightarrow & \faces(\gamma)
\\
\delta_0 & \mapsto & \delta_0^* \ := \ \delta_0^\perp \cap \gamma
\\
\gamma_0^\perp \cap \delta =: \gamma_0^* & \mapsfrom & \gamma_0
\end{eqnarray*}

Given any two corresponding collections 
$\rq{\mathfrak{S}} \subseteq \faces (\delta)$ 
and $\rq{\Phi} \subseteq \faces(\gamma)$,
we infer from~\cite[Lemma~4.3]{BeHa0} that 
the following two compatibility statements 
are equivalent.
\begin{enumerate}
\item
Any two $\rq{\sigma}_1, \rq{\sigma}_2 \in \rq{\mathfrak{S}}$ 
admit an $L$-invariant separating linear form.
\item
For any two $\rq{\tau}_1, \rq{\tau}_2 \in \rq{\Phi}$ one has 
$Q^0(\rq{\tau}_1)^{\circ} \cap Q^0(\rq{\tau}_2)^{\circ} \ne \emptyset$.
\end{enumerate}

Given a compatible collection 
$\rq{\mathfrak{S}} \subseteq \faces(\delta)$, 
the images $P(\rq{\sigma})$, 
where $\rq{\sigma} \in \rq{\mathfrak{S}}$ 
generate a fan $\Sigma$ in $N_\QQ$,
which we call the {\em image fan\/} 
of $\rq{\mathfrak{S}}$.
Moreover, we call the fan $\rq{\Sigma}$ generated 
by a compatible $\rq{\mathfrak{S}}$ a 
{\em projectable fan\/} in $F_\QQ$.
Note that for a maximal compatible collection 
$\rq{\mathfrak{S}} \subseteq \faces(\delta)$, 
the image fan consists precisely 
of the images $P(\rq{\sigma})$, 
where $\rq{\sigma} \in \rq{\mathfrak{S}}$,
whereas $\rq{\Sigma} = \rq{\mathfrak{S}}$ 
need not hold in general.

\begin{construction}
\label{constr:neatemb}
Let $(R,\mathfrak{F},\Phi)$ be a bunched ring 
and $(E \topto{Q} K,\gamma)$ its projected cone,
consider the action of $H = \Spec(\KK[K])$ 
on $\b{X} = \Spec(R)$ and set
$$
\rq{X} \ := \ \rq{X}(R,\mathfrak{F},\Phi),
\qquad\qquad
X \ := \ X(R,\mathfrak{F},\Phi).
$$
Define a $K$-grading on the polynomial ring
$\KK[T_1, \ldots, T_r]$ by 
$\deg(T_i) := \deg(f_i)$.
This gives an $H$-equivariant closed 
embedding into $\b{Z} := \KK^r$,
namely
$$ 
\b{\imath} \colon 
\b{X} \ \to \ \b{Z},
\qquad
z \ \mapsto \ (f_1(z), \ldots, f_r(z)).
$$
Let $\rq{\Theta} \subseteq \faces(\gamma)$
be any compatible collection with 
$\rel(\Phi) \subseteq \rq{\Theta}$.
Then the restriction $\rq{\imath}$ 
of $\b{\imath}$ to $\rq{X}$ defines
a commutative diagram
$$ 
\xymatrix{
{\rq{X}} \ar[r]^{\rq{\imath}} \ar[d]_{\quot H}^{p_X}
&
{\rq{Z}} \ar[d]^{\quot H}_{p_Z}
\\
X  \ar[r]_{\imath}
&
Z
}
$$
with  the toric varieties $\rq{Z}$ and $Z$ 
defined by the projectable fan $\rq{\Sigma}$
associated to the collection $\rq{\mathfrak{S}}$
corresponding to $\rq{\Theta}$ and the image 
fan $\Sigma$ respectively.
\end{construction}

\begin{proposition}
\label{neatembex}
In the setting of Construction~\ref{constr:neatemb},
the following statements hold.
\begin{enumerate}
\item 
$p_Z \colon \rq{Z} \to Z$ 
is a toric Cox construction, and
we have $\rq{X} = \b{\imath}^{-1}(\rq{Z})$.
In particular,
$\rq{\imath} \colon \rq{X} \to \rq{Z}$
is an $H$-equivariant closed embedding.
\item
For $1 \le i \le r$, let $D^i_Z \subseteq Z$ 
be the invariant prime divisor corresponding to 
$e_i \in E$.
Then $D^i_X = \imath^*(D^i_Z)$ holds,
and, as a consequence, $W_X = \imath^{-1}(W_Z)$.
\item
The induced morphism $\imath \colon X \to Z$ of the 
quotient varieties is a neat closed embedding.
\end{enumerate}
\end{proposition}

\begin{proof}[Proof of Construction~\ref{constr:neatemb}
and Proposition~\ref{neatembex}]
In order to see that $p_Z \colon \rq{Z} \to Z$ 
is a toric Cox construction, we only have to 
check that the images $P(e_i) \in N$ are 
pairwise different primitive lattice vectors 
in $N$ and that the rays through them occur in 
the image fan $\Sigma$.
But this is guaranteed by 
$\rel(\Phi) \subseteq \rq{\Theta}$ 
and the Properties~\ref{bringdef}~(i) and~(iii),
part~(b) of $\Phi$.

By construction, $\rq{Z}$ is the union of all
subsets $\b{Z}_{\gamma_0} = \b{Z}_{T^u}$, 
where $\gamma_0 \in \rq{\Theta}$ and 
$u \in \gamma_0^\circ$. 
Thus $\rel(\Phi) \subseteq \rq{\Theta}$ ensures 
$\rq{X} \subseteq \b{\imath}^{-1}(\rq{Z})$.
To show equality, we may assume, 
by suitably enlarging $\rq{Z}$,
that $\rq{\Theta}$ is 2-maximal.
Then $\rq{Z} \subseteq \b{Z}$ 
is the associated ($H$,2)-maximal set.
To see that a given point $z \in \b{X} \cap \rq{Z}$
lies in $\b{\imath}(\rq{X})$, 
we may assume that $H \mal z$ is closed in 
$\rq{Z}$.
Then $\omega_H(z)$ belongs to $\rq{\Theta}$.
Moreover, $\omega_H(z)$ is a projected 
$\mathfrak{F}$-face, and thus, 
by 2-maximality, belongs to $\rel(\Phi)$.
This implies $z \in \b{\imath}(\rq{X})$.

Together, this proves Construction~\ref{constr:neatemb} 
and the first assertion of Proposition~\ref{neatembex}.
The second one follows from
the commutative diagram of 
Construction~\ref{constr:neatemb} 
and 
$$
p_X^* D_X^i 
\ = \ 
\div(f_i)
\ = \ 
\rq{\imath}^* \div(T_i)
\ = \ 
\rq{\imath}^*p_Z^* D_Z^i.
$$
To see that $\imath \colon X \to Z$ 
is a neat embedding, we just note
that the isomorphisms 
$\b{D}_X \colon K \to \Cl(X)$ 
and $\b{D}_Z \colon K \to \Cl(Z)$ of~Proposition~\ref{isodescr} 
fit into a commutative 
diagram
$$ 
\xymatrix{
& 
K 
\ar[dl]_\cong
\ar[dr]^\cong
& 
\\
{\Cl(X)}
& & 
{\Cl(Z)} 
\ar[ll]^{\imath^*}
}
$$
\end{proof}

\begin{corollary}
\label{veryneatemb}
If the toric ambient variety $Z$ of 
Construction~\ref{constr:neatemb}
is $\QQ$-factorial, then 
$p_Z^{-1}(X) \to X$ is a universal 
torsor for $X$ and the 
closure of $p_Z^{-1}(X)$ in 
$\b{Z}$ is the total coordinate space 
of $X$.
\end{corollary}

\begin{proof}[Proof of Theorem~\ref{thm:bring2var}, part~2]
Proposition~~\ref{isodescr} gives the isomorphism
of divisor class groups.
Taking $\rq{\Theta} = \rel(\Phi)$ 
in~\ref{constr:neatemb} gives the desired 
neat embedding $\imath \colon X \to Z$. 
In order to see that $p_X \colon \rq{X} \to X$
is a universal torsor, note that, by normality of 
$\rq{X}$ it suffices to prove this for 
the restriction $p_X \colon W_{\b{X}} \to W_X$.
This in turn follows Theorem~\ref{ringcalc} applied 
to any neat embedding $\imath \colon X \to Z$ as  
provided by Construction~\ref{constr:neatemb}.
Using once more normality, we obtain that
$R$ is the Cox ring of $X$.
\end{proof}

\section{Basic geometric properties}
\label{sec:geomprop}

In this section, we indicate how to read off
basic geometric properties of the 
variety associated to a bunched ring
from its defining data. 
We treat an explicit example, and 
at the end of the section, we characterize
the class of varieties that arise from bunched 
rings. 
Many statements are direct generalizations 
of those of~\cite{BeHa1}.
However our proofs are different, as we mostly 
don't use neat embeddings.

In the sequel,
$(R,\mathfrak{F},\Phi)$ 
is a bunched ring with projected
cone $(E \topto{Q} K, \gamma)$ and
admissible system of generators 
$\mathfrak{F} = \{f_1, \ldots, f_r\}$.
We consider the action of 
$H := \Spec(\KK[K])$ 
on $\b{X} := \Spec(R)$ and set
$$
\rq{X} 
\ := \ 
\rq{X}(R,\mathfrak{F},\Phi),
\qquad\qquad
X 
\ := \ 
X(R,\mathfrak{F},\Phi).
$$
Recall from Proposition~\ref{isodescr}
that we associated to any $h_w \in Q(R)^*_w$
a divisor $D(h_w)$ and this defines
a canonical isomorphism
$\b{D}_X \colon K \to \Cl(X)$.
In particular, we may canonically identify 
$K_\QQ^0 \cong \QQ \otimes K$ with the 
rational divisor class group $\Cl_\QQ(X)$.

In our first statement, we describe some cones 
of divisors. Recall that a divisor on a variety
is called {\em movable\/} if it has a positive multiple
with base locus of codimension at least two 
and it is called {\em semiample\/} if it has a base 
point free multiple.

\begin{proposition}
\label{divcones}
The cones of effective, movable, 
semiample and ample divisor classes 
of $X$ in $\Cl_\QQ(X) = K^0_\QQ$ 
are given as
$$
\Eff(X) 
\ = \
Q^0(\gamma),
\qquad
\Mov(X) 
\ = \  
\bigcap_{\gamma_0 \text{ facet of } \gamma} Q^0(\gamma_0),
$$
$$
\SAmple(X) 
\ = \
\bigcap_{\tau \in \Phi} \tau,
\qquad
\Ample(X) 
\ = \
\bigcap_{\tau \in \Phi} \tau^\circ.
$$
\end{proposition}

\begin{lemma}
\label{lem:zeroes}
In the setting of Proposition~\ref{isodescr},
consider the divisor $D(h_w)$ on $X$, 
a section $g \in \Gamma(X,\mathcal{O}(D(h_w)))$,
the sets 
$$ 
Z(g) 
\ := \
\Supp(\div(g) + D(h_w)),
\qquad\qquad
X_g \ := \ X \setminus Z(g),
$$
and let 
$f := p_X^*(g)h_w \in \Gamma(\b{X},\mathcal{O})_w$
denote the homogeneous function corresponding to $g$.
Then we have
$$ 
Z(g) \ = \ p_X(V(\rq{X};f)),
\qquad\qquad
p_X^{-1}(Z(g) \cap W_X) 
\ = \
V(\b{X}; f) \cap W_{\b{X}}.
$$
Moreover, if the open subset $X_g \subseteq X$ 
is affine, then its inverse image is 
given as $\pi_X^{-1}(X_g) = \b{X}_f$.
\end{lemma}

\begin{proof}
For the first equation, by surjectivity of 
$p_X$, it suffices to show $p_X^{-1}(X_g) = W$, 
where 
\begin{eqnarray*}
W
& := &  
\rq{X} \setminus p_X^{-1}(p_X( V(\rq{X},f))).
\end{eqnarray*}
Consider
$g^{-1} \in \Gamma(X_g,\mathcal{O}(D(h_w^{-1})))$.
The corresponding function 
$f' = p_X^*(g^{-1})h_w^{-1}$ 
on $p_X^{-1}(X_g)$ satisfies $f'f=1$,
which implies $p_X^{-1}(X_g) \subseteq W$.
Now consider 
$f^{-1} \in \Gamma(W,\mathcal{O})_{w}$.
The corresponding section $g' = f^{-1}h_w$ 
in $\Gamma(\pi(W),\mathcal{O}(D(h_w^{-1}))$
satisfies $g'g = 1$, which implies
$\pi(W) \subseteq X_g$.

The second equation follows directly from the first
one and the fact that on $W_{\b{X}} = p_X^{-1}(W_X)$
the $p_X$-fibers are precisely the $H$-orbits and 
hence $V(\b{X};f) \cap W_{\b{X}}$ is saturated with
respect to $p_X$.

Finally, let $X_g$ be affine. The first equation gives 
us $p_X^{-1}(X_g) \subseteq \b{X}_{f}$.
Moreover, $p_X^{-1}(X_g)$ is affine, and thus 
its boundary $\b{X} \setminus p_X^{-1}(X_g)$ 
is of pure codimension one. 
Inside $W_{\b{X}}$, this boundary coincides 
with $V(\b{X}; f)$.
Since $\b{X} \setminus W_{\b{X}}$ is small,
the assertion follows.
\end{proof}

\begin{proof}[Proof of Proposition~\ref{divcones}]
We have to characterize in terms of $w_0 \in K^0_\QQ$,
when a divisor $D(h_w) \in \WDiv(X)$ 
as constructed in Proposition~\ref{isodescr}
is effective, movable, semiample or ample.
The description of the effective cone $\Eff(X)$ 
is clear by Proposition~\ref{isodescr}~(ii).

For the description of the moving cone $\Mov(X)$, 
we may replace $X$ with $W_X$.
By Proposition~\ref{isodescr}~(ii)
and Lemma~\ref{lem:zeroes}, the base locus of 
$D(h_w)$ has a component of codimension 
one if and only if there is an $f_i \in \mathfrak{F}$ 
dividing all $f \in R_w$.
Thus, $D(h_w)$ has small stable base locus if and only if 
we have 
\begin{eqnarray*}
w^0
& \in & 
\bigcap_{i = 1}^r 
\cone(w_1, \ldots, w_{i-1},w_{i+1}, \ldots, w_r).
\end{eqnarray*}

For the description of the semiample cone, 
note that by  Lemma~\ref{lem:zeroes}
the divisor $D(h_w)$ is semiample if and only if 
$\rq{X}$ is contained in the union of 
all $\b{X}_f$, where $f \in R_{nw}$ and $n>0$. 
This in turn is equivalent to 
$$ 
w^0 
\ \in \ 
\bigcap_{z \in \rq{X}} \omega(z)
\ = \ 
\bigcap_{\tau \in \Phi} \tau.
$$ 

We come to the ample cone.
Suppose first that 
$w^0 \in \tau^\circ$ holds for all 
$\tau \in \Phi$.
Then, for some $n>0$, we may write
$\rq{X}$ as the union of $p_X$-saturated 
open subsets $\b{X}_f$ with functions 
$f \in R_{nw}$.
For the corresponding sections
$g \in \Gamma(X;\mathcal{O}(D(h^n_w)))$, 
Lemma~\ref{lem:zeroes} gives
$p_X(\b{X}_f) = X_g$.
In particular, all $X_g$ are affine 
and thus $D(h_w)$ is ample.

Now suppose that $D(h_w)$ is ample.
We may assume that $D(h_w)$ is even
very ample.
Then, for every 
$f^u = f_1^{u_1} \cdots f_r^{u_r}$ 
with $Q(u) = w$ the corresponding 
section $g^u \in \Gamma(X;\mathcal{O}(D(h_w)))$  
defines an affine set $X_{g^u} \subseteq X$,
and these sets cover $X$.
By  Lemma~\ref{lem:zeroes},
we obtain saturated inclusions 
$\b{X}_{f^u} \subseteq \rq{X}$ 
and hence a saturated inclusion 
$\rq{X}(w) \subseteq \rq{X}$,
see Proposition~\ref{H2quots}
for the notation.
By ($H$,2)-maximality of $\rq{X}(w)$,
we obtain $\rq{X}(w) = \rq{X}$
and hence $w \in \tau^\circ$ 
for every $\tau \in \Phi$.
\end{proof}

\begin{corollary}
\label{amplechambers}
If we vary the $\mathfrak{F}$-bunch 
$\Phi$, then the ample cones of the 
quasiprojective ones among the 
resulting varieties 
$X = X(R,\mathfrak{F},\Phi)$ 
are precisely the cones 
\begin{eqnarray*}
\lambda^\circ
& \subseteq & 
\bigcap_{\gamma_0 \text{ facet of } \gamma} Q^0(\gamma_0)^\circ
\end{eqnarray*}
with $\lambda \in \Lambda(\b{X})$.
If $\SAmple(X) \preceq \SAmple(X')$ 
holds, then $\rq{X}' \subseteq \rq{X}$ 
induces a projective morphism $X' \to X$,
which is an isomorphism in codimension one. 
\end{corollary}

In order to investigate local properties of 
$X$, we first observe that $X$ comes by
construction with a decomposition into
locally closed pieces.

\begin{construction}
\label{decomposconstr}
To any $\mathfrak{F}$-face $\gamma_0 \preceq \gamma$,
we associate a locally closed subset, namely
$$ 
\b{X}(\gamma_0)
\ := \ 
\{z \in \b{X}; \; f_i(z) \ne 0 \Leftrightarrow e_i \in \gamma_0\}
\ \subseteq \ 
\b{X}.
$$
These sets are pairwise disjoint 
and cover the whole $\b{X}$.
Taking the pieces defined by relevant 
$\mathfrak{F}$-faces, one obtains
a constructible subset
$$
\t{X} 
\ := \
\bigcup_{\gamma_{0} \in \rel(\Phi)} \b{X}(\gamma_{0})
\ \subseteq \
\rq{X},
$$
which is precisely the union of all closed 
$H$-orbits of $\rq{X}$.
The images of the pieces inside $\t{X}$ form
a decomposition of $X$ 
into pairwise disjoint locally closed 
pieces:
$$ 
X 
\ = \ 
\bigcup_{\gamma_{0} \in \rel(\Phi)} X(\gamma_{0}),
\qquad
\text{where }
X(\gamma_0)
\ := \ 
p_X(\b{X}(\gamma_0)).
$$
\end{construction}

\begin{example}
\label{decompostoric}
If we have $R = \KK[T_1, \ldots, T_r]$
and $\mathfrak{F} = \{T_1, \ldots, T_r\}$,
then $X$ is the toric variety arising from 
the image fan $\Sigma$ associated to 
$\rel(\Phi)$, 
and for any $\gamma_0 \in \rel(\Phi)$, 
the piece $X(\gamma_0) \subseteq X$ 
is precisely the toric orbit corresponding 
to the cone $P(\gamma_0^*) \in \Sigma$.
\end{example}

\begin{proposition}
\label{piecedecompos}
For any $\gamma_0 \in \rel(\Phi)$, 
the associated piece $X(\gamma_0)$
of the decomposition~\ref{decomposconstr}
has the following descriptions.
\begin{enumerate}
\item 
For any neat embedding $X \subseteq Z$ 
as constructed  in~\ref{constr:neatemb}, 
the piece $X(\gamma_0)$ is the intersection 
of $X$ with the toric orbit of $Z$ 
corresponding to $\gamma_0$.
\item
In terms of the prime divisors $D_X^i \subseteq X$ 
defined by the generators $f_i \in \mathfrak{F}$, 
the piece $X(\gamma_0)$ is given as 
\begin{eqnarray*}
X(\gamma_0)
& = &
\bigcap_{e_i \not\in \gamma_0} D^i_X
\ \setminus \ 
\bigcup_{e_j \in \gamma_0} D^j_X
\end{eqnarray*}
\item 
In terms of the open subsets $X_{\gamma_i} \subseteq X$ 
and  $\b{X}_{\gamma_i} \subseteq \b{X}$ defined by the 
relevant faces $\gamma_i \in \rel(\Phi)$, 
we have 
\begin{eqnarray*}
X(\gamma_0)
&  = & 
X_{\gamma_{0}} 
\ \setminus \ 
\bigcup_{\gamma_{0} \prec \gamma_{1} \in \rel(\Phi)} 
X_{\gamma_{1}},
\\
p_X^{-1}(X(\gamma_0))
&  = &
\b{X}_{\gamma_{0}} 
\ \setminus \ 
\bigcup_{\gamma_{0} \prec \gamma_{1} \in \rel(\Phi)} 
\b{X}_{\gamma_{1}}.
\end{eqnarray*}
\end{enumerate}
\end{proposition}

\begin{proof}[Proof of Construction~\ref{decomposconstr} 
and Proposition~\ref{piecedecompos}]
Obviously, $\b{X}$ is the union of the locally 
closed $\b{X}(\gamma_0)$, where $\gamma_0 \preceq \gamma$ 
runs through the $\mathfrak{F}$-faces.
To proceed, recall from Lemma~\ref{bunch22max}
and Proposition~\ref{H2quots} that 
$\rq{X} \subseteq \b{X}$ 
is the ($H$,2)-maximal subset given by the 
2-maximal collection
$$ 
\Psi
\ = \
\{Q^0(\gamma_0); \; \gamma_0 \in \rel(\Phi)\}
\ = \ 
\{\omega(z); \; H \mal z \text{ closed in } \rq{X}\}.
$$

Given $z \in \t{X}$, 
we have $\omega_H(z) = Q^0(\gamma_0) \in \Psi$.
Thus, Lemma~\ref{bunch22max} tells us that $H \mal z$ 
is closed in $\rq{X}$.
Conversely, if $H \mal z$ is closed in $\rq{X}$, 
consider the $\mathfrak{F}$-face
$$ 
\gamma_0 
\ := \ 
\cone(e_i; \; f_i(z) \ne 0)
\ \preceq \ 
\gamma.
$$
Then we have $z \in \b{X}(\gamma_0)$ and 
$Q^0(\gamma_0) = \omega_H(z) \in \Psi$;
see Lemma~\ref{bunch22max}. 
Hence $\gamma_0$ is a relevant face.
This implies $z \in \t{X}$.

All further statements are most easily seen 
by means of a neat embedding $X \subseteq Z$ as 
constructed in~\ref{constr:neatemb}. 
Let $\b{X} \subseteq \b{Z} = \KK^r$ denote the closed 
$H$-equivariant embedding arising from $\mathfrak{F}$.
Then, $\b{X}$ intersects precisely the 
$\b{Z}(\gamma_0)$, where $\gamma_0$ is an
$\mathfrak{F}$-face, and in these cases 
we have
\begin{eqnarray*}
\b{X}(\gamma_0) 
& = & 
\b{Z}(\gamma_0) \cap \b{X}.
\end{eqnarray*}
As mentioned in 
Example~\ref{decompostoric},
the images $Z(\gamma_0) = p_Z(\b{Z}(\gamma_0))$, 
where $\gamma_0 \in \rel(\Phi)$, 
are precisely the toric orbits of $Z$.
Moreover, we have
$$ 
\rq{X} \ = \ \rq{Z} \cap \b{X},
\qquad\qquad
\t{X} \ = \ \t{Z} \cap \b{X}.
$$
Since $p_Z$ separates $H$-orbits along $\t{Z}$,
we obtain $X(\gamma_0) = Z(\gamma_0) \cap X$
for every $\gamma_0 \in \rel(\Phi)$.
Consequently, the  $X(\gamma_0)$, where 
$\gamma_0 \in \rel(\Phi)$, are pairwise disjoint 
and form a decomposition 
of $X$ into locally closed pieces.
Finally, using $X(\gamma_0) = Z(\gamma_0) \cap X$
and $D^i_X = \imath^*(D^i_Z)$,
we obtain Assertions~\ref{piecedecompos}~(ii)
and~(iii) directly from the 
corresponding representations of the 
toric orbit $Z(\gamma_0)$.
\end{proof}

We say that a divisor class $[D] \in \Cl(X)$ 
is Cartier at a point $x \in X$ if some representative
$D$ is principal near $x$.
Similarly, we say that $[D] \in \Cl(X)$ is
$\QQ$-Cartier at $x$ if some nonnegative multiple
of a representative $D$ is principal near~$x$.

\begin{proposition}
\label{divisors}
Consider a relevant face 
$\gamma_{0} \in \rel(\Phi)$,
a point $x \in X(\gamma_{0})$,
and let $w \in K$.
\begin{enumerate}
\item 
The class $\b{D}_{X}(w)$ 
is Cartier at $x$ if and only if 
$w \in Q(\lin(\gamma_{0}) \cap E)$ 
holds.
\item 
The class $\b{D}_{X}(w)$ is 
$\QQ$-Cartier at $x$ if and 
only if 
$w^0 \in \lin(Q^0(\gamma_{0}))$ 
holds.
\end{enumerate}
\end{proposition}

\begin{lemma}
\label{lem:homfkt}
Consider a relevant face 
$\gamma_{0} \in \rel(\Phi)$,
a point $z \in \b{X}(\gamma_{0})$.
Then, for any $w \in K$,
the following statements 
are equivalent.
\begin{enumerate}
\item
There is an invertible homogeneous function 
$f \in \Gamma(H \mal z,\mathcal{O})_w$.
\item
One has $w \in Q(\lin(\gamma_{0}) \cap E)$.
\end{enumerate}
\end{lemma}

\begin{proof}
We have to determine the set 
$K_z \subseteq K$ of degrees $w \in K$ 
admitting an invertible homogeneous 
function on $H \mal z$.
Recall that the orbit $H \mal z$ is closed 
in $\b{X}_{\gamma_0}$, and thus the 
homogeneous functions on 
$H \mal z$ are precisely the 
restrictions of the homogeneous functions
on $\b{X}_{\gamma_0}$.
Moreover, we have $z \in \b{X}(\gamma_0)$,
and thus $f_i(z) \ne 0$ if and only if 
$e_i \in \gamma_0$.
Thus, the invertible homogenous functions
on $H \mal z$ are generated by the products
$f_1^{u_i} \cdots f_r^{u_r}$,
where $u \in \lin((\gamma_0) \cap E)$.
Consequently, we obtain
$K_z = Q(\lin((\gamma_0) \cap E)) \subseteq K$.
\end{proof}

\begin{proof}[Proof of Proposition~\ref{divisors}]
The divisor $D := D(h_w)$ is Cartier at $x$ 
if and only if there is a neighbourhood 
$U \subseteq X$ of $x$ and an invertible
section $g \in \Gamma(U,\mathcal{O}(D))$.
By Proposition~\ref{isodescr}, the latter
is equivalent to the existence of 
an invertible homogeneous function 
$f \in \Gamma(H \mal z,\mathcal{O})_w$
on the closed orbit 
$H \mal z \subseteq p_X^{-1}(x)$.
Now, Lemma~\ref{lem:homfkt} gives the desired
equivalence, and the first statement is proven.
The second is an immediate consequence
of the first.
\end{proof}

\begin{corollary}
\label{picarddescr}
Inside the divisor class group $\Cl(X) \cong K$, 
the Picard group of $X$ is given by
\begin{eqnarray*}
\Pic(X) 
& = & 
\bigcap_{\gamma_{0} \in \cov(\Phi)} 
              Q(\lin(\gamma_{0}) \cap E).
\end{eqnarray*}
\end{corollary}

\begin{proof}
A divisor $\b{D}(w)$ is Cartier if and only 
if it is Cartier along any $X(\gamma_0)$, 
where $\gamma_0 \in \rel(\Phi)$.
Since we have $\cov(\Phi) \subseteq \rel(\Phi)$,
and for any $\gamma_0 \in \rel(\Phi)$,
there is a $\gamma_1 \in \cov(\Phi)$
with $\gamma_1 \preceq \gamma_0$, 
it suffices to take the intersection over
$\cov(\Phi)$.
\end{proof}

A point $x \in X$ is factorial ($\QQ$-factorial) 
if and only if every Weil divisor is Cartier 
($\QQ$-Cartier) at $x$. 
Thus, Proposition~\ref{divisors} has the following 
application to singularities.

\begin{corollary}
\label{singularities}
Consider a relevant face $\gamma_{0} \in \rel(\Phi)$
and point $x \in X(\gamma_{0})$.
\begin{enumerate}
\item 
The point $x$ is factorial if and only if 
$Q$ maps $\lin(\gamma_{0}) \cap E$ onto $K$.
\item 
The point $x$ is $\QQ$-factorial if and only if 
 $Q^0(\gamma_{0})$ is of full dimension.
\end{enumerate}
\end{corollary}

\begin{corollary}
The variety $X$ is $\QQ$-factorial if and only 
if $\Phi$ consists of cones of full dimension.
\end{corollary}

Whereas local factoriality admits a simple 
combinatorial characterization,
smoothness is difficult in general.
Nevertheless, we have the following 
statement.

\begin{proposition}
\label{smoothchar}
Suppose that $\rq{X}$ is smooth,
let $\gamma_{0} \in \rel(\Phi)$,
and $x \in X(\gamma_{0})$.
Then $x$ is a smooth point if and only if
$Q$ maps $\lin(\gamma_{0}) \cap E$ onto $K$.
\end{proposition}

\begin{proof}
The ``only if'' part is clear by
Corollary~\ref{singularities}.
Conversely, if $Q$ maps $\lin(\gamma_{0}) \cap E$ 
onto $K$, then Lemma~\ref{lem:homfkt} 
says that the fibre $p_X^{-1}(x)$ 
consists of a single free $H$-orbit.
Consequently, $H$ acts freely over an open 
neighbourhood $U \subseteq X$ of $x$.
Thus $x$ is smooth.
\end{proof}

\begin{corollary}
\label{singulemb}
Let $X \subseteq Z$ be a neat embedding into 
a toric variety $Z$ as constructed in~\ref{constr:neatemb},
and let $x \in X$.
\begin{enumerate}
\item
The point $x$  is a factorial 
($\QQ$-factorial) point of $X$ 
if and only if it is a smooth 
($\QQ$-factorial) point of $Z$.
\item
If $\rq{X}$ is smooth, then $x$  
is a smooth point of $X$ 
if and only if it is a smooth 
 point of $Z$.
\end{enumerate}
\end{corollary}

As an immediate consequence of
general results on quotient singularities, 
see~\cite{Bo} and~\cite{HR}, 
one obtains the following.

\begin{proposition}
Suppose that $\rq{X}$ is smooth. 
Then $X$ has at most rational 
singularities. 
In particular, $X$ is Cohen-Macaulay.
\end{proposition}

In the case of a complete intersection 
$\b{X} \subseteq \b{Z}$,
we obtain a simple description of 
the canonical divisor class;
the proof given in~\cite[Theorem~9.1]{BeHa1}
works without changes, and therefore is 
omitted.

\begin{proposition}\label{candivdescr}
Suppose that the relations of $\mathfrak{F}$ 
are generated by $K$-homogeneous polynomials 
$g_{1}, \ldots, g_{d} \in \KK[T_{1}, \ldots, T_{r}]$,
where $d := r - \rank(K) - \dim(X)$.
Then, in $K \cong \Cl(X)$,
the canonical divisor class of $X$ is given by
\begin{eqnarray*}
D^{c}_{X}
& = &
\sum_{j=1}^{d} \deg(g_{j})
- 
\sum_{i=1}^{r} \deg(f_{i}).
\end{eqnarray*}
\end{proposition}

A variety is called ($\QQ$-)Gorenstein if 
(some multiple of) 
its anticanonical divisor is Cartier. 
Moreover, it is called ($\QQ$-)Fano
if (some multiple of) its anticanonical 
class is an ample Cartier divisor.

\begin{corollary}
In the setting of Proposition~\ref{candivdescr}, 
the following statements hold.
\begin{enumerate}
\item 
$X$ is $\QQ$-Gorenstein if and only if
\begin{eqnarray*} 
\sum_{i=1}^{r }\deg(f_i) \ - \ \sum_{j=1}^{d} \deg(g_{j})
& \in &
\bigcap_{\tau \in \Phi} \lin(\tau),
\end{eqnarray*}
\item 
$X$ is Gorenstein if and only if
\begin{eqnarray*}
\sum_{i=1}^{r }\deg(f_i) \ - \ \sum_{j=1}^{d} \deg(g_{j})
& \in &
\bigcap_{\gamma_{0} \in \cov(\Phi)} 
              Q(\lin(\gamma_{0}) \cap E).
\end{eqnarray*}
\item  
$X$ is $\QQ$-Fano if and only if we have
\begin{eqnarray*}
\sum_{i=1}^{r }\deg(f_i) \ - \ \sum_{j=1}^{d} \deg(g_{j})
& \in &
\bigcap_{\tau \in \Phi} \tau^{\circ},
\end{eqnarray*}
\item  
$X$ is Fano if and only if we have
\begin{eqnarray*}
\sum_{i=1}^{r }\deg(f_i) \ - \ \sum_{j=1}^{d} \deg(g_{j})
& \in &
\bigcap_{\tau \in \Phi} \tau^{\circ}
\ \cap \ 
\bigcap_{\gamma_{0} \in \cov(\Phi)} 
              Q(\lin(\gamma_{0}) \cap E).
\end{eqnarray*}
\end{enumerate}
\end{corollary}

\begin{example}
\label{ex:delpezzo}
Set $K := \ZZ^2$ and consider the $K$-grading of 
$\KK[T_1, \ldots, T_5]$ defined by $\deg(T_i) := \deg(f_i)$, 
where $w_i$ is the $i$-th column of the matrix
\begin{eqnarray*}
Q
& := & 
\left[
\begin{array}{rrrrr}
1 & -1 & 0 & -1 & 1
\\
1 & 1 & 1 & 0 & 2
\end{array}
\right]
\end{eqnarray*}
Then this $K$-grading descends to a $K$-grading of the
following factorial residue algebra
\begin{eqnarray*}
R 
& := & 
\KK[T_1, \ldots, T_5] \, / \, \bangle{T_1T_2 + T_3^2 + T_4T_5}.
\end{eqnarray*}
The classes $f_i \in R$ of $T_i \in \KK[T_1, \ldots, T_5]$,
where $1 \le i \le 5$, form a system $\mathfrak{F} \subset R$ 
of pairwise nonassociated homogeneous prime generators 
of $R$ and
$$ 
\Phi \ := \ \{\tau\},
\qquad
\tau \ := \ \cone(w_2,w_5)
$$
defines an $\mathfrak{F}$-bunch with associated 
projected cone $(E \topto{Q} K_0,\gamma)$, where 
$E = \ZZ^5$, and $\gamma = \cone(e_1, \ldots, e_5)$.
The covering collection $\cov(\Phi)$ consists of
$$ 
\cone(e_1,e_4),
\quad
\cone(e_2,e_5),
\quad
\cone(e_1,e_2,e_3),
\quad
\cone(e_3,e_4,e_5).
$$

The corresponding variety $X = X(R,\mathfrak{F},\Phi)$
is a $\QQ$-factorial projective surface. 
It has a single singularity, namely the point in the 
piece $X(\gamma_0)$ for $\gamma_0 = \cone(e_2,e_5)$.
The Picard group of $X$ is of index $3$ in $\Cl(X)$,
the canonical class of $X$ is Cartier and $X$ is 
Fano.
\end{example}

\begin{example}
Let $K := \ZZ \oplus \ZZ/3\ZZ$, and consider the 
$K$-grading of the polynomial ring 
$\KK[T_1, \ldots, T_6]$ given by
$$ 
\deg(T_1) \ = \ (1,\b{1}),
\quad
\deg(T_2) \ = \ (1,\b{2}),
\quad
\deg(T_3) \ = \ (1,\b{1}),
$$
$$
\deg(T_4) \ = \ (1,\b{2}),
\quad
\deg(T_5) \ = \ (1,\b{1}),
\quad
\deg(T_6) \ = \ (1,\b{2}).
$$
Similarly as in the preceding example,
this $K$-grading descends to a $K$-grading of 
the factorial residue algebra
\begin{eqnarray*}
R 
& := & 
\KK[T_1, \ldots, T_6] \, / \, \bangle{T_1T_2 + T_3T_4 + T_5T_6}.
\end{eqnarray*}
The classes $f_i \in R$ of $T_i \in \KK[T_1, \ldots, T_6]$,
where $1 \le i \le 5$, form a system 
$\mathfrak{F} \subset R$ 
of pairwise nonassociated homogeneous prime generators 
of $R$ and $\Phi \ := \ \{\QQ_{\ge 0}\}$ defines an 
$\mathfrak{F}$-bunch with associated projected cone 
$(E \topto{Q} K_0,\gamma)$ with $E = \ZZ^6$ and
$K_0 = \ZZ$, where $Q$ sends $e_i$ and we have
$\gamma = \cone(e_1, \ldots, e_6)$.

The resulting variety $X = X(R,\mathfrak{F},\Phi)$
is $\QQ$-factorial, projective, of dimension four
and has divisor class group $\ZZ \oplus \ZZ/3\ZZ$.
Its singular locus is of dimension two; it comprises,
for example, the surface $X(\gamma_0)$ for 
$\gamma_0 = \cone(e_1,e_3,e_5)$. 
The Picard group of $X$ is free cyclic and sits 
in the divisor class group of $X$ as 
\begin{eqnarray*}
(3 \mal \ZZ, \b{0}) & \subseteq & \ZZ \oplus  \ZZ/3\ZZ.
\end{eqnarray*}
\end{example}

Finally, we figure out the class of varieties that
arise from bunched rings.
We say that a variety $X$ is {\em $A_2$-maximal\/} if it 
is $A_2$ and admits no open embedding 
$X \subsetneq X'$ into an $A_2$-variety $X'$
such that $X' \setminus X$ is small in $X'$. 
For example, every projective variety is 
$A_2$-maximal.

\begin{theorem}
\label{thm:var2bring}
Let $X$ be a normal variety with 
$\Gamma(X,\mathcal{O}^*) = \KK^*$ and
finitely generated divisor class group
$\Cl(X)$ and finitely generated Cox 
ring $\mathcal{R}(X)$.
Then the following statements are equivalent.
\begin{enumerate}
\item
$X \cong X(R,\mathfrak{F},\Phi)$ 
holds with some bunched ring $(R,\mathfrak{F},\Phi)$.
\item
$X$ is $A_2$-maximal and every $\Cl(X)^0$-homogeneous 
$f \in \mathcal{R}(X)^*$ is constant.
\end{enumerate}
\end{theorem}

\begin{proof}
Suppose that $X$ is as in~(ii).
Set $K := \Cl(X)$ and 
$R :=  \mathcal{R}(X)$.
Then any $f \in R^*$
is necessarily $K^0$-homogeneous
and hence by assumption constant.
Moreover, by Proposition~\ref{gencoxprops},
the Cox ring $R$ is factorially 
$K$-graded.
Consider the total coordinate space 
$\b{X} = \Spec(R)$ with its action of 
$H = \Spec(\KK[K])$.
Then the universal torsor $p_X \colon \rq{X} \to X$ 
is a good quotient for the action of $H$.

Since $X$ is $A_2$-maximal, 
the subset $\rq{X} \subseteq \b{X}$ must 
be ($H$,2)-maximal. Thus $\rq{X}$ 
arises from a 2-maximal collection $\Psi$ of orbit
cones.
Choose any system $\mathfrak{F}$ of $R$ 
of pairwise nonassociated $K$-prime 
generators.
By Proposition~\ref{gencoxprops},
we have a small complement $\b{X} \setminus p_X^{-1}(X')$,
where $X' \subseteq X$ denotes the set of smooth 
points, and $H$ acts freely on $p_X^{-1}(X')$.
This ensures that $\mathfrak{F}$ is admissible.
Denoting by $\Phi$ the collection of 
minimal members of $\Psi$, we obtain 
the desired bunched ring $(R,\mathfrak{F}, \Phi)$.

Now suppose that~(i) holds, i.e.,
let $X = X(R,\mathfrak{F}, \Phi)$.
By construction, $X$ is an $A_2$-variety with 
finitely generated Cox ring.
Thus, we only have to show that $X$ is 
$A_2$-maximal. For this, suppose we have 
an open embedding $X \subseteq X'$ with 
small boundary into an $A_2$-variety. 
Replacing, if necessary, $X'$ with its
normalization, we may assume that $X'$ 
is normal.
Then $X$ and $X'$ share the same Cox ring
$R$ and they occur as good quotients 
of open subsets 
$\rq{X} \subseteq \rq{X}'$ of their common 
total
coordinate space. By $H$-maximality of 
$\rq{X}$, we obtain $\rq{X} = \rq{X}'$
and thus $X = X'$.
\end{proof}

\begin{corollary}
Let $X$ be a normal $A_2$-maximal variety with 
$\Gamma(X,\mathcal{O}) = \KK$ and finitely 
generated Cox ring. Then $X \cong X(R,\mathfrak{F},\Phi)$ 
holds with some bunched ring $(R,\mathfrak{F},\Phi)$.
\end{corollary}

\begin{corollary}
Let $X$ be a normal projective variety with finitely 
generated Cox ring. 
Then $X \cong X(R,\mathfrak{F},\Phi)$ 
holds with some 
bunched ring $(R,\mathfrak{F},\Phi)$.
\end{corollary}

\section{Ambient modification}
\label{sec:ambmod}

Here, we investigate the behaviour of Cox rings
under modifications induced from toric
modifcations.
As a preparation, we first 
consider modifications $Z_1 \to Z_0$ of 
toric varieties, their lifting properties 
to the respective Cox constructions, and 
give an explicit description of the situation
in terms of coordinates.

Let $\Sigma_0$ be a 
fan in a lattice $N$,
and let $v_1, \ldots, v_r \in N$ 
denote the primitive generators 
of the rays of $\Sigma_0$.
Suppose that $v_1, \ldots, v_r$ generate $N_\QQ$ as a 
vector space and, for some $2 \le d \le r$, we have
$$ 
\sigma_0 
\ := \
\cone(v_1, \ldots, v_d)
\ \in \
\Sigma_0.
$$

\begin{definition}
\label{def:index}
Let $v_\infty \in \sigma_0^\circ$ be a primitive 
lattice vector. 
We define the {\em index\/} of $v_\infty$ in 
$\sigma_0$ to be the minimal number 
$m_\infty \in \ZZ_{\ge 1}$ such that
there are nonnegative integers $a_1, \ldots, a_d$ 
with
\begin{eqnarray*}
m_\infty v_\infty 
& = & 
a_1v_1 + \cdots + a_rv_d.
\end{eqnarray*}
\end{definition}

The star of $\sigma_0$ 
in $\Sigma_0$ is defined as 
$\Star_{\Sigma_0}(\sigma_0) 
:= 
\{\sigma \in \Sigma_0; \; \sigma_0 \preceq \sigma\}$,
and the stellar subdivision
of $\Sigma_0$ at a primitive lattice vector 
$v_\infty \in \sigma_0^\circ$ is the fan
\begin{eqnarray*}
\Sigma_1
& := & 
\left(
\Sigma_0 
\ \setminus \ 
\Star_{\Sigma_0}(\sigma_0)
\right)
\ \cup \
\{\tau + \cone(v_\infty); \; 
\tau \precneqq \sigma \in \Star_{\Sigma_0}(\sigma_0)\}.
\end{eqnarray*}

We will consider the Cox constructions 
$P_0 \colon \rq{\Sigma}_0 \to \Sigma_0$ 
and 
$P_1 \colon \rq{\Sigma}_1 \to \Sigma_1$.
Recall that the fans $\rq{\Sigma}_0$ 
and $\rq{\Sigma}_1$ live in the lattices
$$
F_0 \ := \ \bigoplus_{i=1}^r \ZZ e_i,
\qquad
F_1 \ := \ \bigoplus_{i=1}^r \ZZ e_i \; \oplus \; \ZZ e_\infty
$$ 
and consist of certain faces of
the positive orthants $\delta_0 \subseteq (F_0)_\QQ$
and $\delta_1 \subseteq (F_1)_\QQ$. 
More precisely, as lattice homomorphisms,
the projection maps are given by
$$
P_0 \colon F_0 \to N, \quad e_i \mapsto v_i,
\qquad \qquad
P_1 \colon F_1 \to N, \quad e_i \mapsto v_i, \ e_\infty \mapsto v_\infty;
$$ 
note that in order to have $e_\infty \mapsto v_\infty$
in the Cox construction of $\Sigma_1$, 
it is necessary that $v_\infty \in N$ is a primitive 
lattice vector. 
The fans $\rq{\Sigma}_0$ 
and $\rq{\Sigma}_1$ are given by 
$$ 
\rq{\Sigma}_0
\ = \ 
\{\delta_0' \preceq \delta_0; \; P_0(\delta_0') 
\subseteq \sigma \in \Sigma_0\},
\qquad
\rq{\Sigma}_1
\ = \ 
\{\delta_1' \preceq \delta_1; \; P_1(\delta_1') 
\subseteq \sigma \in \Sigma_1\}.
$$

Now, fix nonnegative integers $a_1, \ldots, a_d$ 
as in Definition~\ref{def:index},
and define lattice homomorphisms
$$
G \colon F_1 \to F_0, 
\qquad
e_i \mapsto e_i,
\quad
e_\infty \mapsto a_1e_1 + \cdots + a_de_d,
$$
$$
G_1 \colon F_1 \to F_1, 
\qquad
e_i \mapsto e_i,
\quad
e_\infty \mapsto m_\infty e_\infty.
$$
Then all the lattice homomorphisms defined so far 
fit into the following commutative diagram
$$ 
\xymatrix{
&
F_1 \ar[dl]_{G_1} \ar[dr]^{G}
& 
\\
F_1 
\ar[d]_{P_1}
& &
F_0 
\ar[d]^{P_0}
\\
N 
\ar[rr]_{\id}
& &
N
}
$$
Note that $G$ defines a map of 
the fans $\rq{\Sigma}_1$ and 
$\rq{\Sigma}_0$ and, as well, 
of the respective fans of faces 
$\b{\Sigma}_1$ and $\b{\Sigma}_0$ 
of the positive orthants 
$\delta_1$ and $\delta_0$.

\begin{proposition}
\label{prop:liftdiag}
Let $Z_i, \rq{Z}_i$ and $\b{Z}_i$ be the 
toric varieties associated to the fans 
$\Sigma_i, \rq{\Sigma}_i$ and $\b{\Sigma}_i$
respectively.
Then the maps of the fans just
defined give rise to a commutative 
diagram of toric morphisms
$$ 
\xymatrix{
&
{\b{Z}_1}  \ar[dr]^{\b{\pi}} \ar[dl]_{\b{\pi}_1}
& 
\\
{\b{Z}_1}
&
{\rq{Z}_1} \ar[u] \ar[dl]_{\rq{\pi}_1} \ar[dr]^{\rq{\pi}}
& 
{\b{Z}_0}
\\
{\rq{Z}_1} \ar[d]_{p_1} \ar[u]
&
& 
{\rq{Z}_0} \ar[u] \ar[d]^{p_0}
\\
Z_1 
\ar[rr]_{\pi}
& &
Z_0
}
$$
where 
$\pi \colon Z_1 \to Z_0$ properly contracts an 
invariant prime divisor,
$p_i \colon \rq{Z}_i \to Z_i$ are universal 
torsors, 
$\b{\pi} \colon \b{Z}_1 \to \b{Z}_0$ is the quotient
for a $\KK^*$-action
and 
$\b{\pi}_1 \colon \b{Z}_1 \to \b{Z}_1$ is the quotient 
for an action of the group $C_{m_\infty}$ of 
$m_\infty$-th roots of the unity.
\end{proposition}

As announced, we will also look at this in terms 
of coordinates. 
For this, recall that, given a point $y$ of any 
affine $\KK^*$-variety $Y$, one says that 
\begin{itemize}
\item
the limit $\lim_{t \to 0} t \mal y$ exists,
if the morphism $\KK^* \to Y$, $t \mapsto t \mal y$
admits a continuation to $\KK$,
\item
the limit $\lim_{t \to \infty} t \mal y$ exists,
if the morphism $\KK^* \to Y$, $t \mapsto t^{-1} \mal y$
admits a continuation to $\KK$.
\end{itemize}
Moreover, one defines the {\em plus\/} and {\em minus cells\/} 
of the affine $\KK^*$-variety $Y$ to be the closed subsets
$$ 
Y^+ \ := \ \{y \in Y; \; \lim_{t \to 0} t \mal y \text{ exists}\},
\qquad
Y^- \ := \ \{y \in Y; \; \lim_{t \to \infty} t \mal y \text{ exists}\}.
$$  
The fixed point set $Y^0 \subseteq Y$ is given by 
$Y^0 = Y^+ \cap Y^-$. 
We denote by $Y^s := Y \setminus (Y^+ \cup Y^-)$ the 
union of all closed $\KK^*$-orbits.

\begin{lemma}
\label{liftmapprops}
Consider the coordinates
$z_1, \ldots, z_r$ on $\b{Z}_0 = \KK^r$ 
corresponding to 
$e_1, \ldots, e_r \in F_0$
and 
$z_1, \ldots, z_r, z_\infty$
on 
$\b{Z}_1 = \KK^{r+1}$ 
corresponding to 
$e_1, \ldots, e_r, e_\infty \in F_1$.
\begin{enumerate}
\item
The $\KK^*$-action on $\b{Z}_1$ and its good quotient 
$\b{\pi} \colon \b{Z}_1 \to \b{Z}_0$ are given by
\begin{eqnarray*}
t \mal (z_1, \ldots, z_r, z_\infty)
& = & 
(t^{-a_1}z_0, \ldots, t^{-a_d}z_d, z_{d+1}, \ldots, z_r, tz_\infty),
\\
\b{\pi}(z_1, \ldots, z_r, z_\infty)
& =  &
(z_\infty^{a_1}z_1, \ldots, z_\infty^{a_d}z_d, z_{d+1}, \ldots, z_r).
\end{eqnarray*}
\item
Plus and minus cells and the fixed point
set of the $\KK^*$-action are given by
\begin{eqnarray*}
\b{Z}_1^{+}
&  = &
V(\b{Z}_1; z_1, \ldots, z_d),
\\
\b{Z}_1^{-}
 & = &
V(\b{Z}_1;z_\infty),
\\
\b{Z}_1^{0}
& = &
V(\b{Z}_1; z_1, \ldots, z_d, z_\infty).
\end{eqnarray*}
In particular, $\b{Z}_1^- = \b{p_1^{-1}(E)}$, 
where $E \subseteq Z_1$ is
the exceptional divisor of $\pi \colon Z_1 \to Z_0$.
\item
The images of $\b{Z}_1^{0}$ and $\b{Z}_1^-$
under $\b{\pi} \colon \b{Z}_1 \to \b{Z}_0$ 
both equal 
$V(\b{Z}_0; z_1, \ldots, z_d)$,
we have 
$\b{Z}_1^+ \subseteq \b{Z}_1 \setminus \rq{Z}_1$ 
and
\begin{eqnarray*}
\b{Z}_1^s \ \cap \ \rq{Z}_1
& = & 
\b{Z}_1^s \ \cap \ \b{\pi}^{-1}(\rq{Z}_0)
\end{eqnarray*}
\item
The $C_{m_\infty}$-action on $\b{Z}_1$ and its good quotient 
$\b{\pi}_1 \colon \b{Z}_1 \to \b{Z}_1$ are given by
\begin{eqnarray*}
\zeta \mal (z_1, \ldots, z_r, z_\infty)
& =  &
(z_1, \ldots, z_r,\zeta z_\infty),
\\
\b{\pi}(z_1, \ldots, z_r, z_\infty)
& =  &
(z_1, \ldots, z_r,z_\infty^{m_\infty}).
\end{eqnarray*}
\end{enumerate}
\end{lemma}

\begin{proof}
Everything is obvious except, maybe, the displayed 
equation of the third assertion.
To verify it, take any face $\delta_1' \preceq \delta_1$ 
such that the associated distinguished point $x_{\delta_1'}$
lies in $\b{Z}_1^s$.
Then we have $e_\infty \not\in \delta_1'$ and 
$\{e_1, \ldots, e_d\} \not\subseteq \delta_1'$. 
Using this, we see
\begin{eqnarray*}
x_{\delta_1'} \ \in \ \rq{Z}_1
& \iff & 
P_1(\delta_1') \ \subseteq \ \tau_1 \ \in \Sigma_1
\\
& \iff & 
P_1(\delta_1') \ \subseteq \ \tau_0 \ \in \Sigma_0
\\
& \iff & 
P_0 \circ G(\delta_1') \ \subseteq \ \tau_0 \ \in \Sigma_0
\\
& \iff & 
G(\delta_1') \ \in \rq{\Sigma}_0
\\
& \iff &
x_{\delta_1'} \ \in \ \b{\pi}^{-1}(\rq{Z}_0).
\end{eqnarray*}
\end{proof}

Now we turn to the effect of the toric
modifications $\pi \colon Z_1 \to Z_0$ as just studied
on subvarieties $X_1 \subseteq Z_1$ 
and $X_0 \subseteq Z_0$ with $\pi(X_1) = X_0$.
A first step is to formulate conditions
on the ambient modification such that it 
preserves neat embeddings.

Let $Z$ be an irreducible variety, $X \subseteq Z$
a closed irreducible subvariets, and suppose that
$X$ and $Z$ are both smooth in codimension one.
Given an irreducible hypersurface $E \subseteq Z$ 
such that $D := X \cap E$ 
is an irreducible hypersurface in $X$, we say that
{\em some local equation for $E \subseteq Z$ restricts to a local 
equation for $D \subseteq X$}, if there are an open subset
$W \subseteq Z$ and an $f \in \mathcal{O}(W)$ with 
$E \cap W = \div(f)$ and 
$W \cap D \ne \emptyset$ 
such that on $U := W \cap X$ we have 
$D \cap U = \div(h)$ for the restriction 
$h := f_{\vert U}$.

\begin{definition}
\label{def:ambmod}
Let $\pi \colon Z_1 \to Z_0$ be the toric morphism 
arising from a stellar subdivision 
of simplicial fans, $T_1$ the big torus of $Z_1$,
and $E = \b{T_1 \mal z} \subseteq Z_1$ 
the exceptional divisor.
Let $X_1 \subseteq Z_1$ and
$X_0 \subseteq Z_0$ be
irreducible subvarieties such that
$X_1$ is smooth in codimension one
and $\pi(X_1) = X_0$ holds.
We say that $\pi \colon Z_1 \to Z_0$ is a 
{\em neat ambient modification\/} for 
$X_1 \subseteq Z_1$ and $X_0 \subseteq Z_0$ if 
$D := X_1 \cap E$ is an irreducible hypersurface
in $X_1$ intersecting $T_1 \mal z$,
some local equation for $E \subseteq Z_1$  
restricts to a local equation for $D \subseteq X_1$, 
and $X_0 \cap \pi(E)$ is of codimension 
at least two in $X_0$.
\end{definition}

Note that for a neat ambient modification
$\pi \colon Z_1 \to Z_0$ the subvarieties 
$X_1 \subseteq Z_1$ and $X_0 \subseteq Z_0$ 
intersect the big torus orbits in $Z_1$ and $Z_0$ 
respectively.
Moreover, $X_0$ inherits the property of being 
smooth in codimension one from $X_1$.

\begin{proposition}
\label{neat2neat}
Let $\pi \colon Z_1 \to Z_0$ be a neat ambient
modification for irreducible subvarieties 
$X_1 \subseteq Z_1$ and $X_0 \subseteq Z_0$,
both being smooth in codimension one.
Then the following statements are equivalent.
\begin{enumerate}
\item
$X_0 \subseteq Z_0$ is neatly embedded.
\item
$X_1 \subseteq Z_1$ is neatly embedded.
\end{enumerate}
\end{proposition}

\begin{proof}
We denote by $T_0$ and $T_1$ the acting tori of the 
ambient toric varieties $Z_0$ and $Z_1$ respectively.

We prove the implication ``(i)$\Rightarrow$(ii)''.
In order to verify Property~\ref{def:neatemb}~(i)
for $X_1 \subseteq Z_1$, let 
$E_1 = \b{T_1 \mal z_1}$ be any invariant prime divisor
of $Z_1$.
If $E_1 = E$ holds, then the definition of an ambient 
modification guarantees that $D_1 = X_1 \cap E_1$ is as wanted. 

For $E_1 \ne E$, note first that 
$X_1$ intersects $T_1 \mal z_1$, because 
$X_0$ intersects the one-codimensional orbit
$T_0 \mal \pi(z_1)$.
Thus, $D_1 = X_1 \cap E_1$ intersects 
$T_1 \mal z_1$.
Since $Z_1$ is $\QQ$-factorial, $E_1$ is locally the 
zero set of a function.
Consequently, $D_1$ is of pure codimension one 
in $X_1$.
Moreover, 
$$ 
D_1 \setminus E 
\ = \
X_1 \cap E_1 \setminus E 
\ = \
\pi^{-1}(X_0 \cap \pi(E_1) \setminus \pi(E))
$$
is irreducible.
Since $ X_1 \cap E$ is irreducible 
and not contained in $E_1$, we see that 
there are no components of $D_1$ inside 
$E$. Thus, $D_1$ is irreducible.

We verify Property~\ref{def:neatemb}~(ii)
for $X_1 \subseteq Z_1$.
Let $\kappa \colon X_1 \to X_0$ be the 
restriction of $\pi \colon Z_1 \to Z_0$.
Then we have the canonical push forward 
homomorphisms
$$ 
\pi_* \colon \Cl(Z_1) \ \to \ \Cl(Z_0),
\qquad \qquad   
\kappa_* \colon \Cl(X_1) \ \to \ \Cl(X_0),
$$
sending the classes of the exceptional 
divisors $E \subseteq Z_1$ and 
$D = X_1 \cap E \subseteq X_1$ to zero.
Denoting by  
$\imath_i \colon X_i \to Z_i$ the 
embeddings, Property~\ref{def:neatemb}~(i)
allows us to work with the pullback maps
$\imath_i^* \colon \Cl(Z_i) \to \Cl(X_i)$
as in Remark~\ref{weilpull} and obtain
a commutative diagram with exact rows
$$ 
\xymatrix{
0 \ar[r]
&
{\ZZ \mal [E]} \ar[r] \ar[d]_{\imath_1^*}
&
{\Cl(Z_1)} \ar[r]^{\pi_*}  \ar[d]_{\imath_1^*}
&
{\Cl(Z_0)}  \ar[r] \ar[d]_{\imath_0^*}^{\cong}
&
0
\\
0 \ar[r]
&{\ZZ \mal [D]} \ar[r]
&
{\Cl(X_1)}  \ar[r]_{\kappa_*}
&
{\Cl(X_0)}  \ar[r]
&
0
}
$$
Note that $\ZZ \mal [E] \subseteq \Cl(Z_1)$
as well as $\ZZ \mal [D] \subseteq \Cl(X_1)$
are nontrivial and free.
Since some local equation $f \in \mathcal{O}(W)$ 
of $E \subseteq Z_1$
restricts to a local equation 
$h \in \mathcal{O}(U)$
of $D \subseteq X_1$,
we obtain 
$$ 
D \cap U
\ = \ 
\div(h)
\ = \ 
\imath_1^* (\div(f))
\ = \ 
(\imath_1^* (E)) \cap U
$$
on the open set $U = W \cap X_1$. 
Sine $E$ and $D$ are irreducible hypersurfaces, 
we can conclude $D = \imath_1^*(E)$.
Thus, $\imath_1^* \colon \ZZ[E] \to \ZZ[D]$ 
is surjective, and hence an isomorphism.
Then, by the Five Lemma, also
$\imath_1^* \colon \Cl(Z_1) \to \Cl(X_1)$ 
is an isomorphism.

We turn to the implication ``(ii)$\Rightarrow$(i)''. 
In order to verify 
Property~\ref{def:neatemb}~(i) for $X_0 \subseteq Z_0$, 
let $E_0 \subseteq Z_0$ be any invariant prime divisor.
Then $E_0 = \pi(E_1)$ holds with a unique 
invariant prime divisor $E_1 \subseteq Z_1$.
Moreover, we have 
\begin{eqnarray*}
X_0 \cap E_0 \setminus \pi(E)
& = & 
\pi (X_1 \cap E_1 \setminus E).
\end{eqnarray*}
Since $X_1 \cap E_1 \setminus E$ 
is irreducible, 
the same holds for $X_0 \cap E_0 \setminus \pi(E)$.
Thus, in order to see that 
$X_0  \cap E_0$ is irreducible,
we only have to show that it has no 
irreducible components inside $X_0 \cap \pi(E)$. 
Since $Z_0$ is $\QQ$-factorial,
some multiple $nE_0$ is locally principal.
Consequently, all components of $X_0 \cap E_0$ 
are hypersurfaces in $X_0$ and none of them
can be contained in the small set $X_0 \cap \pi(E)$.  

We establish Property~\ref{def:neatemb}~(ii) 
for $X_0 \subseteq Z_0$.
Having already verified~\ref{def:neatemb}~(i), 
we may use the pullback map 
$\imath_0^* \colon \Cl(Z_0) \to \Cl(X_0)$.
As before, let
$\kappa \colon X_1 \to X_0$ be the 
restriction of $\pi \colon Z_1 \to Z_0$.
Then we have a commutative diagram
$$ 
\xymatrix{
0 \ar[r]
&
{\ZZ \mal [E]} \ar[r] \ar[d]_{\imath_1^*}^{\cong}
&
{\Cl(Z_1)} \ar[r]^{\pi_*}  \ar[d]_{\imath_1^*}^{\cong}
&
{\Cl(Z_0)}  \ar[r] \ar[d]_{\imath_0^*} 
&
0
\\
0 \ar[r]
&{\ZZ \mal [D]} \ar[r]
&
{\Cl(X_1)}  \ar[r]_{\kappa_*}
&
{\Cl(X_0)}  \ar[r]
&
0
}
$$
Again, the rows of this diagram are exact
sequences.
Using, e.g., the Five Lemma, we see that 
$\imath_0^* \colon \Cl(Z_0) \to \Cl(X_0)$ is an isomorphism.
\end{proof}

\begin{corollary}
Let $Z_1,Z_0$ be $\QQ$-factorial toric varieties
with $\Gamma(Z_i, \mathcal{O}^*) = \KK^*$,
and let $\pi \colon Z_1 \to Z_0$ be a neat ambient 
modification 
for $X_1 \subseteq Z_1$ and $X_0 \subseteq Z_0$.
Denote by $p_i \colon \rq{Z}_i \to Z_i$ the toric
Cox constructions and set $\rq{X}_i := p_i^{-1}(X_i)$.
\begin{enumerate}
\item
If $X_0 \subseteq Z_0$ is neatly embedded,
$X_1$ is normal with $\Gamma(X_1, \mathcal{O}^*) = \KK^*$
and $Z_1$ is smooth, then $p_1 \colon \rq{X}_1 \to X_1$ is 
a universal torsor of $X_1$.
\item
If $X_1 \subseteq Z_1$ is neatly embedded 
$X_0$ is normal with $\Gamma(X_0, \mathcal{O}^*) = \KK^*$
and $Z_0$ is smooth, then $p_0 \colon \rq{X}_0 \to X_0$ is 
a universal torsor of $X_0$.
\end{enumerate}
\end{corollary}

\begin{proof}
The statements are a direct consequence of 
Proposition~\ref{neat2neat} and 
Corollary~\ref{cor:smooth2tosor}.
\end{proof}

We investigate the effect of a neat ambient 
modification on the total coordinate space
and the Cox ring.
Let $\pi \colon Z_1 \to Z_0$ 
be the toric morphism 
associated to a stellar subdivision of 
simplicial fans.
Then Proposition~\ref{prop:liftdiag} provides
a commutative diagram
$$ 
\xymatrix{
&
{\b{Z}_1}  \ar[dr]^{\b{\pi}} \ar[dl]_{\b{\pi}_1}
& 
\\
{\b{Z}_1}
&
{\rq{Z}_1} \ar[u] \ar[dl]_{\rq{\pi}_1} \ar[dr]^{\rq{\pi}}
& 
{\b{Z}_0}
\\
{\rq{Z}_1} \ar[d]_{p_1}^{/ H_1} \ar[u]
&
& 
{\rq{Z}_0} \ar[u] \ar[d]^{p_0}_{/ H_0}
\\
Z_1 
\ar[rr]_{\pi}
& &
Z_0
}
$$
where $p_i \colon \rq{Z}_i \to Z_i$ 
are universal torsors, 
$\b{\pi} \colon \b{Z}_1 \to \b{Z}_0$
is the good quotient for a $\KK^*$-action 
on $\b{Z}_1$,
and 
$\b{\pi}_1 \colon \b{Z}_1 \to \b{Z}_1$
is the quotient for the action of the group 
$C_{m_\infty}$ of $m_\infty$-th toots of unity,
where $m_\infty$ denotes the index of the 
stellar subdivision.

Now, suppose that $\pi \colon Z_1 \to Z_0$ is 
a neat ambient modification 
for two irreducible subvarieties $X_1 \subseteq Z_1$ and 
$X_0 \subseteq Z_0$.
By definition, we have $X_0 = \pi(X_1)$
and both, $X_1$ and $X_0$, are smooth in 
codimension one.
Consider the restriction $\kappa \colon X_1 \to X_0$ 
of $\pi \colon Z_1 \to Z_0$,
the inverse images $\rq{X}_i := p_i^{-1}(X_i)$,
their closures $\b{X}_i \subseteq \b{Z}_i$,
the inverse image
$\rq{Y}_1 := \b{\pi}_1^{-1}(\rq{X}_1)$, 
its closure $\b{Y}_1 \subseteq \b{Z}_1$ 
and the restrictions 
$\kappa$, $\rq{\kappa}$ and 
$\b{\kappa}$ 
of $\pi$, $\rq{\pi}$ and $\b{\pi}$ 
as well as 
$\rq{\kappa}_1$ and 
$\b{\kappa}_1$ 
of $\rq{\pi}_1$ and $\b{\pi}_1$ 
respectively.

\begin{lemma}
\label{lem:induceddiag}
In the above setting, $\b{Y}_1 \subseteq \b{Z}_1$ 
is $C_{m_\infty}$-invariant and $\KK^*$-invariant.
If $\b{T}_0 \subseteq \b{Z}_0$ denotes the 
big torus orbit, then one has
$$
\b{Y}_1
\ = \
\b{\pi}_1^{-1}(\b{X}_1)
\ = \  
\b{\b{\pi}^{-1}(\b{X}_0 \cap \b{T}_0)}.
$$
If one of the embeddings $X_i \subseteq Z_i$ 
is neat, then 
$\b{Y}_1$ is irreducible, and one has a 
commutative diagram 
$$
\xymatrix{
&
{\b{Y}_1}  
\ar[dr]_{\b{\kappa}}^{\quot \KK^*}
\ar[dl]^{\b{\kappa}_1}_{/ C_{m_\infty}}
& 
\\
{\b{X}_1}
&
{\rq{Y}_1} \ar[u] \ar[dl]^{\rq{\kappa}_1} \ar[dr]_{\rq{\kappa}}
& 
{\b{X}_0}
\\
{\rq{X}_1} \ar[d]^{p_1}_{/ H_1} \ar[u]
&
& 
{\rq{X}_0} \ar[u] \ar[d]_{p_0}^{/ H_0}
\\
X_1 
\ar[rr]_{\kappa}
& &
X_0
}
$$
\end{lemma}

\begin{proof}
Let $T_i \subseteq Z_i$ denote the respective 
big torus orbits.
The first statement follows from the fact
that $\b{Y}_1 \subseteq \b{Z}$ is the closure 
of the following $C_{m_\infty}$-invariant and 
$\KK^*$-invariant subset
$$ 
\b{\pi}_1^{-1}(p_1^{-1}(\pi^{-1}(X_0 \cap T_0)))
\ = \ 
\b{\pi}^{-1}(p_0^{-1}(X_0 \cap T_0))
\ \subseteq \
\b{Y}_1.
$$
If one of the embeddings $X_i \subseteq Z_i$
is neat, then Proposition~\ref{neat2neat} 
ensures that $X_0 \subseteq Z_0$ is neat.
Theorem~\ref{ringcalc} then guarantees
that $p_0^{-1}(X_0 \cap T_0)$ is irreducible.
This gives the second claim.
The commutative diagram is then obvious.
\end{proof}

This statement shows in particular
that $\b{Y}_1$ 
can be explicitly computed from either $\b{X}_1$,
or $\b{X}_0$ and vice versa.
Cutting down the minus cell $\b{Z}_1^-$ 
of the $\KK^*$-action on $\b{Z}_1$, 
we obtain the minus cell
for the induced $\KK^*$-action on $\b{Y}_1$:
\begin{eqnarray*} 
\b{Y}_1^- & = & \b{Y}_1 \ \cap \ \b{Z}_1^-.
\end{eqnarray*}
Since $\b{\pi}_1^{-1}(\b{Z}_1^-) = \b{Z}_1^-$ 
holds 
and $\b{\pi}_1^{-1}$ equals  the identical map
along $\b{Z}_1^-$, we see that 
the group $H_1$, acting on $\b{Z}_1$,
leaves $\b{Y}_1^-$ invariant.

\begin{definition}
\label{def:controlled}
We say that the neat ambient modification
$\pi \colon Z_1 \to Z_0$ 
for $X_1 \subseteq Z_1$ and $X_0 \subseteq Z_0$ 
is {\em controlled\/} if the minus
cell $\b{Y}_1^-$ is $H_1$-irreducible,
i.e., $H_1$ acts transitively on the collection 
of irreducible components of $\b{Y}_1^-$.
\end{definition}

Note that in the case of a torus $H_1$, being controlled 
just means that the minus cell $\b{Y}_1^-$ is irreducible.
For our first main result on modifications, we recall the 
following notation: 
given a toric variety $Z$ with toric prime divisors 
$E_1, \ldots, E_r \subseteq Z_0$ and Cox construction
$p \colon \rq{Z} \to Z$ and a closed subvariety 
$X \subseteq Z$, we set 
$$
W_{Z} \ := \ Z \setminus \bigcup_{i \ne j} E_i \cap E_j,
\qquad\qquad
W_{\b{Z}} \ := \ p^{-1}(W_{Z}),
$$ 
$$
W_{X} \ :=  W_{Z} \cap \rq{X}
\qquad\qquad
W_{\b{X}} \ := \ W_{\b{Z}} \cap \rq{X}.
$$

\begin{theorem}
\label{controlled2small}
Let $Z_1,Z_0$ be $\QQ$-factorial toric varieties 
with $\Gamma(Z_i,\mathcal{O}^*) = \KK^*$, 
and let $\pi \colon Z_1 \to Z_0$ be
a neat ambient modification for 
$X_1 \subseteq Z_1$ and $X_0 \subseteq Z_0$,
where $X_0 \subseteq Z_0$ is neatly embedded 
and $W_{X_1}$ is normal with 
$\Gamma(W_{X_1}, \mathcal{O}^*) = \KK^*$.
Suppose that $\b{X}_0 \setminus \rq{X}_0$
is of codimension at least two in $\b{X}_0$
and that $\pi \colon Z_1 \to Z_0$ is controlled.
\begin{enumerate}
\item
The complement $\b{X}_1  \setminus \rq{X}_1$
is of codimension at least two in $\b{X}_1$.
Moreover, if $\rq{X}_1$ is normal, 
then so is $X_1$ and $\rq{X}_1 \to X_1$ is 
a universal torsor of $X_1$.
\item
The normalization of $X_1$ has finitely generated 
Cox ring and it has the $H_1$-equivariant 
normalization of $\b{X}_1$ as its total coordinate 
space.
\end{enumerate}
\end{theorem}

Note that if the embedding $X_0 \subseteq Z_0$ 
is obtained as in Construction~\ref{constr:neatemb},
then any  neat controlled ambient modification
$\pi \colon Z_1 \to Z_0$ for normal subvarieties
$X_1 \subseteq Z_1$ and $X_0 \subseteq Z_0$
fullfills the assumption of Theorem~\ref{controlled2small}.
Before giving the proof, we note an immediate 
consequence.

\begin{corollary}
\label{cor:contr2small}
Let $Z_1,Z_0$ be $\QQ$-factorial toric varieties 
with $\Gamma(Z_i,\mathcal{O}^*) = \KK^*$
and $\pi \colon Z_1 \to Z_0$ 
a neat controlled ambient modification
for $X_1 \subseteq Z_1$ and $X_0 \subseteq Z_0$,
where $X_0 \subseteq Z_0$ is normal,
neatly embedded 
and $\b{X}_0 \setminus \rq{X}_0$ 
is of codimension at least two in $\b{X}_0$.
If $\b{Y}_1$ is normal, then $X_1$ is normal
with universal torsor $\rq{X}_1 \to X_1$ 
and it has $\b{X}_1$ together with the 
induced $H_1$-action 
as its total coordinate space.
\end{corollary}

\begin{proof}
The fact that $\b{Y}_1$ is normal implies that
$\b{X}_1 = \b{Y}_1 / C_{m_\infty}$ and 
$X_1 = \rq{X}_1 / H_1$ are normal.
Thus, we may apply Theorem~\ref{controlled2small}. 
\end{proof}

The proof of Theorem~\ref{controlled2small}
uses the following observation on 
$\KK^*$-varieties $Y$; recall, that we denote
the fixed point set by $Y^0$ and the set 
of closed one-dimensional orbits by $Y^s$.

\begin{lemma}
\label{lem:smallcompl}
Let $Y$ be an irreducible affine $\KK^*$-variety
with quotient $\pi \colon Y \to Y \quot \KK^*$
such that $\pi(Y^0)$ is small in $Y \quot \KK^*$
and $Y \setminus Y^-$ is affine.
Then, for any open invariant $V \subseteq Y$ 
with $Y^s \setminus V$ small in $Y$,
the following statements are equivalent.
\begin{enumerate}
\item 
The complement $Y \setminus V$ is small in $Y$.
\item
The complement $Y \setminus V$ comprises no
component of $Y^-$.
\end{enumerate}
\end{lemma}

\begin{proof}
Suppose that~(i) holds. Since $Y \setminus Y^-$
is affine, every component of the minus cell
$Y^-$ is a hypersurface.
In particular, the small complement $Y \setminus V$ 
cannot contain any component of $Y^-$.

Now suppose that~(ii) holds. The complement 
of $V$ in $Y$ can be decomposed as follows
into locally closed subsets:
\begin{eqnarray*}
Y \setminus V
& = &
Y^s \setminus V
\ \cup \
Y^- \setminus V
\ \cup \
Y^+ \setminus V.
\end{eqnarray*}
By assumption, the first two sets are small
in $Y$.
So, it suffices to show that $Y^+$ is small
in $Y$.
Since $Y \setminus Y^-$ is affine, the same holds
for $(Y \setminus Y^-)/\KK^*$.
Consequently, in the commutative diagram
$$ 
\xymatrix{
Y \setminus Y^- 
\ar@{}[r]|\subseteq
\ar[d]_{/ \KK^*}^{\pi_-}
&
Y  
\ar[d]^{\quot \KK^*}_{\pi}
\\
(Y \setminus Y^-)/\KK^*
\ar[r]
&
Y \quot \KK^*
}
$$
the induced (projective) morphism of 
quotient spaces is finite.
Thus, since $\pi(Y^0)$ is small in $Y \quot \KK^*$,
we see that $\pi_-(Y^+)$ is small in $(Y \setminus Y^-)/\KK^*$.
Consequently, $Y^+ \setminus Y^-$ must be small 
in $Y \setminus Y^-$. This gives the assertion.
\end{proof}

\begin{proof}[Proof of Theorem~\ref{controlled2small}]
Only the assertion that
$\b{X}_1 \setminus \rq{X}_1$ is small,
needs a proof;
the remaining statements are direct 
applications of 
Corollary~\ref{cor:qfact2tosor}.
In order to verify that
$\b{X}_1 \setminus \rq{X}_1$ is small,
it suffices to show that 
$\b{Y}_1 \setminus \rq{Y}_1$
is small.
The idea is to 
apply Lemma~\ref{lem:smallcompl}
to the $\KK^*$-variety $\b{Y}_1$.
Let us check the assumptions.

Since $X_0 \subseteq Z_0$ is neatly 
embedded and $p_0$ is equidimensional, 
the complement $\rq{X}_0 \setminus W_{\b{X}_0} 
= 
p_0^{-1}(X_0 \setminus W_{X_0})$ 
is small in $\rq{X}_0$.
Denote by $D \subseteq X_1$ the 
exceptional divisor.
Then Lemma~\ref{liftmapprops}~(iii) 
tells us that the image of the fixed point 
set under the quotient map
$\b{\kappa} \colon \b{Y}_1 \to \b{X}_0$
satisfies 
$$ 
\b{\kappa}(\b{Y}_1^{0}) \cap \rq{X}_0
\ \subseteq \ 
p_0^{-1}(\pi(D))
\ \subseteq \ 
\rq{X}_0 \setminus W_{\b{X}_0}.
$$
As observed before, the latter set is small,
and, moreover $\b{X}_0 \setminus \rq{X}_0$ 
is small by assumption.
Consequently $\b{\kappa}(\b{Y}_1^{0})$
is small in $\b{X}_0$. 
Since, in the notation of Lemma~\ref{liftmapprops},
we have 
$\b{Y}_1^- = V(\b{Y}_1,z_\infty)$, the set 
$\b{Y}_1 \setminus \b{Y}_1^-$ 
is affine.
Moreover, again by Lemma~\ref{liftmapprops}~(iii), 
we have
$$ 
\b{Y}_1^s \cap \rq{Y}_1
\ = \
\b{Y}_1 \cap \b{Z}_1^s \cap \rq{Z}_1
\ = \
\b{Y}_1 \cap \b{Z}_1^s \cap \b{\pi}^{-1}(\rq{Z}_0)
\ = \
\b{Y}_1^s \cap \b{\kappa}^{-1}(\rq{X}_0).
$$
Since $\kappa \colon \b{Y}_1^s \to \kappa(\b{Y}_1^s)$
has constant fibre dimension (one) and 
$\b{X}_0 \setminus \rq{X}_0$ is small,
we can conclude that 
$\b{Y}_1^s \setminus \rq{Y}_1$
is small.
Since $\b{Y}_1^-$ intersects $\rq{Y}_1$,
and the ambient modification 
$Z_1 \to Z_0$ is controlled,
we obtain that every component of  
$\b{Y}_1^-$ intersects $\rq{Y}_1$.
Thus, we may apply Lemma~\ref{lem:smallcompl},
and obtain that $\b{Y}_1 \setminus \rq{Y}_1$ 
is small.
\end{proof}

Our second main result on modifications treats the 
case of contractions, i.e., we start with a 
neatly embedded subvariety $X_1 \subseteq Z_1$.

\begin{theorem}
\label{controlled2small2}
Let $\pi \colon Z_1 \to Z_0$ the toric morphism
associated to a stellar subdivision of simplicial 
fans and suppose $\Gamma(Z_i,\mathcal{O}^*) = \KK^*$.
Let $X_1 \subseteq Z_1$ be neatly embedded 
with  $W_{X_1}$ normal and
$\Gamma(W_{X_1}, \mathcal{O}^*) = \KK^*$.
Set $X_0 := \pi(X_1)$ and suppose that 
$\b{X}_1 \setminus \rq{X}_1$ is of codimension 
at least two in $\b{X}_1$.
\begin{enumerate}
\item 
The embedding $X_0 \subseteq Z_0$ is neat and 
$\pi \colon Z_1 \to Z_0$ is a neat controlled ambient 
modification for 
$X_1 \subseteq Z_1$ and $X_0 \subseteq Z_0$.
\item
The complement $\b{X}_0  \setminus \rq{X}_0$
is of codimension at least two in $\b{X}_0$.
If  $\rq{X}_0$ is normal, then so is $X_0$ 
and $\rq{X}_0 \to X_0$ is 
a universal torsor of $X_0$.
\item
The normalization of $X_0$ has finitely generated 
Cox ring and it has 
the $H_0$-equivariant normalization of $\b{X}_0$
as its total coordinate space.
\end{enumerate}
\end{theorem}

Again we state a consequence before
proving the result; note that, below,  
if the neat embedding $X_1 \subseteq Z_1$ 
is obtained by Construction~\ref{constr:neatemb} 
and the index of the ambient modification
equals one, then all assumptions are satisfied.

\begin{corollary}
Let $Z_1,Z_0$ be $\QQ$-factorial toric varieties 
with $\Gamma(Z_i,\mathcal{O}^*) = \KK^*$, and let
$\pi \colon Z_1 \to Z_0$ be
a neat ambient modification
for $X_1 \subseteq Z_1$ and $X_0 \subseteq Z_0$,
where $X_1 \subseteq Z_1$ is normal,
neatly embedded 
and $\b{X}_1 \setminus \rq{X}_1$ 
is of codimension at least two in $\b{X}_1$.
Then, if $\b{Y}_1$ is normal, $X_0$ is normal
with universal torsor $\rq{X}_0 \to X_0$ 
and $\b{X}_0$ together with the induced $H_0$-action 
as its total coordinate space.
\end{corollary}

\begin{proof}
The fact that $\b{Y}_1$ is normal implies that
$\b{X}_0 = \b{Y}_1 / \KK^*$ and 
$X_0 = \rq{X}_0 / H_0$ are normal.
Thus, we may apply Theorem~\ref{controlled2small2}. 
\end{proof}

Similarly as for our first main result,
we need also for proving the second one
an observation on $\KK^*$-actions.

\begin{lemma}
\label{lem:kstar2}
Let $Y$ be an irreducible affine $\KK^*$-variety
admitting nontrivial homogeneous 
global regular functions of positive degree.
Then, for any irreducible component
$B \subseteq Y^0$, we have 
$B \subseteq  \b{Y^+ \setminus Y^0}$.
\end{lemma}

\begin{proof}
Since there are functions of positive degree,
we have $Y \ne Y^-$.
Thus, we have the following commutative 
diagram with the induced projective 
morphism $\varphi$ of quotient spaces
$$
\xymatrix{
Y \setminus Y^- 
\ar@{}[r]|\subseteq
\ar[d]_{/ \KK^*}^{\pi_-}
&
Y  
\ar[d]^{\quot \KK^*}_{\pi}
\\
(Y \setminus Y^-)/\KK^*
\ar[r]_{\quad \varphi}
&
Y \quot \KK^*
}
$$
Now, $\pi$ maps $Y^0$ bijectively 
onto $\pi(Y^0)$.
Hence $\pi(B)$ is an irreducible component 
of $\pi(Y^0)$.
Thus, there is an irreducible component
$A \subseteq \pi_-^{-1}(\varphi^{-1}(\pi(B)))$
with $\pi(A) = \pi(B)$.
For this component, we have 
$A \subseteq Y^+ \setminus Y^0$ and 
$B \subseteq \b{A}$.
\end{proof}

\begin{proof}[Proof of Theorem~\ref{controlled2small2}]
Let $E \subseteq Z_1$ be the exceptional 
divisor.
Since $X_1 \subseteq Z_1$ is neatly 
embedded, $D := X_1 \cap E$ is an 
irreducible hypersurface in $X_1$. 
Its inverse image with respect to
$p_1 \colon \rq{X}_1 \to X_1$
is the irreducible hypersurface
$$
p_1^{-1}(D)  
\ = \ 
\rq{X}_1 \cap \b{Z}_1^- 
\ = \ 
\b{\pi}_1(\rq{Y}_1 \cap \b{Y}_1^-)
\ \subseteq \
\rq{X}_1.  
$$
Since $\b{\pi}_1$ equals the identity along 
$\b{\pi}_1^{-1}(\b{Z}_1^-) = \b{Z}_1^-$,
the set $\rq{Y}_1 \cap \b{Y}_1^-$ is
$H_1$-irreducible.
By Lemma~\ref{liftmapprops}~(ii), the minus cell
$\b{Y}_1^-$ is the zero set of a function
and hence of pure codimension one.
Since $\b{X}_1 \setminus \rq{X}_1$ is small,
the same holds for  $\b{Y}_1 \setminus \rq{Y}_1$,
and we can conclude that also $\b{Y}_1^-$ is 
$H_1$-irreducible.

In order to obtain that $\pi \colon Z_1 \to Z_0$ 
is a neat in controlled ambient modification,
we still have to show that $D = X_1 \cap E$ 
has a small image $\pi(D)$ in $X_0$.
Using Lemma~\ref{liftmapprops}~(ii) 
we observe
$$
\pi(D)
\ = \
\pi(X_1 \cap E)
\ = \
\pi\left(p_1\left(\b{\pi}_1\left(  
\rq{Y}_1  \cap \b{Y}_1^- 
\right)\right)\right)
\ \subseteq \
p_0\left(\b{\pi}\left(\b{Y}_1^- \right) \cap \rq{X}_0 \right).
$$
Thus, we have to show that
$\b{\pi}(\b{Y}_1^-) = \b{\pi}(\b{Y}_1^0)$
is of codimension at least two 
in $\b{X}_0$.
For this, it suffices to show that 
$\b{Y}_1^0$ is of codimension 
at least three in $\b{X}_0$.
But this follows from the facts that,
according to Lemma~\ref{lem:kstar2},
any component of $\b{Y}_1^0$
is a proper subset of a component 
of $\b{Y}_1^+$, and that we have
$\b{Y}_1^+ \subseteq \b{Y}_1 \setminus \rq{Y}_1$,
where the latter set is small in $\b{Y}_1$.

Let us see, why $\b{X}_0 \setminus \rq{X}_0$ 
is small in $\b{X}_0$.
First note that we have 
$\b{\pi}(\rq{Y}_1) \subseteq \rq{X}_0$.
Thus, if $\b{X}_0 \setminus \rq{X}_0$ would 
contain a divisor, then also 
$\b{Y}_1 \setminus \rq{Y}_1$
and hence $\b{X}_1 \setminus \rq{X}_1$
must contain a divisor, a contradiction.
All remaining statements follow from
Corollary~\ref{cor:qfact2tosor}.
\end{proof}

\section{Combinatorial contraction}

We take a closer look at $\QQ$-factorial 
projective varieties with finitely 
generated Cox ring.
Corollaries~\ref{amplechambers} and~\ref{singularities} 
tell us that all those sharing a common 
Cox ring
have the same moving cone and their
semiample cones are the full dimensional 
cones of a fan subdivision of this common 
moving cone.
Moreover, any two of them differ only by a 
small birational transformation, i.e., 
a birational map, which is an isomorphism
in codimension one.

Conversely, if any two $\QQ$-factorial 
projective varieties with finitely generated 
Cox ring differ only by 
a small birational transformation
then their Cox rings are isomorphic,
and hence, up to isomorphism, both of them show 
up in a common picture as above. 
In this section, we link all these 
pictures to a larger one by studying 
certain divisorial contractions.
The following concepts will be crucial.

\begin{definition}
Let $X$ be $\QQ$-factorial projective variety 
with finitely generated Cox ring.
\begin{enumerate}
\item
We call a class $[D] \in \Cl(X)$ 
{\em combinatorially contractible\/}
if it generates an extremal ray 
of the effective cone of $X$ and, 
for some representative $D$, 
all vector spaces 
$\Gamma(X,\mathcal{O}(nD))$,
where $n > 0$, are of dimension one.
\item
We say that the variety $X$ is
{\em combinatorially minimal\/} 
if it has no combinatorially contractible
divisor classes.
\end{enumerate}
\end{definition}

As we will see, combinatorial minimality of 
$X$ is equivalent to saying that the 
moving cone of $X$ equals its effective 
cone. 
Moreover, a surface $X$ turns out to be 
combinatorially minimal if and only if 
its semiample cone equals its effective cone;
for example, $\PP_1 \times \PP_1$ is 
combinatorially minimal, whereas the first
Hirzebruch surface is not.
The following theorem shows that one may 
reduce any $X$ in a very controlled process to a 
combinatorially minimal one.

\begin{theorem}
\label{combcontr}
For every $\QQ$-factorial projective variety 
$X$ with finitely generated Cox ring arises 
from a combinatorially minimal one $X_0$
via a finite sequence
$$ 
\xymatrix{
X = X_n' \ar@{-->}[r]
&
X_{n}  \ar[r]
&
X'_{n-1}  \ar@{-->}[r]
&
X_{n-1}  \ar[r]
& 
\quad \ldots \quad
\ar[r]
& 
X_0' = X_0
}
$$
where $\xymatrix@1@!{X_i' \ar@{-->}[r] & X_{i}}$
is a small birational transformation and 
$\xymatrix@1@!{X_i \ar[r] & X_{i-1}'}$
comes from a neat controlled ambient modification
of $\QQ$-factorial projective toric varieties.
\end{theorem}

The crucial part of the proof of this statement 
is an explicit contraction criterion for 
varieties arising from bunched rings, see~
Proposition~\ref{combcontract}.
There, we also explicitly describe the contracting 
ambient modification, which in turn allows 
computation of the Cox rings in the contraction steps.

Let $(R,\mathfrak{F}, \Phi)$ be a 
bunched ring with
$\mathfrak{F} = \{f_1, \ldots, f_r\}$
and
$(E \topto{Q} K,\gamma)$
as its projected cone.
The diagonalizable group 
$H := \Spec(\KK[K])$ acts on 
$\b{X} := \Spec(R)$.
Defining a $K$-grading on 
$\KK[T_1, \ldots, T_r]$ by
$\deg(T_i) := w_i =: \deg(f_i)$, we obtain
a graded epimorphism
$$ 
\KK[T_1, \ldots, T_r] \ \to \ R,
\qquad
T_i \ \mapsto \ f_i.
$$
This gives rise to an $H$-equivariant closed 
embedding $\b{X} \to \b{Z}$, 
where $\b{Z} = \KK^r$. 
Note that the actions of $H$ on $\b{X}$
and $\b{Z}$ define two collections 
of orbit cones, namely $\Omega_H(\b{X})$ 
and $\Omega_H(\b{Z})$.
These collections in turn give rise to 
two GIT-fans, $\Lambda(\b{X})$ and 
$\Lambda(\b{Z})$.
Since 
$\Omega_H(\b{X}) \subseteq  \Omega_H(\b{Z})$
holds, 
$\Lambda(\b{Z})$ refines $\Lambda(\b{X})$.

Let us briefly recall how to construct 
neat embeddings in terms of $\Lambda(\b{Z})$.
The map $Q \colon E \to K$,
$e_i \mapsto w_i$ defines two mutually dual 
exact sequences of vector spaces 
$$ 
\xymatrix{
0 \ar[r]
&
L_\QQ \ar[r]
&
F_\QQ \ar[r]^{P}
&
N_\QQ \ar[r]
&
0
\\
0
&
K^0_\QQ \ar[l]
&
E_\QQ \ar[l]^{Q}
&
M_\QQ \ar[l]
&
0  \ar[l]
}
$$
As earlier, write $e_1, \ldots, e_r \in F$ for the 
dual basis, set $\delta = \gamma^\vee$ and 
$v_i := P(e_i) \in N$.
Note that dealing with projective $X$ and $Z$ amounts 
to requiring that the image $Q^0(\gamma) \subseteq K^0_\QQ$
is pointed.
There is an injection 
\begin{eqnarray*}
\Lambda(\b{Z})
& \longrightarrow &
\{\text{maximal compatible } 
\rq{\mathfrak{S}} 
\subseteq \faces(\delta)\}
\\
\eta 
& \mapsto & 
\{\gamma^*_0; \; \lambda^\circ \subseteq Q^0(\gamma_0)^\circ\}
\end{eqnarray*}
and the maximal projectable collections 
$\rq{\mathfrak{S}}$ assigned to
the $\eta \in \Lambda(\b{Z})$ are precisely those
having a polytopal image fan $\Sigma$ in $N$ with 
its rays among $\varrho_i := \QQ_{\ge 0} \mal v_i$, where 
$1 \le i \le r$.
In fact, $\Sigma$ is the normal fan of the fiber 
polytope $(Q^0)^{-1}(w^0) \cap \gamma$ for any 
$w^0 \in \eta^\circ$.
As a consequence of Construction~\ref{constr:neatemb}
and Proposition~\ref{neatembex}, we note
the following.

\begin{remark}
\label{assocemb}
Suppose that the variety $X$ associated 
to $(R,\mathfrak{F}, \Phi)$ is $\QQ$-factorial 
and projective. 
Then, for any $\eta \in \Lambda(\b{Z})$ with
$\eta^\circ \subseteq \Ample(X)$, we obtain 
a neat embedding $X \to Z$ into the toric 
variety $Z$ arising from the image fan of 
$\eta$. Moreover, $Z$ is $\QQ$-factorial and 
projective.
We call the neat embeddings $X \to Z$ constructed
this way {\em associated neat embeddings}. 
\end{remark}

We are ready to introduce the ``bunch theoretical''
formulation of combinatorial contractibility.
In the above setup, it is the following.

\begin{definition}
We say that a weight 
$w_i = \deg(f_i) \in K$ 
is {\em exceptional\/}, 
if $\QQ_{\ge 0} \mal w_i$ 
is an extremal ray of 
$Q^0(\gamma)$ and 
for any $j \ne i$ we have 
$w_j \not \in \QQ_{> 0} \mal w_i$.
\end{definition}

\begin{example}
Consider the two gradings of 
the polynomial ring 
$\KK[T_1,\ldots,T_4]$
assigning $T_i$ the degree 
$w_i \in \ZZ^2$ as follows;
the shaded area indicates~$Q^0(\gamma)$.
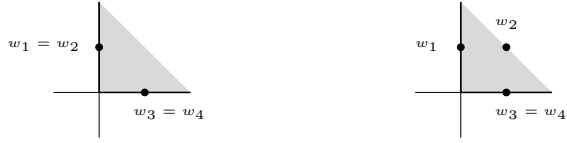
\begin{center}
\begin{picture}(0,0)%
\includegraphics{p1p1hirz.pstex}%
\end{picture}%
\setlength{\unitlength}{1243sp}%
\begingroup\makeatletter\ifx\SetFigFont\undefined%
\gdef\SetFigFont#1#2#3#4#5{%
  \reset@font\fontsize{#1}{#2pt}%
  \fontfamily{#3}\fontseries{#4}\fontshape{#5}%
  \selectfont}%
\fi\endgroup%
\begin{picture}(3648,2745)(3136,-2323)
\put(5626,-1861){\makebox(0,0)[lb]{\smash{{\SetFigFont{5}{6.0}{\familydefault}{\mddefault}{\updefault}{\color[rgb]{0,0,0}$w_3=w_4$}%
}}}}
\put(3151,-511){\makebox(0,0)[lb]{\smash{{\SetFigFont{5}{6.0}{\familydefault}{\mddefault}{\updefault}{\color[rgb]{0,0,0}$w_1=w_2$}%
}}}}
\end{picture}%

\qquad \qquad
\qquad \qquad
\begin{picture}(0,0)%
\includegraphics{p1p1hirz2.pstex}%
\end{picture}%
\setlength{\unitlength}{1243sp}%
\begingroup\makeatletter\ifx\SetFigFontNFSS\undefined%
\gdef\SetFigFontNFSS#1#2#3#4#5{%
  \reset@font\fontsize{#1}{#2pt}%
  \fontfamily{#3}\fontseries{#4}\fontshape{#5}%
  \selectfont}%
\fi\endgroup%
\begin{picture}(2748,2745)(4036,-2323)
\put(5626,-1861){\makebox(0,0)[lb]{\smash{{\SetFigFontNFSS{5}{6.0}{\familydefault}{\mddefault}{\updefault}{\color[rgb]{0,0,0}$w_3=w_4$}%
}}}}
\put(4051,-511){\makebox(0,0)[lb]{\smash{{\SetFigFontNFSS{5}{6.0}{\familydefault}{\mddefault}{\updefault}{\color[rgb]{0,0,0}$w_1$}%
}}}}
\put(5626,-61){\makebox(0,0)[lb]{\smash{{\SetFigFontNFSS{5}{6.0}{\familydefault}{\mddefault}{\updefault}{\color[rgb]{0,0,0}$w_2$}%
}}}}
\end{picture}%

\end{center}
Then in the first one $w_1$ is not 
exceptional, but in the second one it is.
Note that these configurations define
$\PP_1 \times \PP_1$ and the first 
Hirzebruch surface.
\end{example}

\begin{lemma}
\label{subdiv2chambers}
Consider a chamber $\eta_1 \in \Lambda(\b{Z})$, 
the corresponding maximal projectable fan
$\rq{\Sigma}_1$ and the associated 
image fan $\Sigma_1$. 
Suppose that 
$$ 
\dim(\eta_1) 
\ = \ 
\dim(K_\QQ),
\qquad\qquad
\eta_1^\circ 
\ \subseteq \ 
\bigcap_{j=1}^r \cone(e_l; \; l \ne j)^\circ
$$
holds, i.e., the fan $\Sigma_1$ is simplicial and all 
$\varrho_1, \ldots, \varrho_r$ belong to $\Sigma_1$,
and fix any index $i$ with $1 \le i \le r$. 
Then the following statements are equivalent.
\begin{enumerate}
\item
There is a polytopal simplicial fan $\Sigma_0$
such that $\Sigma_1$ is the stellar subdivision 
of $\Sigma_0$ at the vector $v_i$.
\item
The weight $w_i$ is exceptional, there is a 
chamber $\eta_0 \in \Lambda(\b{Z})$ of full dimension,
with $w^0_i \in \eta_0$ and 
$\eta_0 \cap \eta_1$ is a facet of both, 
$\eta_0$ and $\eta_1$.
\end{enumerate}
If one of these statements holds, 
then one may choose $\Sigma_0$ of~(i) 
to be the image fan of the maximal projectable fan 
associated to $\eta_0$ of~(ii),
and, moreover, one has
\begin{eqnarray*}
\eta_0 \cap \eta_1
& = & 
\eta_0 \cap \cone(w^0_1,\ldots,w^0_{i-1},w^0_{i+1}, \ldots, w^0_r).
\end{eqnarray*}
\end{lemma}

\begin{proof}
For the sake of a simple notation, we 
rename  the index $i$ in question to $\infty$
and renumber our vectors in $N$ to 
$v_1, \ldots, v_s, v_\infty$.

Suppose that~(i) holds. Then we obtain the following 
maximal projectable subfans of the fan of faces 
of the positive orthant 
$\delta = \cone(e_1, \ldots, e_s, e_\infty)$ 
in $F_\QQ$:
\begin{eqnarray*}
\rq{\Sigma}_0 
& := & 
\{\delta_0 \preceq \cone(e_1, \ldots, e_s); \; 
P(\delta_0)  \in \Sigma_0\},
\\
\rq{\Sigma}_1 
& := & 
\{\delta_0 \preceq \cone(e_1, \ldots, e_s, e_\infty); \; 
P(\delta_0)  \in \Sigma_1\},
\\
\rq{\Sigma}_{01} 
& := & 
\{\delta_0 \preceq \cone(e_1, \ldots, e_s, e_\infty); \; 
P(\delta_0) \subseteq \sigma \in \Sigma_0\}.
\end{eqnarray*}
Note that $\rq{\Sigma}_0$ and $\rq{\Sigma}_{01}$ have 
$\Sigma_0$ as their image fan, whereas $\rq{\Sigma}_1$ 
has $\Sigma_1$ as its image fan.
Let $\rq{\Theta}_0, \rq{\Theta}_1$ and $\rq{\Theta}_{01}$
denote the compatible collections of faces of 
$\gamma \subseteq E_\QQ$ corresponding to 
$\rq{\Sigma}_0, \rq{\Sigma}_1$ and $\rq{\Sigma}_{01}$. 
The associated chambers in $\Lambda(\b{Z})$ are given by 
$$ 
\eta_0 
\ = \ 
\bigcap_{\rq{\theta}_0 \in \rq{\Theta}_0} Q^0(\rq{\theta}_0),
\qquad
\eta_1 
\ = \ 
\bigcap_{\rq{\theta}_1 \in \rq{\Theta}_1} Q^0(\rq{\theta}_1),
\qquad
\eta_{01} 
\ = \ 
\bigcap_{\rq{\theta}_{01} \in \rq{\Theta}_{01}} Q^0(\rq{\theta}_{01}).
$$ 
Moreover, as $\rq{\Sigma}_0$ and $\rq{\Sigma}_1$ 
are subfans of $\rq{\Sigma}_{01}$, 
Proposition~\ref{H2quots} tells us
that $\eta_{01}$ is a face of both,
$\eta_0$ and $\eta_1$.

Since the fan $\Sigma_0$ is simplicial, the map $F$ is 
injective along the cones of $\rq{\Sigma}_0$, which
in turn implies that $\eta_0$ is of full dimension.
Moreover, since $\QQ_{\ge 0} \mal e_\infty \not\in \rq{\Sigma}_0$ 
holds, we obtain
$$ 
w^0_\infty \ \in \ \eta_0,
\qquad\qquad
\eta_0
\ \not \subseteq \
\cone(w^0_1, \ldots, w^0_s).
$$
In particular, we see that $w^0_\infty$ is exceptional.
The only thing, which remains to be verified is
that $\eta_{01}$ is of codimension one in $K^0_\QQ$.
For this, let $\sigma \in \Sigma_0$ be any cone
with $\varrho_\infty \subseteq \sigma$.
Then $\sigma$ is generated by some 
$v_{i_1}, \ldots, v_{i_d}$, where $1 \le i_d \le r$.
We set
$$ 
\rq{\sigma}_0 
\ := 
\ \cone(e_{i_1}, \ldots, e_{i_d}) 
\ \in \ 
\rq{\Sigma}_0,
\qquad
\rq{\sigma}_{01} 
\ := 
\ \cone(e_{i_1}, \ldots, e_{i_d},e_\infty) 
\ \in \ 
\rq{\Sigma}_{01},
$$
Using this notation, we obtain the following decompositions
of the chamber $\eta_{01}$ and its relative interior:
\begin{eqnarray*}
\eta_{01}
& = & 
\bigcap_{\varrho_\infty \not \subseteq \sigma \in \Sigma_0^{\max}} 
Q^0(\rq{\sigma}_{0}^*)
\ \cap \
\bigcap_{\varrho_\infty \subseteq \sigma \in \Sigma_0^{\max}} 
Q^0(\rq{\sigma}_{01}^*),
\\
\eta_{01}^\circ
& = & 
\bigcap_{\varrho_\infty \not \subseteq \sigma \in \Sigma_0^{\max}} 
Q^0(\rq{\sigma}_{0}^*)^\circ
\ \cap \
\bigcap_{\varrho_\infty \subseteq \sigma \in \Sigma_0^{\max}} 
Q^0(\rq{\sigma}_{01}^*)^\circ.
\end{eqnarray*}

For any $\sigma \in \Sigma_0$,
the cone $Q^0(\rq{\sigma}_0^*)$ 
is of full dimension
in $K^0_\QQ$.
For a cone $\sigma \in \Sigma_0$  with 
$\varrho_\infty \subseteq \sigma$, 
the cone $Q^0(\rq{\sigma}^*_{01})$ 
must be of codimension one in $K^0_\QQ$,
because it is not of full dimension
and generated by one ray less
than $Q^0(\rq{\sigma}_0^*)$.

Now, there is a unique minimal $\kappa \in \Sigma_0$ 
comprising $\varrho_\infty$.
For any further $\sigma \in \Sigma_0$ 
comprising $\varrho_\infty$, 
we have 
$Q^0(\rq{\sigma}^*_{01}) \subseteq Q^0(\rq{\kappa}^*_{01})$.
Thus, the above decompositions of $\eta_{01}$ and 
$\eta_{01}^\circ$ show that $\eta_{01}$
is of codimension one.

Suppose that~(ii) holds. 
Let $\rq{\Sigma}_0$ denote the 
maximal projectable fan 
corresponding to $\eta_0$,
and let $\Sigma_0$ be its image 
fan in $N$.
Note that $\Sigma_0$ is simplicial, 
since $\eta_0$ is of full dimension.

Since $w^0_\infty$ is exceptional
and belongs to $\eta_0$, we have 
$\eta_0 \not \subseteq \cone(w^0_1, \ldots, w^0_s)$.
Moreover,
$\QQ_{\ge 0} \mal e_\infty \in \rq{\Sigma}_1$ implies
$\eta_1 \subseteq \cone(w^0_1, \ldots, w^0_s)$.
Thus, we obtain 
\begin{eqnarray*}
\eta_0 \ \cap \ \eta_1
& \subseteq & 
\eta_0 \ \cap \ \cone(w^0_1, \ldots, w^0_s).
\end{eqnarray*}
By dimension reasons, this must even be an equality.
Moreover, we note that the hyperplane $H_{01}$ in $K^0_\QQ$ 
separating $\eta_0$ and $\eta_1$ does not contain
$w^0_\infty$.

We show now that any $\varrho_i$ with 
$1 \le i \le s$ must belong to $\Sigma_0$.
Suppose that some $\varrho_i$ does not.
As before, we obtain
\begin{eqnarray*}
\eta_0 
\ \cap \ 
\eta_1
& = & 
\eta_0 
\ \cap \ 
\cone(w^0_1, \ldots,w^0_{i-1},w^0_{i+1}, \ldots, w^0_s, w^0_\infty).
\end{eqnarray*}
The left hand side is contained in the 
separating hyperplane $H_{01}$, and the right
hand side contains $w^0_\infty$. A contradiction.

Now consider the stellar subdivision $\Sigma_2$ 
of $\Sigma_0$ in $\varrho_\infty$. 
Then $\Sigma_2$ is a polytopal simplicial fan.
Thus
\begin{eqnarray*}
\rq{\Sigma}_2 
& := & 
\{\delta_0 \preceq \cone(e_1, \ldots, e_s,e_\infty); \; 
P(\delta_0)  \in \Sigma_2\},
\end{eqnarray*}
arises from a full dimensional chamber 
$\eta_2 \in \Lambda(\b{Z})$.
As seen in the proof of ``(i)$\Rightarrow$(ii)'',
the chambers $\eta_0$ and $\eta_2$ share a common 
facet. This gives
$$
\eta_0 \cap \eta_2 
\ = \ 
\eta_0  \cap  \cone(w^0_1, \ldots, w^0_s)
\ = \ 
\eta_0 \cap \eta_1.
$$
Since the cones $\eta_0$, $\eta_1$ and 
$\eta_2$ all show up in a common fan, we may conclude 
$\eta_1 = \eta_2$.
This implies $\Sigma_1 = \Sigma_2$.
\end{proof}

\begin{proposition}
\label{combcontract}
Let $X_1$ be a $\QQ$-factorial 
projective variety arising 
from a bunched ring $(R,\mathfrak{F},\Phi)$, 
and let $D_1^i \subseteq X_1$ be the prime 
divisor given by $p_{X_1}^* D_1^i = \div(f_i)$ 
for the generator $f_i \in \mathfrak{F}$.
Then the following statements are equivalent.
\begin{enumerate}
\item
There are associated neat embeddings
$X_1 \subseteq Z_1$ and $X_0 \subseteq Z_0$
into $\QQ$-factorial projective toric varieties
and a neat controlled ambient modification 
$Z_1 \to Z_0$
for $X_1 \subseteq Z_1$ and $X_0 \subseteq Z_0$
contracting $D_1^i$.
\item
The weight $w_i := \deg(f_i) \in K$ is 
exceptional and  $w^0_i \in \lambda_0$
holds for a full dimensional 
$\lambda_0 \in \Lambda(\b{X})$ 
having a common facet with 
$\lambda_1 := \SAmple(X_1)$.
\end{enumerate}
If~(ii) holds, then any pair 
$\eta_0, \eta_1 \in \Lambda(\b{Z}_1)$ 
of adjacent fulldimensional chambers with 
with $\eta_i^\circ \subseteq \lambda_i^\circ$ 
defines a neat controlled ambient 
modification $Z_1 \to Z_0$ as in~(i).
Moreover, $X_0$ is a $\QQ$-factorial projective 
variety with finitely generated Cox ring.
\end{proposition}

\begin{proof}
Suppose that~(i) holds.
Then the ambient modification $Z_1 \to Z_0$
comes from a stellar subdivision 
$\Sigma_1 \to \Sigma_0$ at $v_i := P(e_i)$.
According to Lemma~\ref{subdiv2chambers},
the stellar subdivision $\Sigma_0 \to \Sigma_1$ 
corresponds to a pair $\eta_0, \eta_1 \in \Lambda(\b{Z}_1)$ 
of full dimensional chambers having a common facet and 
$w_i \in \eta_0$.
Then we have $\eta_i \subseteq \lambda_i$ with unique 
$\lambda_i \in \Lambda(\b{X}_1)$, and these
$\lambda_i$ are obviously as wanted.

Suppose that (ii) holds. Choose any full dimensional
$\eta_0 \in \Lambda(\b{Z}_1)$ such that 
$\eta_0 \subseteq \lambda_0$ holds and $\eta_0$ 
has a facet $\eta_{01}$ inside 
$\lambda_0 \cap \lambda_1$.
Since $\eta_0 \not\subseteq \cone(w_j; \; j \ne i)$
holds, we have $w_i \in \eta_0$.
Let $\eta_1 \in \Lambda(\b{Z}_1)$ be the fulldimensional 
cone with $\eta_{01} \preceq \eta_1$ and 
$\eta_1 \subseteq \lambda_1$.
Then Lemma~\ref{subdiv2chambers} tells us that
the pair $\eta_0,\eta_1$ defines a stellar
subdivision $\Sigma_1 \to \Sigma_0$ of fans
at $v_i = P(e_i)$.

Let $Z_i$ be the toric varieties associated to 
$\Sigma_i$, 
and let $\pi \colon Z_1 \to Z_0$ denote the 
toric contraction morphism associated to 
the stellar subdivision $\Sigma_1 \to \Sigma_0$.
Let $W_i \subseteq \b{Z}_i$ denote the sets 
of semistable points corresponding to the 
chambers $\eta_i$, and let 
$q_i \colon W_i \to Z_i$ be the quotient maps.
Then we have $W_1 = \rq{Z}_1$ and $q_1 = p_1$.
Moreover, in $Z_0$, we observe
$$ 
q_0(\b{X}_1 \cap W_0)
\ = \ 
\b{\pi(q_1(\b{X}_1 \cap \b{T}_1))}
\ = \ 
X_0,
$$
where $\b{T}_1 \subseteq \b{Z}_1$ denotes 
the big torus. In particular, we see
that $X_0$ is the good quotient space of an 
open subset of $\b{X}_1$ and thus is normal.
As it arises from a full dimensional
chamber, it is moreover $\QQ$-factorial.
Finally, Theorem~\ref{controlled2small2} 
guarantees that $\pi \colon Z_1 \to Z_0$ 
is a neat controlled ambient modification.
\end{proof}

Combining this Proposition with
the descriptions of the moving 
cone and the semiample cone given in 
Proposition~\ref{divcones}, we obtain 
the following characterizations 
of combinatorial minimality.

\begin{corollary}
A $\QQ$-factorial projective variety with 
finitely generated Cox ring is combinatorially 
minimal if and only if its effective cone and 
its moving cone coincide.
\end{corollary}

\begin{corollary}
A $\QQ$-factorial projective surface with 
finitely generated Cox ring is 
combinatorially minimal if and only if 
its effective cone and its semiample
cone coincide.
\end{corollary}

\begin{proof}[Proof of Theorem~\ref{combcontr}]
We may assume that $X$ arises from a bunched ring.
Then any contractible class $[D]$ on $X$ defines 
an exceptional weight.
By a suitable small birational transformation, 
we achieve that the semiample cone of $X$ 
has a common facet with a full dimensional 
chamber containing this exceptional weight.
Now we can use Proposition~\ref{combcontract}
to contract the class $[D]$ and obtain a 
$\QQ$-factorial projective 
variety with finitely generated Cox ring.
Repeating this procedure, we finally arrive 
at a combinatorially minimal variety.
\end{proof}

In order to underline the computational
nature of combinatorial contraction, we 
discuss an explicit example; it provides  
a non-toric surface that contracts to 
the projective plane.

\begin{example}
\label{ex:delpezzo2}
We consider again the $\QQ$-factorial 
projective surface of Example~\ref{ex:delpezzo}.
This time, we call it $X_1$; it arises 
from the bunched ring 
\begin{eqnarray*}
R 
& = & 
\KK[T_1, \ldots, T_5] \, / \, \bangle{T_1T_2 + T_3^2 + T_4T_5},
\end{eqnarray*}
with the system $\mathfrak{F}$ consisting of 
the classes 
$f_i \in R$ of $T_i \in \KK[T_1, \ldots, T_5]$.
The $K$-grading of $R$ is given 
by $\deg(f_i) := w_i$, with 
the columns $w_i$ of the matrix
\begin{eqnarray*}
Q_1
& := & 
\left[
\begin{array}{rrrrr}
1 & -1 & 0 & -1 & 1
\\
1 & 1 & 1 & 0 & 2
\end{array}
\right]
\end{eqnarray*}
and the $\mathfrak{F}$-bunch $\Phi$ 
consists of the ample cone 
$\tau = \cone(w_2,w_5)$. 
The weight $w_4 \in K$ is extremal, 
and thus there must be a contraction.

To determine the contraction explicitly,
we need an associated neat embedding 
into a toric variety as in Remark~\ref{assocemb}. 
The two GIT-fans look as follows;
the shaded area indicates~$Q^0(\gamma)$.
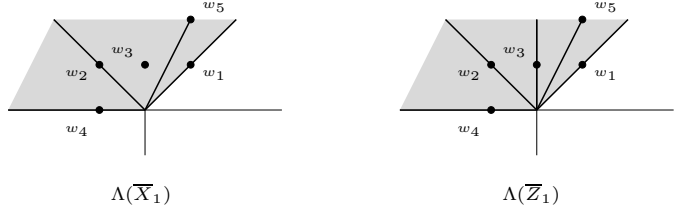
\begin{center}
\begin{picture}(0,0)%
\includegraphics{delptoricgitfan.pstex}%
\end{picture}%
\setlength{\unitlength}{1243sp}%
\begingroup\makeatletter\ifx\SetFigFont\undefined%
\gdef\SetFigFont#1#2#3#4#5{%
  \reset@font\fontsize{#1}{#2pt}%
  \fontfamily{#3}\fontseries{#4}\fontshape{#5}%
  \selectfont}%
\fi\endgroup%
\begin{picture}(5445,4159)(2218,-3338)
\put(6076,-736){\makebox(0,0)[lb]{\smash{{\SetFigFont{5}{6.0}{\familydefault}{\mddefault}{\updefault}{\color[rgb]{0,0,0}$w_1$}%
}}}}
\put(6076,614){\makebox(0,0)[lb]{\smash{{\SetFigFont{5}{6.0}{\familydefault}{\mddefault}{\updefault}{\color[rgb]{0,0,0}$w_5$}%
}}}}
\put(3376,-1861){\makebox(0,0)[lb]{\smash{{\SetFigFont{5}{6.0}{\familydefault}{\mddefault}{\updefault}{\color[rgb]{0,0,0}$w_4$}%
}}}}
\put(4276,-286){\makebox(0,0)[lb]{\smash{{\SetFigFont{5}{6.0}{\familydefault}{\mddefault}{\updefault}{\color[rgb]{0,0,0}$w_3$}%
}}}}
\put(3376,-736){\makebox(0,0)[lb]{\smash{{\SetFigFont{5}{6.0}{\familydefault}{\mddefault}{\updefault}{\color[rgb]{0,0,0}$w_2$}%
}}}}
\put(4276,-3211){\makebox(0,0)[lb]{\smash{{\SetFigFont{7}{8.4}{\familydefault}{\mddefault}{\updefault}{\color[rgb]{0,0,0}$\Lambda(\b{X}_1)$}%
}}}}
\end{picture}%

\qquad \qquad
\begin{picture}(0,0)%
\includegraphics{delpezzogitfan.pstex}%
\end{picture}%
\setlength{\unitlength}{1243sp}%
\begingroup\makeatletter\ifx\SetFigFont\undefined%
\gdef\SetFigFont#1#2#3#4#5{%
  \reset@font\fontsize{#1}{#2pt}%
  \fontfamily{#3}\fontseries{#4}\fontshape{#5}%
  \selectfont}%
\fi\endgroup%
\begin{picture}(5445,4159)(2218,-3338)
\put(6076,-736){\makebox(0,0)[lb]{\smash{{\SetFigFont{5}{6.0}{\familydefault}{\mddefault}{\updefault}{\color[rgb]{0,0,0}$w_1$}%
}}}}
\put(6076,614){\makebox(0,0)[lb]{\smash{{\SetFigFont{5}{6.0}{\familydefault}{\mddefault}{\updefault}{\color[rgb]{0,0,0}$w_5$}%
}}}}
\put(3376,-1861){\makebox(0,0)[lb]{\smash{{\SetFigFont{5}{6.0}{\familydefault}{\mddefault}{\updefault}{\color[rgb]{0,0,0}$w_4$}%
}}}}
\put(4276,-286){\makebox(0,0)[lb]{\smash{{\SetFigFont{5}{6.0}{\familydefault}{\mddefault}{\updefault}{\color[rgb]{0,0,0}$w_3$}%
}}}}
\put(3376,-736){\makebox(0,0)[lb]{\smash{{\SetFigFont{5}{6.0}{\familydefault}{\mddefault}{\updefault}{\color[rgb]{0,0,0}$w_2$}%
}}}}
\put(4276,-3211){\makebox(0,0)[lb]{\smash{{\SetFigFont{7}{8.4}{\familydefault}{\mddefault}{\updefault}{\color[rgb]{0,0,0}$\Lambda(\b{Z}_1)$}%
}}}}
\end{picture}%

\end{center}
We choose as toric ambient variety of $X_1$ 
the toric variety $Z_1$ defined by the 
GIT-cone generated by $w_2$ and $w_3$.
Consider
\begin{eqnarray*}
P_1
& := & 
\left[
\begin{array}{rrrrr}
1  & 0 & -1 & 1  & 0
\\
0  & 1 & -1 & -1 & 0
\\
-1 & 0 & -1 & 0  & 1
\end{array}
\right]
\end{eqnarray*}
Then the rays of the fan $\Sigma_1$ of $Z_1$ 
are the rays through 
the columns $v_1, \ldots, v_5$ and the maximal 
cones of $\Sigma_1$ are explicitly given as
$$ 
\cone(v_1,v_2,v_3), \qquad
\cone(v_1,v_2,v_5), \qquad
\cone(v_1,v_3,v_4), 
$$
$$
\cone(v_1,v_4,v_5), \qquad
\cone(v_2,v_3,v_5), \qquad
\cone(v_3,v_4,v_5).
$$ 
The toric ambient modification 
$\Sigma_1 \to \Sigma_0$ corresponding 
to the chambers $\cone(w_2,w_4)$ 
and  $\cone(w_3,w_2)$ of the 
fan $\Lambda(\b{Z}_1)$ contracts the ray $v_4$.
Thus, the primitive generators of 
$\Sigma_0$ are the vectors 
$v_1,v_2,v_3,v_5$ and its maximal cones 
are 
$$ 
\cone(v_1,v_2,v_3),
\qquad
\cone(v_1,v_2,v_5),
\qquad
\cone(v_2,v_3,v_5),
\qquad
\cone(v_1,v_3,v_5).
$$
Note that the last one contains the vector $v_4$;
in fact, we have $v_4  = 2v_1 + v_3 + 3v_5$.
The index of the stellar subdivision 
$\Sigma_1 \to \Sigma_0$ equals one, and 
the associated toric morphism 
$\b{\pi} \colon \b{Z}_1 \to \b{Z}_0$ 
of the toric total coordinate spaces is 
given by 
\begin{eqnarray*}
(z_1, \ldots, z_5)
& \mapsto & 
(z_1z_4^2,z_2,z_3z_4,z_4^3z_5).
\end{eqnarray*}
According to Theorem~\ref{controlled2small2}, 
the total coordinate space of $X_0$ is given as 
$\b{X}_0 = \b{\pi}(\b{X}_1)$. 
One easily checks that the vanishing  ideal 
$\b{X}_0$ in $\b{Z}_0 = \KK^4$ is 
generated by the polynomial
$$ 
T_1T_2 + T_3^2 + T_4.
$$
In particular, we see that the total coordinate space
$\b{X}_0$ is the affine 3-space; 
it can explicitly be parametrized
by 
$$
\KK^3 \ \to \ \b{X}_0,
\qquad
(w_1,w_2,w_3) \ \to \ (w_1,w_2,w_3,-w_1w_2-w_3^2).
$$
Using this, we obtain that the 
induced action of $H_0 \cong \KK^*$ 
is scalar multiplication.
Thus, $X_0$ is a projective plane.
Note that the map $X_1 \to X_0$ contracts 
the smooth rational curve 
$p_1(V(\rq{X}_0;T_4)) \subseteq X_1$ 
containing the singular point of $X_1$.
\end{example}

\section{Example: Cox rings with a single relation}
\label{sec:singrel}

Here we apply the results on toric ambient 
modifications to Cox rings defined by a single relation.
We give a simple criterion for neat controlled
ambient modifications in this setting,
and we explicitly determine the resulting Cox
rings; this can be used to construct new factorial 
hypersurfaces out of given ones. 
As a concrete example, we perform the minimal
resolution of a singular del Pezzo surface 
via toric ambient resolution.

The setup of this section is the following.
Let $K_0$ be a finitely generated abelian group,
and endow the polynomial ring 
$\KK[T_1, \ldots ,T_r]$ 
with a $K_0$-grading by setting 
$$ 
\deg(T_i) \ := \ w_i, \qquad
\text{where }
w_1, \ldots, w_r \in K_0.
$$
Set $\b{Z}_0 := \KK^r$ and $H_0 := \Spec(\KK[K_0])$,
and suppose that 
$\rq{Z}_0 \to Z_0 = \rq{Z}_0 \quot H_0$ 
are the data of a toric Cox 
construction.
Then $\rq{Z}_0 \to Z_0$ is given by a map 
$\rq{\Sigma}_0 \to \Sigma_0$ of fans
living in lattices $F_0 = \ZZ^r$ and $N_0$.
We assume that  $\rq{\Sigma}_0$ is a 
maximal projectable fan; in particular, 
$\Sigma_0$ cannot be enlarged without adding 
new rays.

Let $f_0 \in \KK[T_1, \ldots, T_r]$
be homogeneous with respect to the 
$K_0$-grading and suppose that it defines 
a factorially $K_0$-graded residue algebra
\begin{eqnarray*}
R 
& := &
\KK[T_1, \ldots, T_r] / \bangle{f_0}
\end{eqnarray*}
such that the classes $T_i + \bangle{f_0}$ 
are pairwise nonassociated $K_0$-prime 
elements in $R$. Consider 
$$
\b{X}_0 \ := \ V(f_0), 
\qquad
\rq{X}_0 \ := \ \b{X}_0 \cap \rq{Z}_0,
\qquad
X_0 \ := \ \rq{X}_0 \quot H \ \subseteq\  Z_0.
$$
Then $X_0$ is a normal variety having 
the $H_0$-variety $\b{X}_0$ as its 
total coordinate space, and $X_0 \to Z_0$ 
is a neat embedding, see Proposition~\ref{neatembex}.

Let $v_1, \ldots, v_r$ be the primitive
lattice vectors in the rays of $\Sigma_0$.
We suppose that for $2 \le d \le r$,
the cone $\sigma_0$ generated by 
$v_1, \ldots, v_d$ belongs to $\Sigma_0$
and consider the stellar subdivision 
$\Sigma_1 \to \Sigma_0$ 
at a vector
\begin{eqnarray*}
v_\infty 
& = &
a_1v_1 + \cdots + a_dv_d.
\end{eqnarray*}
Denote the index of this subdivision by $m_\infty$
and the associated toric modification  by
$\pi \colon Z_1 \to Z_0$.
Then we obtain the strict transform 
$X_1 \subseteq Z_1$ mapping onto 
$X_0 \subseteq Z_0$.
Moreover, as in Section~\ref{sec:ambmod},
we have the commutative diagrams
$$ 
\xymatrix{
&
{\b{Z}_1}  \ar[dr]^{\b{\pi}} \ar[dl]_{\b{\pi}_1}
& 
\\
{\b{Z}_1}
&
{\rq{Z}_1} \ar[u] \ar[dl]_{\rq{\pi}_1} \ar[dr]^{\rq{\pi}}
& 
{\b{Z}_0}
\\
{\rq{Z}_1} \ar[d]_{p_1}^{/ H_1} \ar[u]
&
& 
{\rq{Z}_0} \ar[u] \ar[d]^{p_0}_{/ H_0}
\\
Z_1 
\ar[rr]_{\pi}
& &
Z_0
}
\qquad \qquad
\xymatrix{
&
{\b{Y}_1}  \ar[dr]^{\b{\kappa}} \ar[dl]_{\b{\kappa}_1}
& 
\\
{\b{X}_1}
&
{\rq{Y}_1} \ar[u] \ar[dl]_{\rq{\kappa}_1} \ar[dr]^{\rq{\kappa}}
& 
{\b{X}_0}
\\
{\rq{X}_1} \ar[d]_{p_1}^{/ H_1} \ar[u]
&
& 
{\rq{X}_0} \ar[u] \ar[d]^{p_0}_{/ H_0}
\\
X_1 
\ar[rr]_{\kappa}
& &
X_0
}
$$
where $H_1 = \Spec(\KK[K_1])$ for 
$K_1 = E_1 / M_1$ in analogy to 
$K_0 = E_0 / M_0$;
see Section~\ref{sec:crembvar}
for the notation.
In particular, the Cox ring 
$\KK[T_1, \ldots, T_r,T_\infty]$
of $\b{Z}_1$ comes with a $K_1$-grading.
Recall moreover, that with respect 
to the coordinates corresponding 
to the rays of the fans $\Sigma_i$, 
the map 
$\b{\pi} \colon \b{Z}_1 \to \b{Z}_0$ is given by 
$$ 
\b{\pi}(z_1, \ldots, z_r,z_\infty) 
\ = \
(z_\infty^{a_1}z_1, 
\ldots, z_\infty^{a_d}z_d, 
z_{d+1}, \ldots, z_r).
$$

We want to give an explicit criterion 
for $Z_1 \to Z_0$ to be a neat controlled ambient 
modification for $X_0 \subset Z_0$,
and, in this case, describe the total coordinate
space of $X_1$ explictily. 
For this, consider the grading 
$$ 
\KK[T_1, \ldots, T_r] 
\ = \ 
\bigoplus_{k \ge 0} \KK[T_1, \ldots, T_r]_k,
\quad
\text{where }
\deg(T_i) 
\ := \ 
\begin{cases}  
a_i & i \le d,
\\
0   & i \ge d+1.
\end{cases}
$$
Then we may write $f_0 = g_{k_0} + \cdots + g_{k_m}$
where $k_0 < \cdots < k_m$ and each $g_{k_i}$
is a nontrivial polynomial  having 
degree $k_i$ with respect to this grading.

\begin{definition}
We say that the polynomial $f_0 \in \KK[T_1, \ldots, T_r]$ 
is {\em admissible\/} if 
\begin{enumerate}
\item
the toric orbit $0 \times \TT^{r-d}$ intersects 
$\b{X}_0 = V(f_0)$,
\item
$g_{k_0}$ 
is a $K_1$-prime polynomial in at 
least two variables.
\end{enumerate}
\end{definition}

Note that for the case of a free abelian group $K_1$,
the second condition just means that $g_{k_0}$ is 
an irreducible polynomial.

\begin{proposition}\label{singlerel}
If, in the above setting, the 
polynomial $f_0$ is admissible,
then the following holds.
\begin{enumerate}
\item
The toric morphism $Z_1 \to Z_0$ is a neat controlled  
toric ambient modification for $X_0 \subseteq Z_0$
and $X_1 \subseteq Z_1$.
\item
The strict transform $X_1 \subseteq Z_1$ is a neatly 
embedded normal subvariety with total coordinate 
space $\b{X}_1$ and Cox ring
$$
\mathcal{R}(X_1) 
\ = \ 
\KK[T_1, \ldots, T_r,T_\infty] / 
\bangle{f_1(T_1, \ldots, T_r,\sqrt[m_\infty]{T_\infty})},
$$
where in
$$ 
f_1 
\ := \ 
\frac{f_0(T_\infty^{a_1}T_1, 
\ldots, T_\infty^{a_d}T_d, 
T_{d+1}, \ldots, T_r )}{T_\infty^{k_0}}
\ \in \ 
\KK[T_1, \ldots, T_r,T_\infty]
$$
only powers $T_\infty^{lm_\infty}$ with $l \ge 0$
of $T_\infty$ occur, 
and the notation $\sqrt[m_\infty]{T_\infty}$ means 
replacing 
 $T_\infty^{lm_\infty}$ with $T_\infty^l$ in $f_1$. 
\end{enumerate}
\end{proposition}

\begin{lemma}
\label{singinf}
Suppose that the polynomial ring $\KK[T_1, \ldots, T_r,T_\infty]$ 
is graded by some finitely generated abelian group $K_1$.
Let $f_1 = g T_\infty + h \in \KK[T_1, \ldots, T_r,T_\infty]$
be irreducible with 
$g \in \KK[T_1, \ldots, T_r,T_\infty]$
and a $K_1$-prime $h \in \KK[T_1, \ldots, T_r]$.
Then, for $\b{X}_1 := V(f_1)$
and 
$H_1 := \Spec(\KK[K_1])$,
the intersection $\b{X}_1 \cap V(T_\infty)$
is $H_1$-irreducible, and 
$\Sing(\b{X}_1) \cap V(T_\infty)$ 
is a proper subset of 
$\b{X}_1 \cap V(T_\infty)$.
\end{lemma}

\begin{proof}
We only have to show that  
$A := V(T_\infty, f_1, \partial f_1 / \partial T_1 , 
\ldots, 
\partial f_1 / \partial T_\infty)$ 
is a proper subset of 
$B := V(f_1, T_\infty)$.
A simple calculation gives
\begin{eqnarray*}
\grad (f_1)
& = & 
\left(
\frac{\partial g}{T_1} T_\infty
+ 
\frac{\partial h}{T_1},
\ldots,
\frac{\partial g}{T_r} T_\infty
+ 
\frac{\partial h}{T_r},
\frac{\partial g}{T_\infty} T_\infty + g
\right).
\end{eqnarray*}
Thus, setting $T_\infty = 0$, wee see that 
$A = B$ implies vanishing of $\grad(h)$ along
$V(h)$. This contradicts to the assumption
that $h$ is $K_1$-prime.
\end{proof}

\begin{proof}[Proof of Proposition~\ref{singlerel}]
We first show that $\b{X}_1 \subseteq \b{Z}_1$
is normal.
Recall that  
$\b{\kappa} \colon \b{Y}_1 \to \b{X}_0$
is a good quotient for the $\KK^*$-action on 
$\b{Y}_1$. 
Outside $V(T_\infty)$, this action is free
and  $\b{Y}_1$
locally splits as $\b{X}_0 \times \KK^*$.
In particular, $\b{Y}_1 \setminus V(T_\infty)$
inherits normality from $\b{X}_0$.
Let $\b{T}_0 \subseteq \b{X}_0$
and $\b{T}_1 \subseteq \b{X}_1$ 
be the big tori. Then we have
$$
\b{Y}_1 
\ = \ 
\b{\b{\kappa}^{-1}(\b{X}_0 \cap \b{T}_0)}
\ = \
\b{V(\b{T}_1;\b{\kappa}^*f_0)}
\ = \ 
V(\b{Z}_1; f_1).
$$
Combining this with Lemma~\ref{singinf},
we see that $\b{Y}_1$ is regular in codimension
one. Thus, Serre's criterion shows
that $\b{Y}_1$ is normal.
Moreover, the group $C_{m_\infty}$ of roots of unity 
acts on $\b{Z}_1$ via multiplication on the last 
coordinate, and $\b{Y}_1$ is invariant under this action.
Thus, $f_1$ is $C_{m_\infty}$-homogeneous.
Writing
\begin{eqnarray*}
f_1
& = & 
\sum_{i=0}^m
\frac{g_{k_i}(T_\infty^{a_1}T_1, 
\ldots, T_\infty^{a_d}T_d, 
T_{d+1}, \ldots T_r )}{T_\infty^{k_0}},
\end{eqnarray*}
we see that the $g_{k_0}$-term is invariant 
under $C_{m_\infty}$, and hence the others
must be as well; note that in different $g_{k_i}$,
the variable $T_\infty$ occurs in different
powers.
Thus, only powers $T_\infty^{lm_\infty}$
of $T_\infty$ occur in $f_1$.
Now, $\b{X}_1$ is the quotient of $\b{Y}_1$ by the action 
of $C_{m_\infty}$.
Consequently, $\b{X}_1$ is normal, and we 
obtain
\begin{eqnarray*}
\b{X}_1
&  = & 
V(\b{Z}_1; f_1(T_1, \ldots, T_r,\sqrt[m_\infty]{T_\infty})).
\end{eqnarray*}
Now the assertions drop out.
First, as a quotient of the normal 
variety $\rq{X}_1 \subseteq \b{X}_1$, 
the variety $X_1$ is again normal.
Since $\b{X}_1 \cap V(T_\infty)$ is 
$H_1$-irreducible, the exceptional locus
of $X_1 \to X_0$, as the image of 
$\rq{X}_1 \cap V(T_\infty)$ under 
the quotient $\rq{X}_1 \to X_1$ by $H_1$, 
is irreducible.
Moreover, note that
$$ 
\Gamma(\b{X}_1,\mathcal{O}) / \bangle{T_{\infty \vert \b{X}_1}}
\ \cong \
\KK[T_1, \ldots, T_r,T_\infty] / \bangle{f_1,T_\infty}
\ \cong \ 
\KK[T_1, \ldots, T_r] / \bangle{g_{k_0}}.
$$
We can conclude that the 
restriction of $T_\infty$ to $\b{X}_1$
provides a local equation for 
$\b{X}_1 \cap V(T_\infty)$.
Thus, we verified all conditions
for $Z_1 \to Z_0$ to be a neat 
ambient modification.
Consequently, Corollary~\ref{cor:contr2small}
shows that $\b{X}_1$ is the total
coordinate space of $X_1$.
\end{proof}

Proposition~\ref{singlerel} may as well 
be used to produce multigraded factorial
rings, as the following statement shows.

\begin{corollary}
The ring $R := \KK[T_1, \ldots, T_r,T_\infty] / 
\bangle{f_1(T_1, \ldots, T_r,\sqrt[m_\infty]{T_\infty})}$
in Proposition~\ref{singlerel} is factorially 
$K_1$-graded.
Moreover, if $K_0$ is torsion free, 
then $R$ is even factorial.
\end{corollary}

As an example, we consider a
$\QQ$-factorial projective surface $X_0$ 
with a single singularity,
and resolve this singularity by means of 
toric ambient resolution,
showing thereby that $X_0$ 
is a singular del Pezzo surface.

\begin{example}
\label{ex:delpezzo3}
We consider once more the $\QQ$-factorial 
projective surface arising 
from the bunched ring $(R,\mathfrak{F},\Phi)$ 
of~\ref{ex:delpezzo}.
This time, we call it $X_0$.
Recall that we have 
\begin{eqnarray*}
R 
& = & 
\KK[T_1, \ldots, T_5] \, / \, \bangle{T_1T_2 + T_3^2 + T_4T_5},
\end{eqnarray*}
the system $\mathfrak{F}$ consists of the classes 
$\b{T}_i \in R$ of $T_i \in \KK[T_1, \ldots, T_5]$.
Moreover, the $K_0$-grading of $R$ is given 
by $\deg(\b{T}_i) := w_i$, with 
the columns $w_i$ of the matrix
\begin{eqnarray*}
Q_0
& := & 
\left[
\begin{array}{rrrrr}
1 & -1 & 0 & -1 & 1
\\
1 & 1 & 1 & 0 & 2
\end{array}
\right]
\end{eqnarray*}
and the $\mathfrak{F}$-bunch $\Phi$ 
consists of the ample cone 
$\tau = \cone(w_2,w_5)$. 
As in~\ref{ex:delpezzo2},
we take the toric variety $Z_0$ defined by the 
GIT-cone generated by $w_2$ and $w_3$
as the toric ambient variety of $X_0$.
Then the columns  $v_1, \ldots, v_5$ of 
\begin{eqnarray*}
P_0
& := & 
\left[
\begin{array}{rrrrr}
1  & 0 & -1 & 1  & 0
\\
0  & 1 & -1 & -1 & 0
\\
-1 & 0 & -1 & 0  & 1
\end{array}
\right]
\end{eqnarray*}
generate  the rays of $\Sigma_0$, and 
the maximal cones of $\Sigma_0$ 
are explicitly given as
$$ 
\cone(v_1,v_2,v_3), \qquad
\cone(v_1,v_2,v_5), \qquad
\cone(v_1,v_3,v_4), 
$$
$$
\cone(v_1,v_4,v_5), \qquad
\cone(v_2,v_3,v_5), \qquad
\cone(v_3,v_4,v_5).
$$ 
Note that $X_0$ inherits its singularity 
from $Z_0$; it is the toric singularity 
corresponding to the cone generated by 
$v_1,v_3$ and $v_4$. 
In view of Corollary~\ref{singulemb},
we may obtain a resolution of the singularity 
of $X$ by resolving the ambient singularity.

Resolving the toric ambient singularity
coresponding to $\cone(v_1,v_3,v_4)$
means to successively subdivide this cone
along the interior members of the Hilbert basis.
These are $v_6 := (0,-1,-1)$ and $v_7 := (1,-1,-1)$;
note that the order does not influence the result
on the resolution obtained for $X_0$.

We start with subdividing in $v_6 = (0,-1,-1)$.
Up to renumbering coordinates, we are in the setting 
of Proposition~\ref{singlerel}.
Note that the index of this subdivision is $m_6 = 3$.
The pullback relation on $\b{Z}_1$ is 
$T_1T_2 + T_3^2T_6^3 + T_4T_5$.
Dividing by the $C_3$-action, we obtain
$$ 
f_1
\ := \ 
T_1T_2 + T_3^2T_6 + T_4T_5
\ \in \ 
\KK[T_1, \ldots, T_6].
$$
for the relation in the Cox ring of 
the modified surface $X_1 \subseteq Z_1$. 
This still has a singular point, namely 
the toric fixed point corresponding to the 
cone generated by $v_1,v_4$ and $v_6$.

In the next step, we have to subdivide by $v_7$.
Note that now we have index $m_7=2$.
The pullback relation on $\b{Z}_2$ is 
$T_1T_2 + T_3^2T_6^3 + T_4T_5$.
It doesn't depend on the new variable $T_7$,
hence, dividing by the $C_2$-action, we get 
$$ 
f_2
\ := \ 
T_1T_2 + T_3^2T_6 + T_4T_5
\ \in \ 
\KK[T_1, \ldots, T_7]
$$
for the defining relation of the Cox ring 
of the modified surface $X_2 \subseteq Z_2$.
In order to see that $X_2 \to X_0$ is a 
resolution, we need to know that $\rq{X}_2$
is smooth.
But this is due to the fact that the singular 
locus of $\b{X}_2$ is described by $T_6 = 0$ 
and $w_6$ is exceptional.
Note that the grading of the Cox ring of $X_2$ 
is given by assigning 
to $T_i$ the $i$-th column of the matrix
\begin{eqnarray*}
Q_2
& := & 
\left[
\begin{array}{rrrrrrr}
1 & -1 & 0 & -1 & 1 & 0 & 0 
\\
0 & 1 & 0 & 0 & 1 & 1 & 0
\\
1 & 0  & 1 & 0 & 1 & -1 & 0
\\
0 & 0  & 0 & -1 & 1 & 0 & 1
\end{array}
\right]
\end{eqnarray*}
\end{example}

\end{document}